\documentclass[final,leqno]{siamltex}

\usepackage{amsmath,amsfonts}
\usepackage{color,graphicx}
\usepackage{hyperref}
\usepackage{enumerate}

\newcommand{\Dt}{\Delta {t}}
\newcommand{\R}{\mathbb{R}}
\newcommand{\Ad}{A}
\newcommand{\Wz}{z}
\newcommand{\Adcut}{S}

\newcommand{\Crand}{C_z}
\newcommand{\Cchar}{C_{s}}
\newcommand{\Cr}{C^{\textsc{f}}_{r}}
\newcommand{\Crl}{C^{\textsc{l}}_{r}}
\newcommand{\Ct}{C_{\tau}}
\newcommand{\Ca}{C_{a}}
\newcommand{\HFF}{H^{\textsc{f}}}
\newcommand{\HFL}{H^{\textsc{l}}}
\newcommand{\HLF}{K^{\textsc{f}}}
\newcommand{\HLL}{K^{\textsc{l}}}
\newcommand{\NF}{N^{\textsc{f}}}
\newcommand{\NL}{N^{\textsc{l}}}
\newcommand{\Nmpc}{N_{\textsc{mpc}}}
\newcommand{\rhoF}{m^{\textsc{f}}}
\newcommand{\rhoL}{m^{\textsc{l}}}
\newcommand{\etaF}{\eta^{\textsc{f}}}
\newcommand{\etaL}{\eta^{\textsc{l}}}
\newcommand{\lamF}{\lambda^{\textsc{f}}}
\newcommand{\lamL}{\lambda^{\textsc{l}}}
\newcommand{\thetaF}{\vartheta^{\textsc{f}}}

\newcommand{\RestF}{\Xi_{\varphi}^{\textsc{f}}}
\newcommand{\RestL}{\Xi_{\varphi}^{\textsc{l}}}
\newcommand{\barS}{s}
\newcommand{\noisepar}{\Crand}
\newcommand{\xf}{x}
\newcommand{\xl}{y}
\newcommand{\vf}{v}
\newcommand{\vl}{w}
\newcommand{\target}{x^\tau}
\newcommand{\CclosenessLF}{\mu^{\textsc{l}}}
\newcommand{\CclosenessFF}{\mu^{\textsc{f}}}

\newcommand{\Ccontrol}{\nu}

\graphicspath{{./figs/},{./figsR1/}}

\newtheorem{remark}[theorem]{Remark}

\title{Invisible control of self-organizing agents \\ leaving unknown environments\thanks{This research has received funding from the European Union under grant ERC-STG 306274 (HDSPCONTR) and ERC-ADG 668998 (OCLOC).}}

\author{
Giacomo Albi\thanks{TUM, 
M15 - Chair of Numerical Mathematics, Boltzmann stra\ss e. 3, D-85748 Garching, Munich, Germany. \texttt{albi@ma.tum.de}}   
\and	
Mattia Bongini\thanks{TUM, 
M15 - Chair of Numerical Mathematics, Boltzmann stra\ss e. 3, D-85748 Garching, Munich, Germany. \texttt{bongini@ma.tum.de}} 
\and	
Emiliano Cristiani\thanks{Istituto per le Applicazioni del Calcolo ``M. Picone'', Consiglio Nazionale delle Ricerche, Via dei Taurini 19, I-00185, Rome, Italy.
\texttt{e.cristiani@iac.cnr.it}}
\and	
Dante Kalise\thanks{Johann Radon Institute for Computational and Applied Mathematics, Altenberger Stra\ss e 69, A-4040 Linz, Austria.
\texttt{dante.kalise@ricam.oeaw.ac.at}}	}
		
\begin{document}
\maketitle

\begin{abstract}
In this paper we are concerned with multiscale modeling, control, and simulation of self-organizing agents leaving an unknown area under limited visibility, with special emphasis on crowds.
We first introduce a new microscopic model characterized by an \textit{exploration phase} and an \textit{evacuation phase}. The main ingredients of the model are an alignment term, accounting for the \emph{herding effect} typical of uncertain behavior, and a random walk, accounting for the need to explore the environment under limited visibility. We consider both metrical and topological interactions. Moreover, a few special agents, the \emph{leaders}, not recognized as such by the crowd, are ``hidden'' in the crowd with a special controlled dynamics.
Next, relying on a Boltzmann approach, we derive a mesoscopic model for a continuum density of followers, coupled with a microscopic description for the leaders' dynamics. 
Finally, optimal control of the crowd is studied. It is assumed that leaders aim at steering the crowd towards the exits so to ease the evacuation and limit clogging effects, and locally-optimal behavior of leaders is computed. 
Numerical simulations show the efficiency of the control techniques in both microscopic and mesoscopic settings. We also perform a real experiment with people to study the feasibility of such a bottom-up control technique.
\end{abstract}



\pagestyle{myheadings}
\thispagestyle{plain}
\markboth{G. ALBI, M. BONGINI, E. CRISTIANI, D. KALISE}{INVISIBLE CONTROL OF SELF-ORGANIZING AGENTS}

\section{Introduction}\label{sec:introduction}
Self-organizing systems are not as mysterious as they were when researchers began to deal with them. A number of local interaction rules were identified, such as repulsion, attraction, alignment, self-propulsion, etc., and the large-scale outcome triggered by the combination of them is so far sufficiently understood \cite{vicsek2012PR}. Researchers have also achieved a certain knowledge regarding the \emph{control} of self-organizing systems, namely the methods to dictate a target behavior to the agents still keeping the local rules active all the time. 

In this paper we are concerned with multiscale modeling, control and simulations of self-organizing agents which leave an unknown area. Agents are not informed about the positions of the exits, so that they need to explore the environment first. Agents are assumed not to communicate with each other and they cannot share information directly. Moreover, they are selfish and do not intend to help the other agents.

The guiding application is that of a human crowd leaving an unknown environment under limited visibility. In that case, we claim that people exhibit two opposite tendencies: on the one hand, they tend to spread out in order to explore the environment efficiently, on the other hand they follow the group mates in the hope that they have already found a way out (\emph{herding effect}). Controlling these natural behaviors is not easy, especially in emergency or panic situations. 
Typically, the emergency management actors (e.g., police, stewards) adopt a \emph{top-down} approach: they receive information about the current situation, update the evacuation strategies, and finally inform people communicating directives. It is useful to stress here that this approach can be very inefficient in large environments (where communications to people are difficult) and in panic situations (since ``instinctive'' behavior prevails over the rational one). Let us also mention that in some cases people are simply not prone to follow authority's directives because they are opposing it (e.g., in public demonstrations).

In this paper we explore the possibility of controlling crowds adopting a \emph{bottom-up} approach. Control is obtained by means of special agents (leaders) who are \emph{hidden in the crowd, and they are not recognized by the mass as special}. From the modeling point of view, this is translated by the fact that the individuals of the crowd (followers) interact in the same way with both the other crowd mates and the leaders. 

\paragraph{Goals} 
We first introduce a new microscopic (agent-based) model characterized by an \textit{exploration phase} and an \textit{evacuation phase}. This is achieved by a combination of repulsion, alignment, self-propulsion and random (social) force, together with the introduction of an \textit{exit's visibility area}. We also consider both \emph{metrical} and \emph{topological} interactions so to avoid unnatural all-to-all interactions. A few leaders, not recognized as such by the crowd, are added to the model with a special, controlled dynamics.

Second, relying on a Boltzmann approach, we derive a mesoscopic model for a continuum density of followers, coupled with a microscopic description for the leaders' dynamics. 
This procedure, based on a binary interaction approximation, shows that, in a grazing collision regime, the microscopic-kinetic limit system is the mesoscopic description of the original microscopic dynamics.

Finally, dynamic optimization procedures with short and long-time horizon are considered. It is assumed that leaders aim at steering the crowd toward the exits so to ease the evacuation and limit clogging effects. The goal is achieved by exploiting the herding effect. Locally-optimal behavior of leaders is computed by means of two methods: the model predictive control and a modified compass search. 
Several numerical simulations show the efficiency of the optimization methods in both the microscopic and mesoscopic settings.

Beside numerical simulations, we performed a real experiment with people to study the feasibility of the proposed bottom-up crowd control technique.

\paragraph{Relevant literature} 
This paper falls in several crossing, and often independent, lines of research. 
%
Concerning pedestrian modeling, virtually any kind of models have been investigated so far and several reviews and books are available. For a quick introduction, we refer the reader to the reviews \cite{bellomo2011SR, helbing2001RMP} and the books \cite{cristiani2014book, kachroo2008book}. Some papers deal specifically with \emph{evacuation problems}: a very good source of references is the paper \cite{abdelghany2014EJOR}, where evacuation models both with and without optimal planning search are discussed. The paper \cite{abdelghany2014EJOR} itself proposes a cellular automata model coupled with a genetic algorithm to find a top-down optimal evacuation plan.  Evacuation problems were studied by means of lattice models \cite{cirillo2013PhysA, guo2012TRB}, social force models \cite{parisi2005PhysA}, cellular automata models \cite{abdelghany2014EJOR, wang2015PhysA}, mesoscopic models \cite{agnelli2015M3AS}, and macroscopic models \cite{carrillo1501.07054}. Limited visibility issues were considered in \cite{carrillo1501.07054, cirillo2013PhysA, guo2012TRB}. A real experiment involving people can be found in \cite{guo2012TRB}.

The control of self-organizing agents can be achieved in several ways. For example, the dynamics of all agents can be controlled by means of a single control variable, which remains active at all times \cite{albi1401.7798}. This kind of control is suitable, e.g., to model the influence of television on people. Alternatively, every agent can be controlled by an independent control variable. This is the most effective but also the most ``expensive'' control technique because of the large number of interactions with the agents that is required \cite{bongini2015conditional}. A more parsimonious control technique, called \emph{sparse control}, can be obtained by penalizing the number of these interventions by means of an $L^1$ control cost, so that the control is active only on few agents at every instant \cite{bofo13, bonginijunge2014sparse, caponigro2013MCRF}. If the existing behavioral rules of the agents cannot be redesigned, the control of the system can be obtained by means of external agents. Here the literature splits in two branches: the case of \emph{recognizable} leaders, characterized by the fact that the external agents have an influence on normal agents stronger than the influence exerted by group mates (like, e.g., a celebrity) \cite{albi1405.0736, borzi2015M3ASa, during2009PRSA, motsch2011JSP}; and the case of \emph{invisible} leaders, which are instead characterized by the fact that the external agents are completely anonymous and are perceived as normal group mates \cite{couzin2005N, duan2014SR, han2006JSSC, han2013PLOSONE}. Alternative approaches, based on the control of the surrounding environment rather than of the agents themselves are proposed in \cite{cristiani_obstacles, cristiani2015SIAP, johansson2007}.

Let us focus on the invisible control, also known as \emph{soft control}. In the seminal paper \cite{couzin2005N} a repulsion-alignment-attraction model is considered, for which the authors show that a little percentage of informed agents pointing towards a target is sufficient to steer the whole system to it. 
The papers \cite{han2006JSSC, han2013PLOSONE} deal with a Vicsek-like (pure alignment) model with a single leader. A feedback control for the leader is proposed, which is able to align the whole group along a desired direction from any initial condition. Invisible control was successfully implemented in real experiments involving animals, where leaders come in the form of disguised robots; see, e.g., \cite{halloy2007S}, which deals with cockroaches and \cite{butail2013PLOSONE}, which deals with zebrafishes.

The complexity reduction of microscopic models when the number of interacting agents is large has received deep interest in the last years. The main technique to address the \emph{curse of dimensionality} is to recast the original microscopic dynamics in the form of a PDE, by substituting the influence that the entire population has on a single agent with an averaged one. Several approaches have been used to derive rigorously this procedure, like the BBGKY hierarchy \cite{Tadmor2008}, the mean-field limit \cite{bolley11MMMAS, CCR11}, or the binary interaction approximation \cite{albi2013MMS, PT:13}. Furthermore, only recently this problem has been investigated for optimal control of multi-agent systems \cite{fornasier1306.5913}.

In the case of the control by means of external agents, if the number of leaders is relatively small in comparison to the crowd of followers, it might be convenient to apply the above dimension reduction techniques to the followers' population only \cite{FornasierSolombRossiBongini,fornasier2014PTRSA}. The resulting modeling setting is described by a coupled ODE-PDE system, where we control the microscopic dynamics in order to steer the macroscopic quantities. Nonetheless, from the computational point of view, one has to face the numerical tractability of the optimal control problem, whose solution requires the iterated evaluation of the macroscopic system. Since the underling dynamics are typically non-linear and imply the resolution of high-dimensional integro-differential terms, \emph{ad hoc} approaches have to be developed in order to avoid time-consuming strategies. One possibility is to use a fixed dynamics for the system of ODEs, assigned by the modelers or as feedback control derived by a properly designed functional, see, e.g., \cite{albi2013AML, colombo2012JNS}.

\paragraph{Paper organization} 
The paper is organized as follows. 
In Section \ref{sec:modelguidelines} we introduce the main features of the model, valid at any scale of observation. 
In Section \ref{sec:micromodel} we present in detail the model at the microscopic scale. 
In Section \ref{sec:mesomodel} we derive the associated Boltzmann-type equation and the coupled ODE (micro) - PDE (kinetic) model by means of the \textit{grazing collision limit}, following the path laid out in \cite{toscani2006kinetic}. 
In Section \ref{sec:optimalcontrol} we introduce the control and optimization problem and 
in Section \ref{sec:simulations} we show the results of the numerical tests. 
In Section \ref{sec:theexperiment} we report the results of a real experiment with pedestrians aimed at validating the proposed crowd control technique. 
We conclude the paper sketching some research perspectives.

\section{Model guidelines}\label{sec:modelguidelines}
In this section we describe the model at a general level. 
Hereafter, we divide the population between \textit{leaders}, which are the controllers and behave in some optimal way (to be defined), and \textit{followers}, which represent the mass of agents to be controlled. Followers cannot distinguish between followers and leaders. 

\paragraph{First- vs.\ second-order model} One can notice that walking people, animals and robots are in general able to adjust their velocity almost instantaneously, reducing to a negligible duration the acceleration/deceleration phase. For this reason, a second-order (inertia-based) framework does not appear to be the most natural setting. However, if we include in the model the tendency of the agents to move together and to align with group mates we implicitly assume that agents can perceive the velocity of the others. This makes unavoidable the use of a second-order model, where both positions and velocities are state variables. See \cite[Sect.\ 4.2]{cristiani2014book} for a general discussion on this point. In our model we adopt a mixed approach, describing leaders by a first-order model (since they do not need to align) and followers by a second-order one. In the latter case, the small inertia is obtained by means of a fast relaxation towards the target velocity. 

\paragraph{Metrical vs.\ topological interactions} Interaction of one agent with group mates is said to be \textit{metrical} if it involves only mates within a predefined \textit{sensory region}, regardless of the number of individuals which actually fall in it. Interaction is instead said \textit{topological} if it involves a predefined number of group mates regardless their distance from the considered agent. 
Again, in our model we adopt a mixed approach, assuming short-range interactions to be metrical and long-range ones to be topological.

\paragraph{Isotropic vs.\ anisotropic interactions} Living agents (people, animals) are generally asymmetric, in the sense that they better perceive stimuli coming from their ``front'', rather than their ``back''. In pedestrian models, e.g., interactions are usually restricted to the half-space in front of the person, since the human visual field is approximately 180$^\circ$. Nevertheless, humans can easily see all around them simply turning the head. In the case we are interested in, people have no idea of the location of their target, hence we expect that they often look around to explore the environment and see the behavior of the others. This is why we prefer to adhere to isotropic interactions.

Let us describe the social forces acting on the agents.

\emph{Leaders.}
\begin{itemize}
\item Leaders are subject to an isotropic metrical short-range \textit{repulsion force} directed against all the others, translating the fact that they want to avoid collisions and that a maximal density exists.
\item Leaders are assumed to know the environment and the self-organizing features of the crowd. They respond to an \textit{optimal force} which is the result of an offline optimization procedure, defined as to minimizing some cost functional (see Section \ref{sec:optimalcontrol}).
\end{itemize}

\emph{Followers.}
\begin{itemize}
\item Similarly to leaders, followers respond to an isotropic metrical short-range \textit{repulsion force} directed against all the others.
\item Followers tend to a \textit{desired velocity} which corresponds to the velocity the would follow if they were alone in the domain. This term takes into account the fact that the environment is unknown. Followers describe a \textit{random walk} if the exit is not visible (exploration phase) or a sharp motion toward the exit if the exit is visible (evacuation phase). In addition, we include a self-propulsion term which translates the tendency to reach a given characteristic speed (modulus of the velocity).
\item If the exit is not visible, followers are subject to an isotropic topological \textit{alignment force} with all the others, i.e., they tend to have the same velocity of the group mates. Assuming that the agents' positions are close enough, this corresponds to the tendency to go where group mates go (herding effect).
\end{itemize}


\section{The microscopic model}\label{sec:micromodel}
In this section we introduce the microscopic model for followers and leaders. We denote by $d$ the dimension of the space in which the motion takes place (typically $d=2$), by $\NF$ the number of followers and by $\NL \ll \NF$ the number of leaders. We also denote by $\Omega\equiv\R^d$ the walking area and by $\target\in\Omega$ the target point. To define the target's visibility area, we consider the set $\Sigma$, with $\target\in\Sigma\subset\Omega$, and we assume that the target is completely visible from any point belonging to $\Sigma$ and completely invisible from any point belonging to $\Omega\backslash\Sigma$. 

For every $i=1,\ldots,\NF$, let $(\xf_i(t),\vf_i(t))\in\R^{2d}$ denote position and velocity of the agents belonging to the population of followers at time $t\geq 0$ and, for every  $k=1,\ldots,\NL$, let $(\xl_k(t),\vl_k(t))\in\R^{2d}$  denote position and velocity of the agents among the population of leaders at time $t \geq 0$. Let us also define $\mathbf{\xf}:=(\xf_1,\ldots,\xf_{\NF})$ and $\mathbf{\xl}:=(\xl_1,\ldots,\xl_{\NL})$.

Finally, let us denote by $B_r(x)$ the ball of radius $r>0$ centered at $x\in\Omega$ and by $\mathcal B_\mathcal N(x;\mathbf{\xf},\mathbf{\xl})$ the \emph{minimal} ball centered at $x$ encompassing at least $\mathcal N$ agents, and by $\mathcal N^*$ the actual number of agents in $\mathcal B_\mathcal N(x;\mathbf{\xf},\mathbf{\xl})$. Note that $\mathcal N^*\geq\mathcal N$. 
\begin{remark}
The computation of $\mathcal B_\mathcal N(x;\mathbf{\xf},\mathbf{\xl})$ requires the knowledge of the positions of all the agents, since all the distances $|\xf_i-x|$, $i=1,\ldots,\NF$, and  $|\xl_k-x|$, $k=1,\ldots,\NL$ must be evaluated in order to find the $\mathcal N$ closest agents to $x$.
\end{remark}

\medskip

The microscopic dynamics described by the two populations is given by the following set of ODEs: for $i = 1, \dots, \NF$ and $k = 1, \ldots, \NL$,
\begin{equation}\label{eq:micro}
\left\{
\begin{array}{l}
\dot{x}_i = \vf_i,\\ [1.5mm] 
\dot{v}_i = \Ad(\xf_i,\vf_i) + \sum_{j=1}^{\NF} \HFF(\xf_i,\vf_i,\xf_j,\vf_j;\mathbf{\xf},\mathbf{\xl}) +\sum_{\ell=1}^{\NL} \HFL(\xf_i,\vf_i,\xl_\ell,\vl_\ell;\mathbf{\xf},\mathbf{\xl}),\\ [1.5mm]
\dot{y}_k= w_k = \sum_{j=1}^{\NF} \HLF(\xl_k,\xf_j) +  \sum_{\ell=1}^{\NL} \HLL(\xl_k,\xl_\ell) + u_k.
\end{array}
\right.
\end{equation}
We assume that

\begin{itemize}
\item $\Ad$ is a self-propulsion term, given by the relaxation toward a random direction or the relaxation toward a unit vector pointing to the target (the choice depends on the position), plus a term which translates the tendency to reach a given characteristic speed $\barS \geq 0$ (modulus of the velocity), i.e.,
\begin{align} \label{eq:Ad}
\Ad(x,v) := \theta(x) \Crand(\Wz-v) & + (1-\theta(x))\Ct\left(\frac{\target - x}{|\target - x|} - v\right)+\Cchar(\barS^2-|v|^2)v,
\end{align}
where $\theta:\mathbb{R}^d \rightarrow [0,1]$ is the characteristic function of $\Omega\backslash\Sigma$, $\theta(x) =\chi_{\Omega\backslash\Sigma}(x)$, 
$\Wz$ is a $d$-dimensional random vector with normal distribution $\mathcal N(0,\sigma^2)$, and $\Crand$, $\Ct$, $\Cchar$ are positive constants.

\item
The interactions follower-follower and follower-leader coincide and are equal to 
\begin{equation}
\begin{aligned}\label{eq:HFF}
&\HFF(x,v,x',v';\mathbf \xf,\mathbf \xl) := 
-\Cr R_{\gamma,r}(x,x') + 
\theta(x)\frac{\Ca}{\mathcal N^*}\left(v'-v\right) \chi_{B_\mathcal{N}(x;\mathbf \xf,\mathbf \xl)}(x'),\\
&\HFL(x,v,y,w;\mathbf \xf,\mathbf \xl):=\HFF(x,v,y,w;\mathbf \xf,\mathbf \xl),
\end{aligned}
\end{equation}
for given positive constants $\Cr, \Ca, r,$ and $\gamma$, and
\begin{align*}
R_{\gamma,r}(x,x') & = \begin{cases}
e^{-|x'-x|^\gamma}\frac{x'-x}{|x'-x|} & \text{ if } x'\in B_r(x)\backslash\{x\}, \\
0 & \text{ otherwise,}
\end{cases}
\end{align*}
models a (metrical) repulsive force, while the second term accounts for the (topological) alignment force, which vanishes inside $\Sigma$.
Note that, once the summations over the agents $\sum_j \HFF$ and $\sum_\ell\HFL$ are done, the alignment term models the tendency of the followers to relax toward the average velocity of the $\mathcal{N}$ closer agents.
With the choice $\HFF\equiv \HFL$ the leaders are not recognized by the followers as special. This feature opens a wide range of new applications, including the control of crowds not prone to follow authority's directives.

\item
The interactions leader-follower and leader-leader reduce to a mere (metrical) repulsion, i.e., $\HLF = \HLL = -\Crl R_{\zeta,r}$, where $\Crl>0$ and $\zeta>0$ are in general different from $\Cr$ and $\gamma$, respectively. Note that here the repulsion force is interpreted as a velocity field, while for followers it is an acceleration field. 

\item
$u_k:\R^+\to\R^{d\NL}$ is the control variable, to be chosen in a set of admissible control functions.
\end{itemize}
\begin{remark}
The behaviour of the leaders is entirely encapsulated in the control term $u$ but for a short-range repulsion force. The latter should be indeed interpreted as a force due to the presence of the others, thereby non controllable.
\end{remark}

\section{Formal derivation of a Boltzmann-type equation}\label{sec:mesomodel}

As already mentioned, our main interest in \eqref{eq:micro} lies in the case $\NL \ll \NF$, that is the population of followers exceeds by far the one of leaders. When $\NF$ is very large, a microscopic description of both populations is no longer a viable option. We thus consider the evolution of the distribution of followers at time $t \geq 0$, denoted by $f(t,x,v)$, together with the microscopic equations for the leaders (whose number is still small). To this end, we denote with $\rhoF$ the total mass of followers, i.e.,
\begin{align*}
\rhoF(t) = \int_{\mathbb{R}^{2d}} f(t,x,v) \ dx \ dv,
\end{align*}
which we shall eventually require to be equal to $\NF$. We introduce, for symmetry reasons, the distribution of leaders $g$ and their total mass
\begin{equation}\label{eq:deltaleaders} 
g(t,x,v) = \sum^{\NL}_{k = 1} \delta_{(\xl_k(t),\vl_k(t))}(x,v), 
\qquad 
\rhoL(t) = \int_{\mathbb{R}^{2d}} g(t,x,v) \ dx \ dv = \NL.
\end{equation}

The evolution of $f$ can be then described by a Boltzmann-type dynamics, derived from the above instantaneous control formulation, which is obtained by analyzing the binary interactions between a follower and another follower and the same follower with a leader. The application of standard methods of binary interactions, see \cite{cercignani1994}, shall yield a mesoscopic model for the distribution of followers, to be coupled with the previously presented ODE dynamics for leaders.

To derive the Boltzmann-type dynamics, we assume that, before interacting, each agent has at his disposal the values $\mathbf{\xf}$ and $\mathbf{\xl}$ that he needs in order to perform its movement: hence, in a binary interaction between two followers with state parameter $(x,v)$ and $(\hat{x},\hat{v})$, the value of $\HFF(x,v,\hat{x},\hat{v};\mathbf{\xf},\mathbf{\xl})$ does not depend on $\mathbf{\xf}$ and $\mathbf{\xl}$. In the case of $\HFF$ of the form \eqref{eq:HFF}, this means that the ball $B_\mathcal{N}(x;\mathbf \xf,\mathbf \xl)$ and the value of $\mathcal{N}^*$ have been already computed before interacting.

Moreover, since we are considering the distributions $f$ and $g$ of followers and leaders, respectively, the vectors $\mathbf{\xf}$ and $\mathbf{\xl}$ are derived from $f$ and $g$ by means of the first marginals of $f$ and $g$, $\pi_1 f$ and $\pi_1 g$, respectively, which give the spatial variables of those distribution. Hence, since no confusion arises, we write $\HFF(x,v,\hat{x},\hat{v};\pi_1 f,\pi_1 g)$ in place of  $\HFF(x,v,\hat{x},\hat{v};\mathbf{\xf},\mathbf{\xl})$ to stress the dependence of this term on $f$ and $g$.

We thus consider two followers with state parameter $(x,v)$ and $(\hat{x},\hat{v})$ respectively, and we describe the evolution of their velocities after the interaction according to
\begin{align}
\begin{cases}
\begin{split}\label{eq:follcoord}
v^* & = v + \etaF \left[\theta(x) \noisepar \xi + \Adcut(x,v) + \rhoF\HFF(x,v,\hat{x},\hat{v};\pi_1 f,\pi_1 g) \right],\\ 
\hat{v}^* & = \hat{v} + \etaF \left[\theta(\hat{x}) \noisepar \xi + \Adcut(\hat{x},\hat{v}) + \rhoF\HFF(\hat{x},\hat{v},x,v;\pi_1 f,\pi_1 g) \right],
\end{split}
\end{cases}
\end{align}
where $\etaF$ is the strength of interaction among followers, $\xi$ is a random variables whose entries are i.i.d.\ following a normal distribution with mean $0$, variance $\varsigma^2$ (which shall be related to the variance $\sigma^2$ of the random vector $\Wz$ in \eqref{eq:Ad} 
in such a way that we recover from \eqref{eq:follcoord} a Fokker-Planck-type dynamics for $f$), taking values in a set $\mathcal{B}$, and $\Adcut$ is defined as the deterministic part of the self-propulsion term \eqref{eq:Ad},
\begin{align}\label{eq:Adcut}
\Adcut(x,v)=-\theta(x)\Crand v+(1-\theta(x))\Ct\left(\frac{\target - x}{|\target - x|} - v\right)+\Cchar(\barS^2-|v|^2)v.
\end{align}
We then consider the same follower as before with state parameters $(x,v)$ and a leader agent $(\tilde{x},\tilde{v})$; in this case the modified velocities satisfy
\begin{align}
\begin{cases}
\begin{split}\label{eq:leadcoord}
& v^{**}  = v + \etaL \rhoL \HFL(x,v,\tilde{x},\tilde{v};\pi_1 f,\pi_1 g),\\
& \tilde{v}^*  = \tilde{v}, 
\end{split}
\end{cases}
\end{align}
where $\etaL$ is the strength of the interaction between followers and leaders. Note that \eqref{eq:leadcoord} accounts only the change of the followers' velocities, since leaders are not evolving via binary interactions.

Since we are interested in studying the problem in the widest possible framework, and avoid sticking only to the  previous choice of the functions $\Adcut$, $\HFF$ and $\HFL$, in what follows we assume that 
\medskip
\begin{description}
\item[(Inv)] the systems \eqref{eq:follcoord} and \eqref{eq:leadcoord} constitute invertible changes of variables from $(v,\hat{v})$ to $(v^*,\hat{v}^{*})$ and from $(v,\tilde{v})$ to $(v^{**},\tilde{v}^{*})$, respectively.
\end{description}
\medskip
The time evolution of $f$ is then given by a balance between bilinear gain and loss of space and velocity terms according to the two binary interactions \eqref{eq:follcoord} and \eqref{eq:leadcoord}, quantitatively described by the following Boltzmann-type equation
\begin{align}\label{eq:strongBoltz}
\partial_t f+ v \cdot \nabla_x f = \lamF Q(f,f) + \lamL Q(f,g),
\end{align}
where $\lamF$ and $\lamL$ stand for the interaction frequencies among followers and between followers and leaders, respectively. The interaction integrals $Q(f,f)$ and $Q(f,g)$ are defined as
\begin{align*}
Q(f,f)(t) & = \mathbb{E}\left(\int_{\mathbb{R}^{4d}}\left(\frac{1}{J_{\textsc{f}}} f(t,x_*,v_*)f(t,\hat{x}_*,\hat{v}_*) - f(t,x,v)f(t,\hat{x},\hat{v})\right) \ d\hat{x} \ d\hat{v}\right), \\
Q(f,g)(t) & = \mathbb{E}\left(\int_{\mathbb{R}^{4d}}\left(\frac{1}{J_{\textsc{l}}} f(t,x_{**},v_{**})g(t,\tilde{x}_*,\tilde{v}_*) - f(t,x,v)g(t,\tilde{x},\tilde{v})\right)  \ d\tilde{x} \ d\tilde{v}\right),
\end{align*}
where the couples $(x_*,v_*)$ and $(\hat{x}_*,\hat{v}_*)$ are the pre-interaction states that generates $(x,v)$ and $(\hat{x},\hat{v})$ via \eqref{eq:follcoord}, and $J_{\textsc{f}}$ is the Jacobian of the change of variables given by \eqref{eq:follcoord} (well-defined by (Inv)). Similarly, $(x_{**},v_{**})$ and $(\tilde{x}_*,\tilde{v}_*)$ are the pre-interaction states that generates $(x,v)$ and $(\tilde{x},\tilde{v})$ via \eqref{eq:leadcoord}, and $J_{\textsc{l}}$ is the Jacobian of the change of variables given by \eqref{eq:leadcoord}. Moreover, the expected value $\mathbb{E}$ is computed with respect to $\xi \in \mathcal{B}$.

In what follows, for the sake of compactness, we shall omit the time dependency of $f$ and $g$, and hence of $Q(f,f)$ and $Q(f,g)$ too.

\begin{remark}
If we would have opted for a description of agents as hard-sphere particles, the arising Boltzmann equation \eqref{eq:strongBoltz} would be of Enskog type, see \cite{toscanibellomoenskog}. The relationship between the hard- and soft-sphere descriptions (i.e., where repulsive forces are considered, instead) has been deeply discussed, for instance, in \cite{andersenhard}. In our model, the repulsive force $R_{\gamma,r}$ is not singular at the origin for computational reasons, therefore the parameters $\gamma$ and $r$ have to be chosen properly to avoid arbitrary high density concentrations. 
\end{remark}

\subsection{Notations and basic definitions}

We shall start the analysis of equation \eqref{eq:strongBoltz} by fixing some notation and terminology. First of all, for any $\beta \in \mathbb{N}^d$ we set $|\beta| = \sum^d_{i = 1} \beta_i$, and for any function $h \in C^p(\R^d \times \R^d,\R)$, with $p\geq0$ and any $\beta \in \mathbb{N}^d$ such that $|\beta| \leq p$, we define for every $(x,v) \in \R^d \times \R^d$
\begin{align*}
\partial^{\beta}_v h(x,v) := \frac{\partial^{|\beta|} h}{\partial^{\beta_1}v_1 \cdots \partial^{\beta_d}v_d} (x,v),
\end{align*}
with the convention that if $\beta = (0, \ldots, 0)$ then $\partial^{\beta}_v h(x,v) := h(x,v)$.
\medskip
\begin{definition}[Test functions] We denote with $\mathcal{T}_{\delta}$ the set of compactly supported functions $\varphi$ from $\R^{2d}$ to $\R$ such that for any multi-index $\beta  \in \mathbb{N}^d$ we have:
\begin{enumerate}
\item if $|\beta| < 2$, then $\partial^{\beta}_v \varphi (x, \cdot)$ is continuous for every $x \in \R^d$;
\item if $|\beta| = 2$, then there exists $M > 0$ such that:
\begin{enumerate}
\item $\partial^{\beta}_v \varphi (x, \cdot)$ is uniformly H\"older continuous of order $\delta$ for every $x \in \R^d$ with H\"older bound $M$, that is for every $x \in \R^d$ and for every $v,w  \in \R^d$ it holds
\begin{align*}
\left\| \partial^{\beta}_v \varphi (x, v) - \partial^{\beta}_v \varphi (x, w) \right\| \leq M \left\| v - w \right\|^{\delta};
\end{align*}
\item $\|\partial^{\beta}_v \varphi (x, v)\| \leq M$ for every $(x,v) \in \R^{2d}$.
\end{enumerate}
\end{enumerate}
\end{definition}
\medskip
Notice that $\mathcal{C}^{\infty}_c(\R^d \times \R^d;\R) \subseteq \mathcal{T}_{\delta}$ for every $0 < \delta \leq 1$.

Let $\mathcal{M}_0(\R^{2d})$ denote the set of measures taking values in $\R^{2d}$. For any two function $f$ and $\varphi$ from $\R^{2d}$ to $\R$ for which the integral below is well-defined (in particular, when $f \in \mathcal{M}_0(\R^{2d})$ and $\varphi \in \mathcal{T}_{\delta}$), we set
$$\left \langle f, \varphi \right \rangle := \int_{\R^{2d}} \varphi(x,v) f(x,v) \ dx \ dv.$$

Finally, we say that $(f,\xl_1,\vl_1, \ldots, \xl_{\NL},\vl_{\NL}) \in \mathcal{M}_0(\R^{2d}) \times \R^{2d\NL}$ is \emph{admissible} if the following quantities exist and are finite:
\begin{enumerate}[$(i)$]
\item $\int_{\R^{2d}} \Adcut(x,v) f(x,v) \ dx \ dv,$
\item $\int_{\R^{4d}} \HFF(x,v,\hat{x},\hat{v};\pi_1 f,\mathbf{\xl}) f(x,v) f(\hat{x},\hat{v}) \ dx \ dv \ d\hat{x} \ d\hat{v},$
\item $\sum^{\NL}_{\ell = 1} \int_{\R^{2d}} \HFL(x,v,\xl_{\ell},\vl_{\ell};\pi_1 f,\mathbf{\xl}) f(x,v) \ dx \ dv,$
\item $\int_{\R^{2d}} \HLF(\xl_k,x) f(x,v) \ dx \ dv$ for every $k = 1,\ldots,\NL$,
\item $\sum_{\ell = 1}^{\NL} \HLL(\xl_k,\xl_{\ell})$ for every $k = 1,\ldots,\NL$.
\end{enumerate}

Consequently, we introduce our definition of solution for the combined ODE-PDE system for the dynamics of microscopic leaders and mesoscopic followers.

\begin{definition} \label{def:PDE}
Fix $T > 0$, $\delta > 0$, and $u: [0,T] \rightarrow \mathbb{R}^{d\NL}$. By a \emph{$\delta$-weak solution} of the initial value problem for the equation
\begin{align} \label{eq:micromeso}
\left\{
\begin{array}{l}
\partial_t f + v \cdot \nabla_x f = \lamF Q(f,f) + \lamL Q(f,g), \\
\displaystyle\dot{\xl}_k = \vl_k = \int_{\R^{2d}} \HLF(\xl_k,x) f(x,v) \ dx\ dv + \sum_{\ell = 1}^{\NL} \HLL(\xl_k,\xl_{\ell}) + u_k,
\end{array}
\right.
\end{align}
corresponding to the control $u$ and the initial datum $(f^0, \xl^0_1, \ldots, \xl^0_{\NL}) \in \mathcal{M}_0(\R^{2d}) \times \R^{d\NL}$ in the interval $[0,T]$, we mean any $(f, \xl_1, \ldots, \xl_{\NL}) \in L^2([0,T], \mathcal{M}_0(\R^{2d})) \times C^1([0,T],\R^{d\NL})$ such that:
\begin{enumerate}
\item $f(0,x,v) = f^0(x,v)$ for every $(x,v) \in \R^{2d}$;
\item the vector $(f(t),\xl_1(t),\dot{\xl}_1(t), \ldots, \xl_{\NL}(t),\dot{\xl}_{\NL}(t))$ is admissible for every $t \in [0,T]$;
\item there exists $R_T > 0$ such that $\supp(f(t)) \subset B_{R_T}(0)$ for every $t \in [0,T]$;
\item $\rhoF(t) = \NF$ for every $t \in [0,T]$;
\item $f$ satisfies the weak form of the equation \eqref{eq:strongBoltz}, i.e.,
\begin{align*}
\frac{\partial}{\partial t} \left \langle f, \varphi \right\rangle + \left \langle f, v \cdot \nabla_x \varphi \right\rangle =  \lamF\left\langle Q(f,f),\varphi \right\rangle + \lamL \left\langle Q(f,g),\varphi \right\rangle,
\end{align*}
for all $t \in (0,T]$ and all $\varphi \in \mathcal{T}_{\delta}$, where
\begin{align}
\left\langle Q(f,f),\varphi \right\rangle & = \mathbb{E}\left(\int_{\mathbb{R}^{4d}} \left(\varphi(x,v^{*}) - \varphi(x,v) \right)f(x,v)f(\hat{x},\hat{v}) \ dx \ dv \ d\hat{x} \ d\hat{v}\right), \label{eq:collopF} \\
\left\langle Q(f,g),\varphi \right\rangle & = \mathbb{E}\left(\int_{\mathbb{R}^{4d}} \left(\varphi(x,v^{**}) - \varphi(x,v) \right)f(x,v)g(\tilde{x},\tilde{v}) \ dx \ dv \ d\tilde{x} \ d\tilde{v} \right), \label{eq:collopL}
\end{align}
and where $v^*$ and $v^{**}$ are given by $\eqref{eq:follcoord}$ and $\eqref{eq:leadcoord}$, respectively;
\item for every $k = 1, \ldots, \NL$, $\xl_k$ satisfies in the Carath\'eodory sense
\begin{align*}
\dot{\xl}_k = \vl_k = \int_{\R^{2d}} \HLF(\xl_k,x) f(x,v) \ dx\ dv + \sum\limits_{\ell = 1}^{\NL} \HLL(\xl_k,\xl_{\ell}) + u_k.
\end{align*}
\end{enumerate}
\end{definition}

\subsection{The grazing interaction limit}
In order to obtain a more regular operator than the Boltzmann operator in \eqref{eq:strongBoltz}, we introduce the so called grazing interaction limit. 
This technique, analogous to the so-called grazing collision limit in plasma physics, has been thoroughly studied in \cite{Vil02} and allows, as pointed out in \cite{PT:13}, to pass the binary Boltzmann description introduced in the previous section to the mean-field limit. This process, while somehow losing the microscopic interaction rules \eqref{eq:follcoord}--\eqref{eq:leadcoord}, allows a better understanding of the solution behavior for large times, see \cite{albi1405.0736, toscani2006kinetic}. 
Moreover, the regularized operator of Fokker-Planck-type will directly show the connection with the microscopic model.

In what follows, we shall assume that our agents densely populate a small region but weakly interact with each other. Formally, we assume that the interaction strengths $\etaF$ and $\etaL$ scale according to a parameter $\varepsilon$, the interaction frequencies $\lamF$ and $\lamL$ scale as $1/\varepsilon$, and we let $\varepsilon \rightarrow 0$. In order to avoid losing the diffusion term in the limit, we also scale the variance of the noise term $\varsigma^2$ as $1/\varepsilon$. More precisely, we set
\begin{equation}\label{eq:scaling}
\etaF = \varepsilon, \qquad \etaL = \varepsilon,  \qquad
\lamF = \frac{1}{\varepsilon \rhoF},  \qquad \lamL  = \frac{1}{\varepsilon \rhoL}, \qquad \varsigma^2 = \frac{\sigma^2}{\varepsilon}. 
\end{equation}
We now study the Boltzmann equation \eqref{eq:strongBoltz} to see if simplifications may occur under the above scaling assumptions.

Let us fix $T > 0$, $\delta > 0$, $\varepsilon>0$, a control $u:[0,T] \rightarrow \R^{d\NL}$, and an initial datum $(f^0,\xl^0_1, \ldots, \xl^0_{\NL}) \in \mathcal{M}_0(\R^{2d}) \times \R^{d\NL}$. We consider a $\delta$-weak solution $(f,\xl_1, \ldots, \xl_{\NL})$ of system \eqref{eq:micromeso} with control $u$, initial datum $(f^0,\xl^0_1, \ldots, \xl^0_{\NL})$ and define $g$ as in \eqref{eq:deltaleaders}. Let the scaling \eqref{eq:scaling} hold for the chosen $\varepsilon$.
Following the ideas in \cite{carrillo2010particle,cordier2005kinetic}, we can expand $\varphi(x,v^{*})$ inside \eqref{eq:collopF} in Taylor' series of $v^* - v$ up to the second order, to get
\begin{align*}
&\left\langle Q(f,f),\varphi \right\rangle  = \mathbb{E}\Bigg(\int_{\mathbb{R}^{4d}} \nabla_v\varphi(x,v) \cdot \left(v^{*} - v\right)f(x,v)f(\hat{x},\hat{v}) \ dx \ dv \ d\hat{x} \ d\hat{v} \\
 \quad& + \frac{1}{2}\int_{\mathbb{R}^{4d}} \left[\sum^{d}_{i,j = 1} \partial^{(i,j)}_{v} \varphi(x,v) \left(v^{*} - v\right)_i\left(v^{*} - v\right)_j\right] f(x,v)f(\hat{x},\hat{v}) \ dx \ dv \ d\hat{x} \ d\hat{v}\Bigg) + R_{\varphi}^{\textsc{f}},
\end{align*}
where the remainder $R^{\textsc{f}}_{\varphi}$ of the Taylor expansion has the form
\begin{eqnarray*}
R^{\textsc{f}}_{\varphi} & = \mathbb{E}\Bigg(\frac{1}{2}\int_{\mathbb{R}^{4d}} \left[\sum^{d}_{i,j = 1} \left(\partial^{(i,j)}_{v} \varphi(x,v) - \partial^{(i,j)}_{v} \varphi(x,\overline{v}^{\textsc{f}})\right) \left(v^{*} - v\right)_i\left(v^{*} - v\right)_j\right] \\
& \quad f(x,v)f(\hat{x},\hat{v}) \ dx \ dv \ d\hat{x} \ d\hat{v}\Bigg),
\end{eqnarray*}
with $\overline{v}^{\textsc{f}} = \thetaF v^* + (1 - \thetaF) v$, for some $\thetaF \in [0,1]$. By using the substitution given by the interaction rule \eqref{eq:follcoord}, i.e.
$$v^* - v = \etaF \left[\theta(x) \noisepar \xi + \Adcut(x,v) + \rhoF\HFF(x,v,\hat{x},\hat{v};\pi_1 f,\pi_1 g) \right],$$
we obtain
\begin{eqnarray}
\int_{\mathbb{R}^{4d}} \nabla_v\varphi(x,v) & \cdot & \left(v^{*} - v\right)f(x,v)f(\hat{x},\hat{v}) \ dx \ dv \ d\hat{x} \ d\hat{v} \nonumber \\
& = & \etaF \int_{\mathbb{R}^{4d}} \nabla_v\varphi(x,v) \cdot \left[\theta(x) \noisepar \xi + \Adcut(x,v) + \rhoF\HFF(x,v,\hat{x},\hat{v};\pi_1 f,\pi_1 g)\right] \nonumber \\
& & \qquad\quad f(x,v)f(\hat{x},\hat{v}) \ dx \ dv \ d\hat{x} \ d\hat{v} \nonumber  \\
& = & \etaF\rhoF\left\langle  f ,\nabla_v \varphi \cdot  \theta \noisepar \xi \right\rangle + \etaF \rhoF \left\langle f ,\nabla_v\varphi \cdot  \left(\Adcut + \mathcal{H}^{\textsc{f}}[f]\right) \right\rangle, \label{eq:todomean} 
\end{eqnarray}
having denoted with
\begin{align*}
\mathcal{H}^{\textsc{f}}[f](x,v) & = \int_{\mathbb{R}^{2d}} \HFF(x,v,\hat{x},\hat{v};\pi_1 f,\pi_1 g) f(\hat{x},\hat{v}) \ d\hat{x} \ d\hat{v}.
\end{align*}
By taking the expected value w.r.t.\ $\xi \in \mathcal{B}$ of \eqref{eq:todomean}, the term $\etaF\rhoF\left\langle  f ,\nabla_v \varphi \cdot  \theta \noisepar \xi \right\rangle$ is canceled out, since $\mathbb{E}\left(\xi\right) = 0$.

Writing (whenever possible) $\HFF$, $\HFL$, $\theta$ and $\Adcut$ in place of $\HFF(x,v,\hat{x},\hat{v};\pi_1 f,\pi_1 g)$, $\HFL(x,v,\tilde{x},\tilde{v};\pi_1 f,\pi_1 g)$, $\theta(x)$ and $\Adcut(x,v)$, respectively, the same substitution yields for the second order term
\begin{align*}
\frac{1}{2}\int_{\mathbb{R}^{4d}} & \left[\sum^{d}_{i,j = 1} \partial^{(i,j)}_{v}\varphi(x,v) \left(v^{*} - v\right)_i\left(v^{*} - v\right)_j\right]f(x,v)f(\hat{x},\hat{v}) \ dx \ dv \ d\hat{x} \ d\hat{v} \\
& = \frac{1}{2}(\etaF)^2\int_{\mathbb{R}^{4d}} \Bigg[\sum^{d}_{i,j = 1} \partial^{(i,j)}_{v}\varphi(x,v) \left(\theta \noisepar \xi + \Adcut + \rhoF \HFF\right)_i \\
& \quad \left(\theta \noisepar \xi + \Adcut + \rhoF \HFF\right)_j\Bigg] f(x,v)f(\hat{x},\hat{v}) \ dx \ dv \ d\hat{x} \ d\hat{v}.
\end{align*}
Notice that if we compute the expected value of the expression above, all the cross terms $\xi_i \xi_j$ vanish, since they are drawn independently from each other. Hence we obtain
\begin{align*}
\lamF \left\langle Q(f,f),\varphi \right\rangle & = \lamF\etaF \rhoF \left\langle f ,\nabla_v\varphi \cdot  \left(\Adcut + \mathcal{H}^{\textsc{f}}[f]\right) \right\rangle + \frac{1}{2} \lamF\left(\etaF\right)^2\rhoF \varsigma^2 \left\langle f,(\theta \noisepar)^2 \Delta_v \varphi \right\rangle \\
& \quad + \frac{1}{2} \lamF\left(\etaF\right)^2 \RestF[f,f] + \lamF R^{\textsc{f}}_{\varphi},
\end{align*}
where the notation $\RestF[f,f]$ stands for
\begin{eqnarray*}
\RestF[f,f] & = \int_{\mathbb{R}^{4d}} \Bigg[\sum^d_{i,j = 1}\partial^{(i,j)}_v\varphi \left(\Adcut + \rhoF \HFF\right)_i\left(\Adcut + \rhoF \HFF\right)_j \Bigg]f(x,v)f(\hat{x},\hat{v}) \ dx \ dv \ d\hat{x} \ d\hat{v} .
\end{eqnarray*}
By means of the same computations, we can derive a similar expression for $\left\langle Q(f,g),\varphi \right\rangle$. Indeed, from \eqref{eq:collopL} we have
\begin{align*}
\lamL\left\langle Q(f,g),\varphi \right\rangle & = \lamL\etaL\rhoL \left\langle f ,\nabla_v\varphi\cdot  \mathcal{H}^{\textsc{l}}[g]\right\rangle + \frac{1}{2} \lamL\left(\etaL\right)^2 \RestL[f,g] + \lamL R^{\textsc{l}}_{\varphi},
\end{align*}
having denoted with
\begin{eqnarray*}
\mathcal{H}^{\textsc{l}}[g](x,v) = \int_{\R^{2d}} \HFL(x,v,\tilde{x},\tilde{v};\pi_1 f, \pi_1 g) g(\tilde{x},\tilde{v}) \ d\tilde{x} \ d\tilde{v},
\end{eqnarray*}
\begin{eqnarray*}
\RestL [f,g] & = \left(\rhoL\right)^2 \int_{\mathbb{R}^{4d}} \Bigg[\sum^d_{i,j = 1}\partial^{(i,j)}_{v}\varphi(x,v)\HFL_i \HFL_j \Bigg] f(x,v) g(\tilde{x},\tilde{v}) \ dx \ dv \ d\tilde{x} \ d\tilde{v},
\end{eqnarray*}
\begin{eqnarray*}
R^{\textsc{l}}_{\varphi} & = \mathbb{E}\Bigg(\frac{1}{2}\int_{\mathbb{R}^{4d}} \left[\sum^{d}_{i,j = 1} \left(\partial^{(i,j)}_{v}\varphi(x,v) -\partial^{(i,j)}_{v}\varphi(x,\overline{v}^{\textsc{l}})\right) \left(v^{**} - v\right)_i \left(v^{**} - v\right)_j \right] \\
& \quad f(x,v)g(\tilde{x},\tilde{v}) \ dx \ dv \ d\tilde{x} \ d\tilde{v}\Bigg),
\end{eqnarray*}
with $\overline{v}^{\textsc{l}} = \vartheta^{\textsc{l}} v^{**} + (1 - \vartheta^{\textsc{l}}) v$ for some $\vartheta^{\textsc{l}} \in [0,1]$.

If we now apply the rescaling rules \eqref{eq:scaling}, the following simplifications occur
\begin{equation*}
 \lamF \etaF \rhoF = 1, 
 \quad  \lamL \etaL \rhoL= 1, 
 \quad \lamF\left(\etaF\right)^2 \varsigma^2 \rhoF = \sigma^2, 
 \quad  \lamL\left(\etaL\right)^2 = \varepsilon, 
 \quad \lamL\left(\etaL\right)^2 = \varepsilon.
\end{equation*}
Hence, if we let $\varepsilon \rightarrow 0$ and we assume that, for every $\varphi \in \mathcal{T}_{\delta}$,
\begin{align} \label{eq:limRest}
\lim_{\varepsilon \rightarrow 0} \left(\frac{\varepsilon}{2} \left\|\RestF [f,f]\right\| + \frac{\varepsilon}{2} \left\| \RestL [f,g]\right\| + \frac{1}{\varepsilon}\left\|R^{\textsc{f}}_{\varphi}\right\| + \frac{1}{\varepsilon}\left\|R^{\textsc{l}}_{\varphi}\right\| \right) = 0
\end{align}
holds true, we obtain the weak formulation of a Fokker-Planck-type equation for the followers' dynamics
\begin{align} \label{eq:FokkerPlanck}
\frac{\partial}{\partial t} \left \langle f, \varphi \right\rangle + \left \langle f, v \cdot \nabla_x \varphi \right\rangle =  \left\langle f, \nabla_v \varphi \cdot \mathcal{G}\left[f,g\right] + \frac{1}{2}\sigma^2 (\theta \noisepar)^2 \Delta_v \varphi \right\rangle,
\end{align}
where
\begin{align*}
\mathcal{G}\left[f,g\right] & = \Adcut + \mathcal{H}^{\textsc{f}}[f] + \mathcal{H}^{\textsc{l}}[g].
\end{align*}
We leave the details of the proof of the limit \eqref{eq:limRest} to the following Section \ref{sec:limRest}.

Since $\varphi$ has compact support, equation \eqref{eq:FokkerPlanck} can be recast in strong form by means of integration by parts. Coupling the resulting PDE with the microscopic ODEs for the leaders  $k=1,\ldots,\NL$, we eventually obtain the system
\begin{equation}\label{eq:FokkerPlanckStrong}
\left\{ 
\begin{array}{l}
\partial_t f + v \cdot \nabla_x f = - \nabla_v \cdot \left(\mathcal{G}\left[f,g\right]f\right) + \frac{1}{2} \sigma^2 (\theta \noisepar)^2 \Delta_v f,\\ [1mm]
\displaystyle
\dot{\xl}_k = \vl_k = \int_{\R^{2d}} \HLF(\xl_k,x) f(x,v) \ dx\ dv + \sum_{\ell = 1}^{\NL} \HLL(\xl_k,\xl_{\ell}) + u_k.
\end{array}
\right.
\end{equation}
\begin{remark}
Observe that, assuming $f$ in \eqref{eq:FokkerPlanckStrong} to be  the empirical measure $f= \sum^{\NF}_{i = 1} \delta_{(x_i(t),v_i(t))}$ concentrated on the trajectories $(x_i(t),v_i(t))$ of the microscopic dynamics \eqref{eq:micro}, we recover the original microscopic model itself.
\end{remark}

\begin{remark}
In terms of model hierarchy one could imaging to compute the moments of \eqref{eq:FokkerPlanckStrong} in order to further reduce the complexity. Let us stress that deriving a consistent macroscopic system from the kinetic equation is in general a difficult task, since equilibrium states are difficult to obtain, therefore no closure of the moments equations is possible. For self-organizing models similar to \eqref{eq:FokkerPlanckStrong}, in the noiseless case (i.e.\ $\sigma\equiv0$), a standard way to obtain a closed hydrodynamic system is to assume the velocity distribution to be mono-kinetic, i.e.\ $f(t,x,v)=\rho(t,x)\delta(v-V(t,x))$, and the fluctuations to be negligible, thus computing the moments of \eqref{eq:FokkerPlanckStrong} leads to the following macroscopic system for the density $\rho$ and the bulk velocity $V$,
\begin{equation}
\left\{ 
\begin{array}{l}
\partial_t \rho + \nabla_x\cdot (\rho V) = 0,\\ [2mm]
\partial_t (\rho V)+ \nabla_x\cdot (\rho V\otimes V) = \mathcal{G}_m\left[\rho,\rho^{\textsc{l}},V,V^{\textsc{l}}\right]\rho,\\ [1mm]
\displaystyle
\dot{\xl}_k = \vl_k = \int_{\R^{d}} \HLF(\xl_k,x) \rho(t,x) \ dx+ \sum_{\ell = 1}^{\NL} \HLL(\xl_k,\xl_{\ell}) + u_k,
\end{array}
\right.
\end{equation}
where $\rho^{\textsc{l}}(x,t), V^{\textsc{l}}(x,t)$ represent the leaders' macroscopic density and bulk velocity, respectively, and $\mathcal{G}_m$ the macroscopic interaction operator, see \cite{albi2013AML,carrillo2010particle} for further details. 
For $\sigma>0$, the derivation of a macroscopic system depends highly on the scaling regime between the noise and the interaction terms, see for example \cite{carrillo2009double,carrillo2010particle,karper2015hydrodynamic}. Furthermore, the presence of diffusion operator in model \eqref{eq:FokkerPlanckStrong} depends on the spatial domain, therefore the derivation of a reasonable macroscopic model is not trivial and it is left for further studies.
\end{remark}

\subsection{Estimates for the remainder terms}\label{sec:limRest}
Motivated by the results in \cite{toscani2006kinetic}, we shall estimate the quantity $\left\|R^{\textsc{f}}_{\varphi}\right\|$ as follows: since $\overline{v}^{\textsc{f}} = \vartheta^{\textsc{f}} v^* + (1 - \vartheta^{\textsc{f}}) v$ implies $\left\| v - \overline{v}^{\textsc{f}} \right\| \leq \vartheta^{\textsc{f}} \left\| v - v^* \right\|$, then for every $\varphi \in \mathcal{T}_{\delta}$ it follows that
\begin{align*}
\left\|\partial^{(i,j)}_{v}\varphi(x,v) -\partial^{(i,j)}_{v}\varphi(x,\overline{v}^{\textsc{f}})\right\| & \leq M \left\| v - \overline{v}^{\textsc{f}} \right\|^{\delta}  \leq M \left(\vartheta^{\textsc{f}}\right)^\delta \left\| v - v^* \right\|^{\delta}.
\end{align*}
Hence we get
\begin{align*}
\frac{1}{\varepsilon}\left\|R^{\textsc{f}}_{\varphi}\right\| & \leq \frac{1}{2\varepsilon} M\left(\vartheta^{\textsc{f}}\right)^\delta \mathbb{E}\left(\int_{\R^{4d}} \left\| v^* - v \right\|^{2+\delta} f(x,v) f(\hat{x},\hat{v})  \ dx \ dv \ d\hat{x} \ d\hat{v}\right) \\
& = \frac{1}{2\varepsilon}M \left(\vartheta^{\textsc{f}}\right)^\delta |\etaF|^{2+\delta} \mathbb{E}\left(\int_{\R^{4d}} \left\| \theta \xi + \Adcut + \rhoF \HFF \right\|^{2+\delta} f(x,v) f(\hat{x},\hat{v})  \ dx \ dv \ d\hat{x} \ d\hat{v}\right) \\
& = \frac{1}{2}M \left(\vartheta^{\textsc{f}}\right)^\delta |\varepsilon|^{1+\delta} \mathbb{E}\left(\int_{\R^{4d}} \left\| \theta \xi + \Adcut + \rhoF \HFF \right\|^{2+\delta} f(x,v) f(\hat{x},\hat{v})  \ dx \ dv \ d\hat{x} \ d\hat{v}\right).
\end{align*}
From the inequality
\begin{align*}
\left\| \theta \xi + \Adcut + \rhoF \HFF \right\|^{2+\delta} & \leq 2^{2+2\delta} \left( \left\| \theta \xi  \right\|^{2+\delta} + \left\| \Adcut \right\|^{2+\delta} + \left\| \rhoF \HFF \right\|^{2+\delta}\right),
\end{align*}
and since $\theta\in[0,1]$, we obtain
\begin{align*}
\left\|R^{\textsc{f}}_{\varphi}\right\| & \leq 2^{1+2\delta} M \left(\vartheta^{\textsc{f}}\right)^\delta |\varepsilon|^{1+\delta} \Bigg( \left(\rhoF\right)^2\mathbb{E}\left(\left\|\xi  \right\|^{2+\delta}\right) + \rhoF \int_{\R^{2d}} \left\| \Adcut \right\|^{2+\delta} f(x,v) \ dx \ dv \\
& \quad + \left(\rhoF\right)^{2+\delta}\int_{\R^{4d}} \left\| \HFF \right\|^{2+\delta} f(x,v) f(\hat{x},\hat{v}) \ dx \ dv \ d\hat{x} \ d\hat{v} \Bigg).
\end{align*}
Analogously for $\left\|R^{\textsc{l}}_{\varphi}\right\|$ we get
\begin{align*}
\left\|R^{\textsc{l}}_{\varphi}\right\| & \leq \frac{1}{2}M \left(\vartheta^{\textsc{l}}\right)^\delta |\varepsilon|^{1+\delta} \left(\rhoL\right)^{2+\delta}\int_{\R^{4d}} \left\| \HFL \right\|^{2+\delta} f(x,v) g(\tilde{x},\tilde{v}) \ dx \ dv \ d\tilde{x} \ d\tilde{v}.
\end{align*}
Similar computations performed on the terms $\RestF [f,f]$ and $\RestL [f,g]$ yield the inequalities
\begin{align*}
\frac{\varepsilon}{2} \left\|\RestF [f,f]\right\| & \leq \frac{\varepsilon}{2} M \vartheta^{\textsc{f}}\Bigg( \rhoF \int_{\R^{2d}} \left\| \Adcut \right\|^{2} f(x,v) \ dx \ dv \\
& \quad + \left(\rhoF\right)^2 \int_{\R^{4d}} \left\| \HFF \right\|^{2} f(x,v) f(\hat{x},\hat{v})\ dx \ dv \ d\hat{x} \ d\hat{v} \Bigg), \\
\frac{\varepsilon}{2} \left\|\RestL [f,g]\right\| & \leq \frac{\varepsilon}{2} M \vartheta^{\textsc{l}}\left(\rhoL\right)^2 \int_{\R^{4d}} \left\| \HFL \right\|^{2} f(x,v) g(\tilde{x},\tilde{v}) \ dx \ dv \ d\tilde{x} \ d\tilde{v}.
\end{align*}

Remembering that solutions of \eqref{eq:micromeso} have support uniformly bounded in time and finite constant mass (i.e., $\rhoF(t) = \rhoF(0) < +\infty$), we have proven the following result.
\begin{theorem}\label{th:BoltztoFP}
Fix $T > 0$, $\delta > 0$, $u:[0,T] \rightarrow \R^{d\NL}$, and let, for every $\varepsilon > 0$, $(f^{\varepsilon},\xl^{\varepsilon}_1, \ldots, \xl^{\varepsilon}_{\NL})$ be a $\delta$-weak solution of \eqref{eq:micromeso} corresponding to the control $u$ and the initial condition $(f^0,\xl^0_1, \ldots, \xl^0_{\NL})$, and where the quantities $\etaF, \etaL, \lamF, \lamL,$ and $\varsigma^2$ are rescaled w.r.t.\ $\varepsilon$ according to \eqref{eq:scaling}. Suppose that:
\begin{enumerate}
\item $\mathbb{E}\left(\left\| \xi  \right\|^{2+\delta}\right)$ is finite,
\item the functions $\Adcut, \HFF$ and $\HFL$ are in $L^p_{loc}$ for $p = 2, 2+\delta$.
\end{enumerate}
Then, as $\varepsilon \rightarrow 0$, the solutions $(f^{\varepsilon},\xl^{\varepsilon}_1, \ldots, \xl^{\varepsilon}_{\NL})$ converge pointwise, up to a subsequence, to $(f,\xl_1, \ldots, \xl_{\NL})$, where $f$ satisfies the Fokker-Planck-type equation \eqref{eq:FokkerPlanck} with initial datum $(f^0,\xl^0_1, \ldots, \xl^0_{\NL})$, and for every $k = 1, \ldots, \NL$, $\xl_k$ satisfies 
\begin{align*}
\left\{
\begin{array}{l}
\displaystyle
\dot{\xl}_k = \vl_k = \int_{\R^{2d}} \HLF(\xl_k,x) f(x,v) \ dx \ dv + \sum_{\ell = 1}^{\NL} \HLL(\xl_k,\xl_{\ell}) + u_k, \\
\xl_k(0) =\xl_k^0.
\end{array}
\right.
\end{align*}

\end{theorem}

\begin{remark}
The hypothesis of admissibility in Definition \ref{def:PDE} makes sure that all the integrals considered in the proof of Theorem \ref{th:BoltztoFP} exist and are finite.
\end{remark}

\begin{remark}
The choice of the functions $\Adcut, \HFF$ and $\HFL$ made in \eqref{eq:Adcut} and \eqref{eq:HFF} clearly satisfy the hypotheses of Theorem \ref{th:BoltztoFP} for any $\delta > 0$. Hence, we expect the simulation of system \eqref{eq:micromeso} with the scaling \eqref{eq:scaling} to be a good approximation of the ones of the original microscopic model \eqref{eq:micro} for sufficiently small values of $\varepsilon$.
\end{remark}

\section{Optimal control of the crowd}\label{sec:optimalcontrol}
In this section we discuss how one can (optimally) control the crowd of followers by means of ``invisible'' leaders, whose dynamics are obtained by the minimization of certain cost functionals. It is well known that in case of alignment-dominated models the invisibility of leaders is not a limitation \emph{per se}, but nothing is known about more complex models. In particular, \textit{the coupling between random walk and alignment} gives rise to interesting phenomena, cf.\ \cite{vicsek1995novel}. Note that, in order to fit real behavior (see also Section \ref{sec:theexperiment}), the ratio between random walk force and alignment force should be large enough to assure a complete exploration of the domain, and small enough to catch the herding behavior. By the way, this also narrows the choice of parameters.

\subsection{A preliminary test} 
\begin{figure}[h!]
\centering
\includegraphics[width=0.27\textwidth]{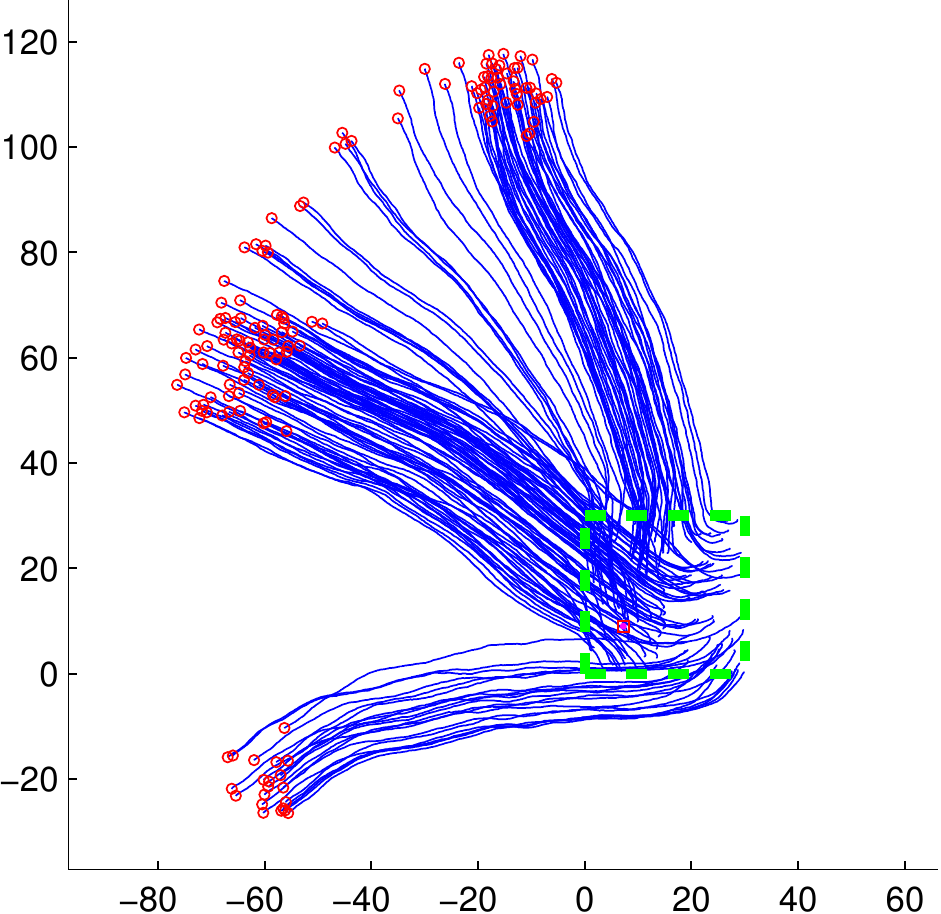}\qquad
\includegraphics[width=0.27\textwidth]{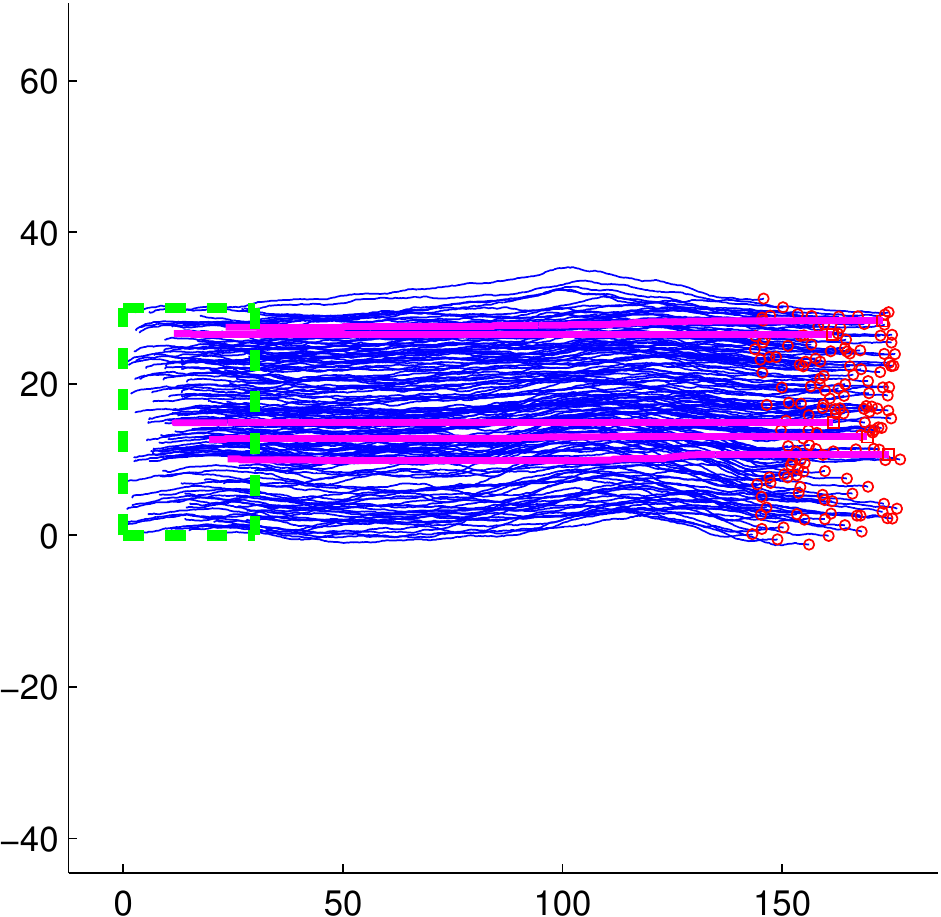}\qquad
\includegraphics[width=0.27\textwidth]{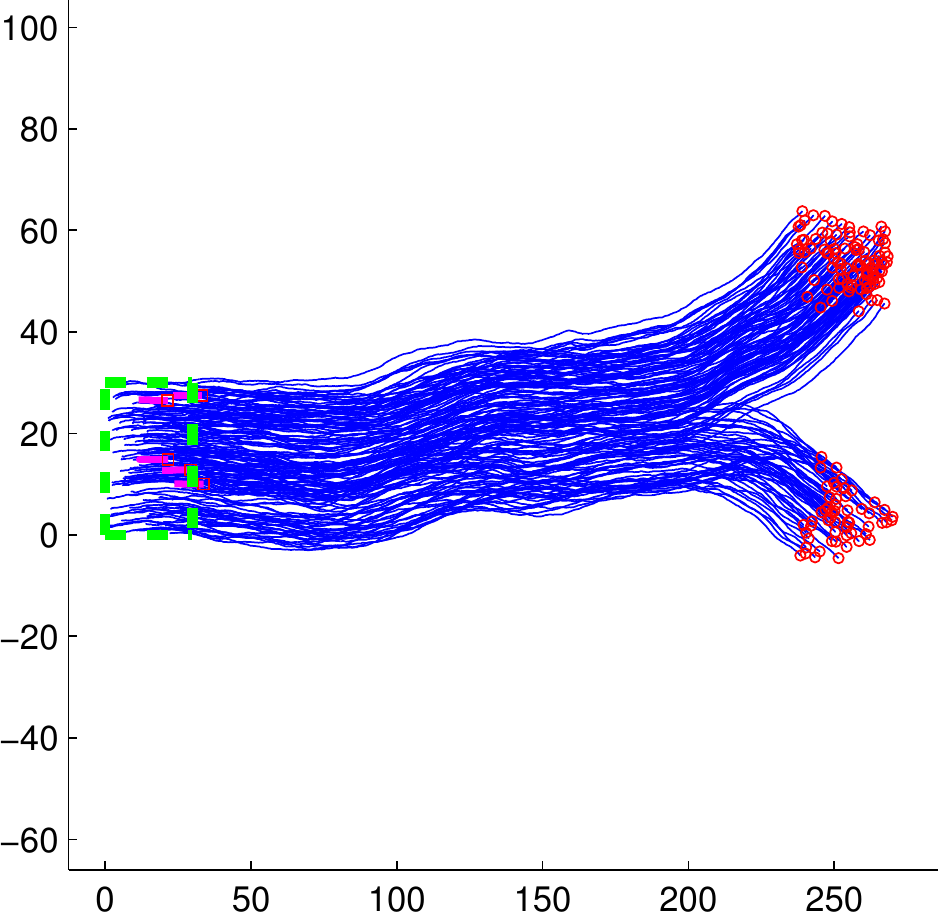}
\caption{Setting 0 (no exit is visible here). Microscopic dynamics. The crowd is initially confined in the green dashed square. Leaders' trajectories are in magenta, followers' are in blue. Final positions of followers are in red. Left: no leaders. Center: 5 leaders moving rightward all the time. Right: 5 leaders moving rightward, disappearing after a short time.}
\label{fig:S0}
\end{figure}

To fix the ideas, we plot a typical outcome of the microsimulator when the crowd is far away from the exit. At the initial time, agents are uniformly distributed in a square, see Fig.\ \ref{fig:S0}. We refer to this situation as Setting 0 (model parameters for this and following simulations are reported in Table \ref{tab:all_parameters}). With no leaders, the group splits in several subgroups, each of them having a common direction of motion (local consensus is reached in very short time), see Fig.\ \ref{fig:S0}-left. With 5 leaders moving rightward, the crowd reaches immediately a consensus and align with leaders, see Fig.\ \ref{fig:S0}-center. If, for some reason, leaders cease to be influential, the crowd tend to split again, see Fig.\ \ref{fig:S0}-right. 
\begin{table}[!h]
\caption{Model parameters.}
\label{tab:all_parameters}
\begin{center}
\begin{tabular}{|c|c|c|c|c|c|c|c|c|c|c|c|c|}
\hline
Setting & $\NL$ & $\NF$ & $\mathcal N$ & $\Cr$ & $\Crl$ & $\Ca$ & $\Crand$ & $\Ct$ & $\Cchar$ & $\barS^2$ & $r=\zeta$ & $\gamma$\\
\hline\hline
0 & 0-5 & 150   & 10 & 2 & 1.5 & 2 & 0.25 & -- & 1 & 0.5 & 0.4 & 1 \\
\hline
1 & 0-3 & 150   & 10 & 2 & 1.5 & 3 & 0.2  & 1  & 1 & 0.5 & 0.4 & 1 \\
\hline
2 & 0-2 & 100   & 10 & 2 & 1.5 & 3 & 0.2  & 1  & 1 & 0.5 & 0.4 & 1 \\
\hline
3 & 0-6 & 30-70 & 10 & 1 & 0    & -- & --  & 1  & 1 & 0.5 & 0.5 & 1 \\
\hline
\end{tabular}
\end{center}
\end{table}

This makes clear how difficult is to control the crowd in this scenario, since leaders have to fight continuously against the natural tendency of the crowd to split in subgroups and move randomly, even after that a consensus is reached, cf.\ \cite{han2013PLOSONE}.

A natural question is about the minimum number of leaders needed to lead the whole crowd to consensus, i.e.\ to align all the agents along a desired direction. Numerical simulations suggests that this number strongly depends on the initial conditions. Table \ref{tab:consensus} gives a rough idea in the case of multiple runs with random initial conditions. 
\begin{table}[!h]
\caption{Number of leaders versus percentage of runs where consensus is reached.}
\label{tab:consensus}
\begin{center}
\begin{tabular}{|c|c|c|c|c|c|}
\hline
\#leaders & 1 & 2 & 3 & 4 & 5\\
\hline
consensus reached & 0\% & 1\% & 27\% & 58\% & 100\% \\
\hline
\end{tabular}
\end{center}
\end{table}

\begin{figure}[h!]
\centering
\includegraphics[width=0.27\textwidth]{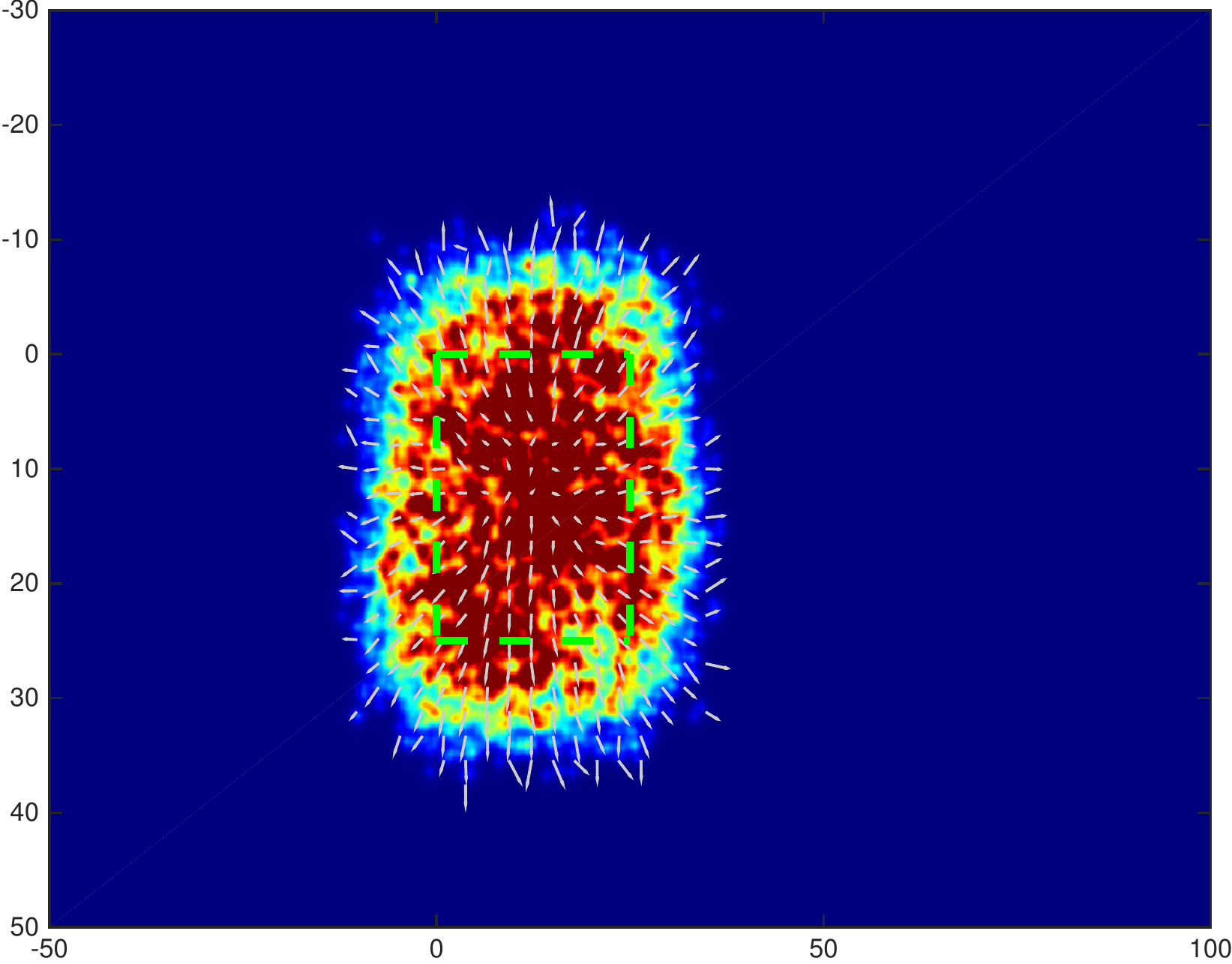}\qquad
\includegraphics[width=0.27\textwidth]{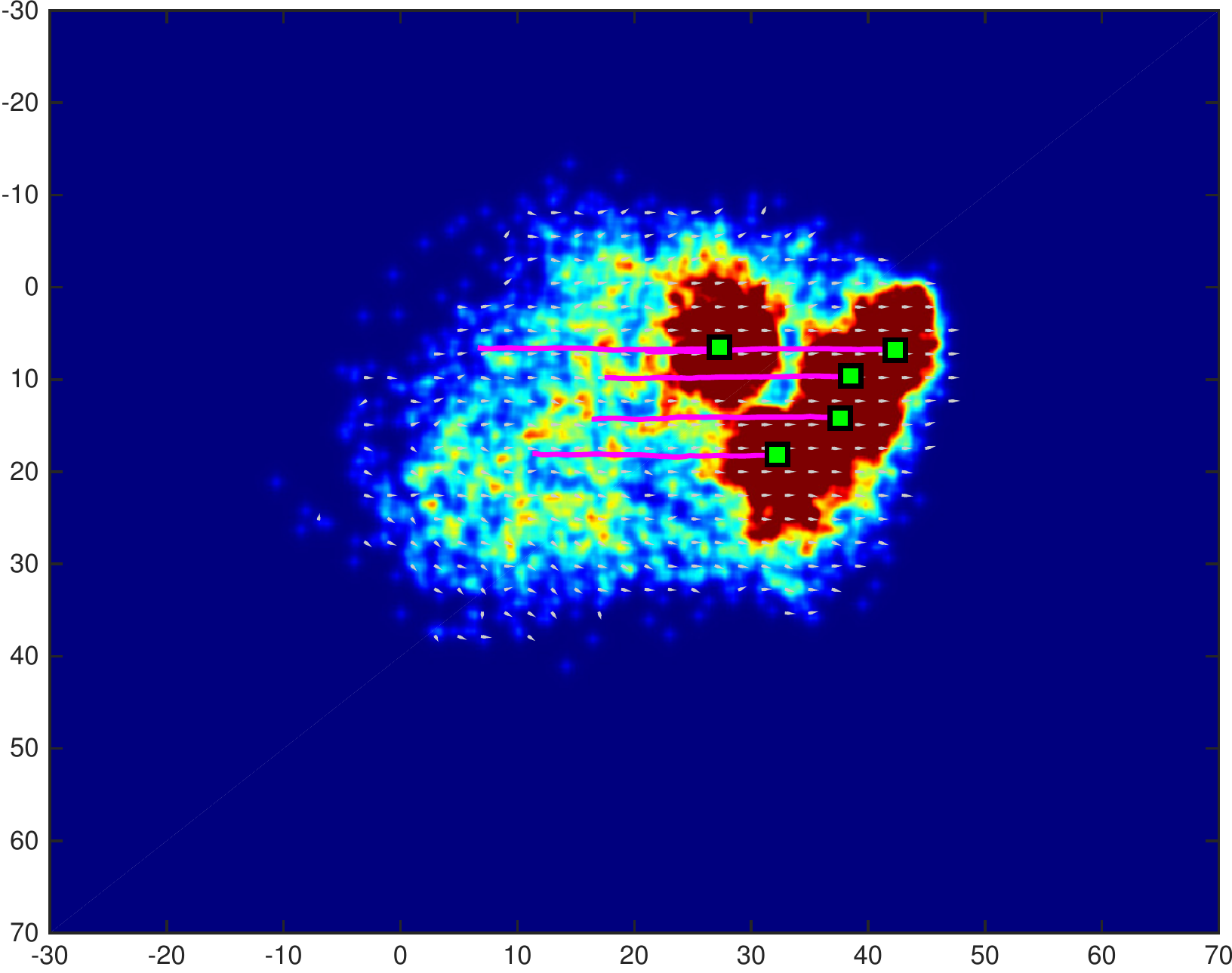}\qquad
\includegraphics[width=0.27\textwidth]{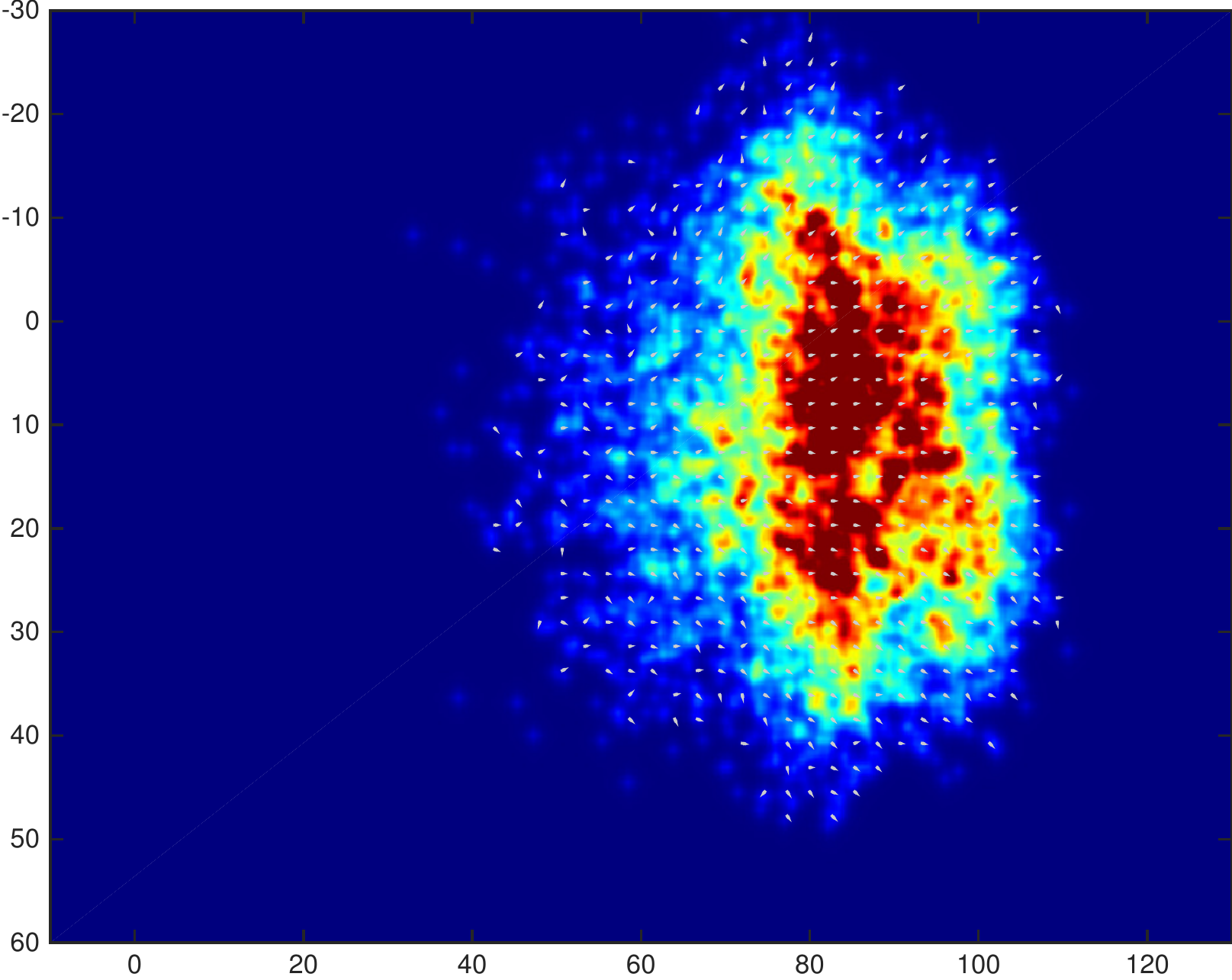}
\caption{Setting 0 (no exit is visible here). Mesoscopic dynamics. The crowd density is initially confined in the green dashed square shown in the left figure. Left: no leaders. Center: 5 leaders moving rightward all the time. Right: 5 leaders moving rightward, disappearing after a short time.}
\label{fig:S0_kin}
\end{figure}

At mesoscopic level we consider the control of the Boltzmann-type dynamics \eqref{eq:micromeso}, where we account for small values of the scaling parameter $\varepsilon$, in order to be sufficiently close to the Fokker-Planck model  \eqref{eq:FokkerPlanckStrong}. 
In this case, the control through few microscopic leaders has to face the tendencies of the continuous density to spread around the domain and to locally align with the surrounding mass.
Moreover, their action is weakened by the type of interaction considered.
Figure \ref{fig:S0_kin} shows the results corresponding to the kinetic density in the Setting 0. Unlike the microscopic case, no splitting occurs where there is no action of the leaders. Instead, the mass smears out around the domain, due to the diffusion term.
\vskip +0.1cm

\subsection{Setting the optimization problem} 
The functional to be minimized can be chosen in several ways. The effectiveness mostly depends on the optimization method which is used afterwards. The most natural functional is the \textit{evacuation time},
\begin{equation}\label{defJ_evactime}
\min\{t>0~|~x_i(t)\notin\Omega \quad \forall i=1,\ldots,\NF\},
\end{equation}
subject to (\ref{eq:micro}) or (\ref{eq:micromeso}) and with $u(\cdot)\in U_{\textup{adm}}$, where $U_{\textup{adm}}$ is the set of admissible controls (including for instance box constraints to avoid excessive velocities).

Another cost functional, more affordable by standard methods \cite{albi1401.7798, borzi2015M3ASa}, is
\begin{equation}\label{defJ}
\begin{split}
l(\mathbf{\xf},\mathbf{\xl},u)=\CclosenessFF\sum_{i=1}^{\NF}\|x_i-\target\|^2+
\CclosenessLF\sum_{i=1}^{\NF}\sum_{k=1}^{\NL}\|x_i-y_k\|^2+
\Ccontrol\sum_{k=1}^{\NL}\|u_k\|^2,
\end{split}
\end{equation}
for some positive constants $\CclosenessFF$, $\CclosenessLF$, and $\Ccontrol$. The first term promotes the fact that followers have to reach the exit while the second forces 
leaders to keep contact with the crowd. 
The last term penalizes excessive velocities. This minimization is performed at every instant (instantaneous control), or along a fixed time frame $[t_i,t_f]$
\begin{equation}\label{eq:finhor}
\underset{u(\cdot)\in U_{\textup{adm}}}{\min} \int\limits_{t_i}^{t_f}l(\mathbf{\xf}(t),\mathbf{\xl}(t),u(t))\,dt,\qquad \text{subject to (\ref{eq:micro}) or (\ref{eq:micromeso})}.
\end{equation}

With regards to the mesoscopic scale, both functionals \eqref{defJ_evactime} and \eqref{defJ} can be considered, however at the continuous level the presence of few invisible leaders does not assure that the whole mass of followers is evacuated. The major difficulty to reach a complete evacuation of the continuous density is mainly due to the presence of the diffusion term and to the invisible interaction with respect to the leaders. Therefore a more appropriate functional is given by the mass evacuated at the final time $T$,
\begin{equation}\label{defJ_evacmass}
\underset{u(\cdot)\in U_{\textup{adm}}}{\min}\left\{\rhoF(T)=\int_{\mathbb{R}^{d}}\int_{\Omega}f(T,x,v)\,dx\,dv\right\},\qquad \text{subject to (\ref{eq:micromeso})}.
\end{equation}

\begin{remark}
We will not require that leaders themselves reach the exit, but only the followers.
\end{remark}

\medskip

\subsubsection{Model predictive control}  A computationally efficient way to address the optimal control problem \eqref{eq:finhor} is by means of a relaxed approach known as model predictive control (MPC) \cite{MayneRawlingsRao2000aa}. We consider a sampling of the dynamics \eqref{eq:micro} at every time interval $\Delta t$, and the following minimization problem over $\Nmpc$ time steps, starting from the current time step $\bar n$,
\begin{equation}
\underset{u(\cdot)\in U_{\textup{adm}}}{\min} 
\sum_{n=\bar n}^ {\bar n+\Nmpc-1} 
l(\mathbf{x}(n\Delta t),\mathbf{y}(n\Delta t),u(n\Delta t))
\end{equation}
generating an optimal sequence of controls $\{u(\bar n\Dt),\ldots,u((\bar n+\Nmpc-1)\Dt)\}$, from which only the first term is taken to evolve the dynamics for a time $\Dt$, to recast the minimization problem over an updated time frame $\bar n\leftarrow\bar n+1$. Note that for $\Nmpc=2$, the MPC approach recovers an instantaneous controller, whereas for $\Nmpc=(t_f-t_i)/\Dt$ it solves the full time frame problem $\eqref{eq:finhor}$. Such flexibility is complemented with a robust behavior, as the optimization is re-initialized every time step, allowing to address perturbations along the optimal trajectory.

\subsubsection{Modified compass search}
When the cost functional is highly irregular and the search of local minima is particularly difficult, it could be convenient to move towards random methods as compass search (see \cite{audet2014MPC} and references therein), genetic algorithms, or particle swarm optimization. In the following we describe a compass search method that works surprisingly well for our problem.

First of all, we consider only piecewise constant trajectories, introducing suitable \emph{switching times} for the leaders' controls. More precisely, we assume that leaders move at constant velocity for a given fixed time interval and when the switching time is reached, a new velocity vector is chosen. Therefore, the control variables are the velocities at the switching times for each leader. Note that controlling directly the velocities rather than the acceleration makes the optimization problem much simpler because minimal control variations have an immediate impact on the dynamics.

Starting from an initial guess, at each iteration the optimization algorithm modifies the current best control strategy found so far by means of small random variations of the current values. Then, the cost functional is evaluated. If the variation is advantageous (the cost decreases), the variation is kept, otherwise it is discarded. The method stops when the strategy cannot be improved further. As initial guess for the strategy, leaders move at constant speed along the direction joining the target with their initial position.

%
%
\section{Numerical tests}\label{sec:simulations}
In what follows we present some numerical tests to validate our modeling framework at the microscopic and mesoscopic level. 

The microscopic model \eqref{eq:micro} is discretized by means of the explicit Euler method with a time step $\Delta t=0.1$.
The evolution of the kinetic density in \eqref{eq:FokkerPlanckStrong} is approximated by means of binary interaction algorithms, which approximates the Boltzmann dynamics  \eqref{eq:strongBoltz} with a meshless Monte-Carlo method for small values of the parameter $\varepsilon$, as presented in \cite{albi2013MMS}.  We choose $\varepsilon=0.02$, $\Delta t = 0.01$ and a sample of $N_s=10000$ particles to reconstruct the kinetic density. This type of approach is inspired by numerical methods for plasma physics and it allows to solve the interaction dynamics with a reduced computational cost compared with mesh-based methods, and an accuracy of $O(N_s^{-1/2})$. For further details on this class of binary interaction algorithms see \cite{albi2013MMS, PT:13}.

Concerning optimization, in the microscopic case we adopt either the compass search with functional \eqref{defJ_evactime} or MPC with functional \eqref{eq:finhor}. In the mesoscopic case we adopt the compass search with functional \eqref{defJ_evacmass}. 

We consider three settings for pedestrians, without and with obstacles, hereafter referred to Setting 1, 2, and 3, respectively. In Setting 1 we set the compass search switching times every 20 time steps, and in Setting 2 every 50, having fixed the maximal random variation to 1 for each component of the velocity. In Setting 1, the inner optimization block of the MPC procedure is performed via a direct formulation, by means of the \textsl{fmincon} routine in MatLab, which solves the optimization problem via an SQP method.

\subsection{Setting 1}
\begin{figure}[h!]
\begin{center}
\includegraphics[width=0.41\textwidth]{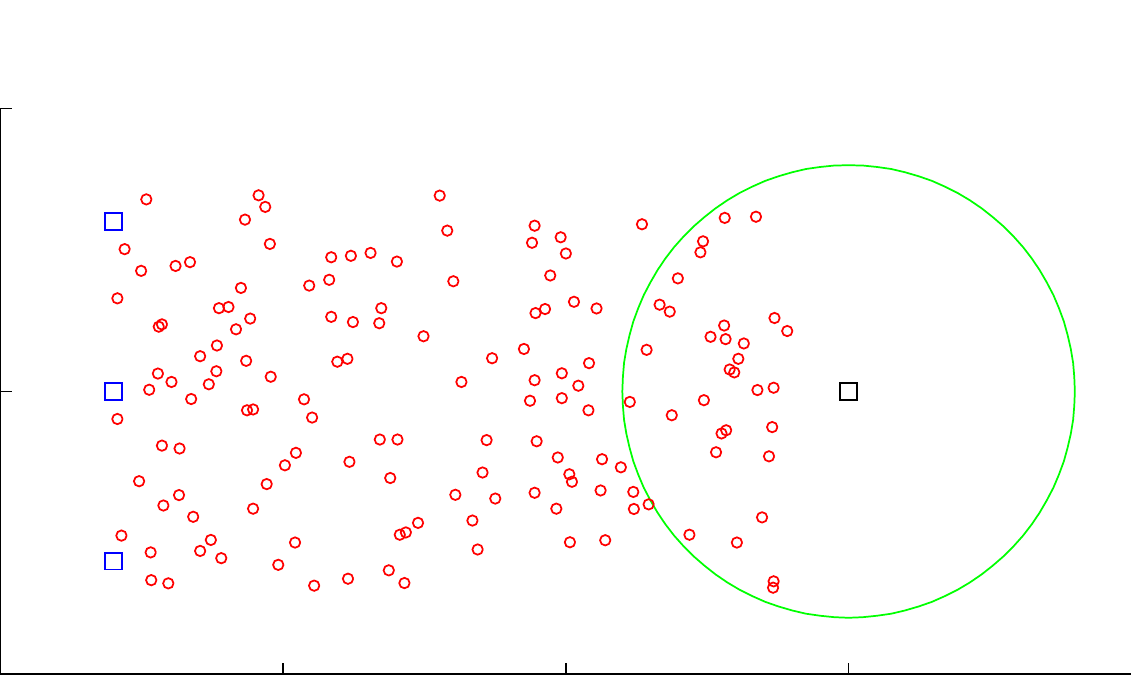}\qquad
\includegraphics[width=0.41\textwidth]{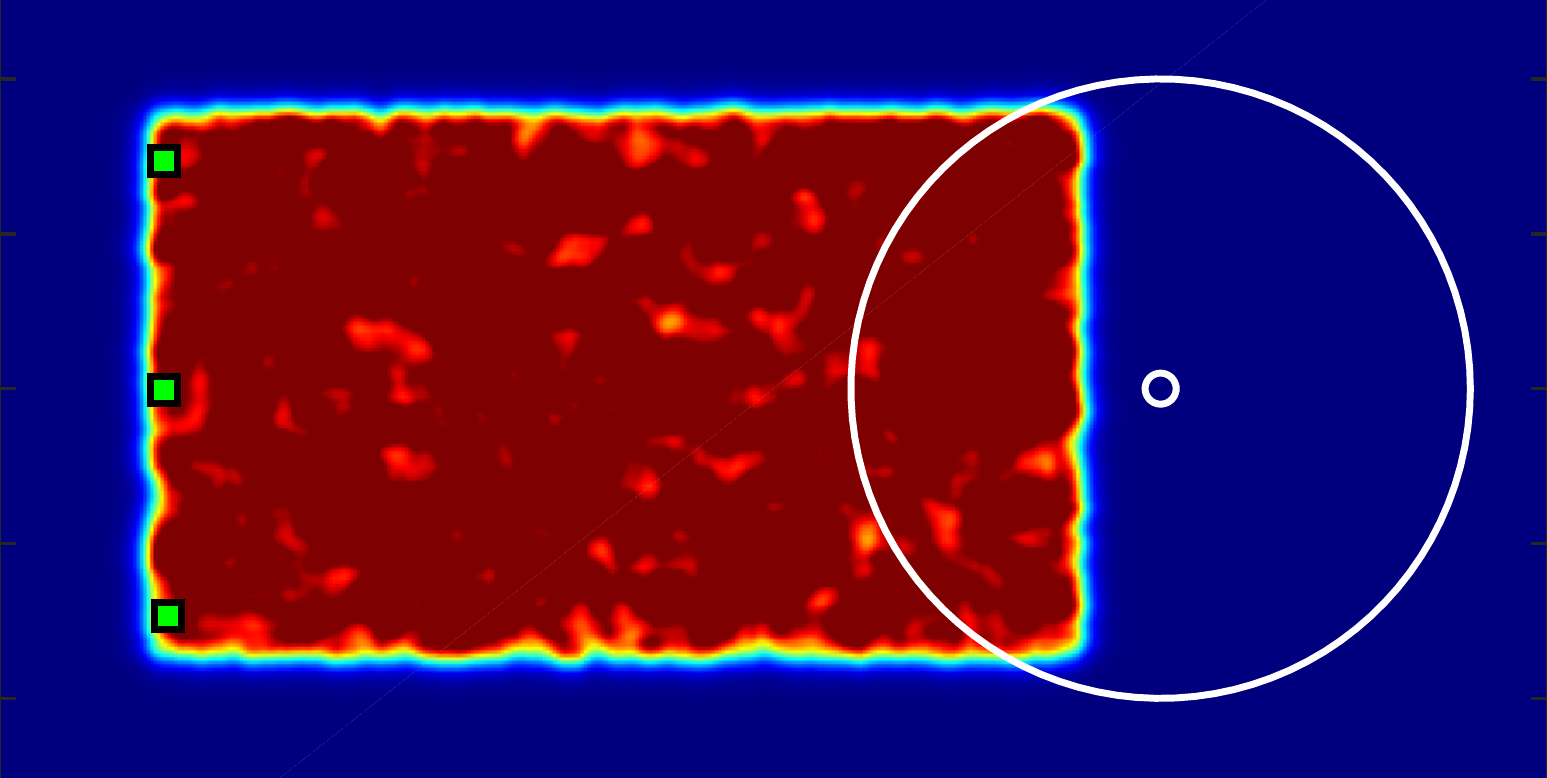}
\end{center}
\caption{Setting 1. Left: initial positions of followers (circles) and leaders (squares). Right: uniform density of followers and the microscopic leaders (squares).}
\label{fig:S1_initial_condition}
\end{figure}

To begin with, let us consider the case of a large room with no obstacles, see Fig.\ \ref{fig:S1_initial_condition}. 
The exit is a point located at $E=(30,10)$ which can be reached from any direction. We set $\Sigma=\{x\in\R^2:|x-E|<4\}$. This simple setting helps elucidating the role and the interplay of the different terms of our model. Followers are initially randomly distributed in the domain $[17,29]\times[6.5,13.5]$ with velocity $(0,0)$. Leaders, if present, are located to the left of the crowd. 
Parameters are reported in Table \ref{tab:all_parameters}.

\subsubsection{Microscopic model}
\begin{figure}[b!]
\centering
\includegraphics[width=0.3\textwidth]{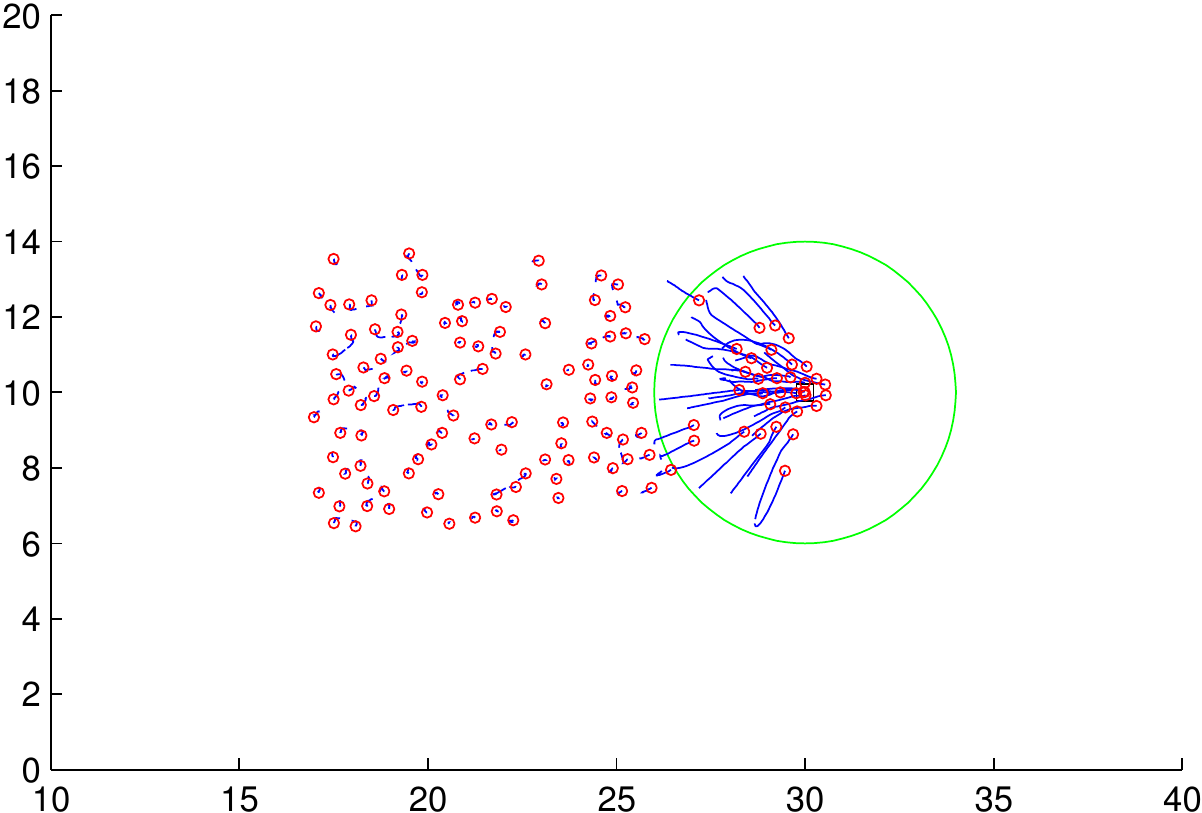}\quad
\includegraphics[width=0.3\textwidth]{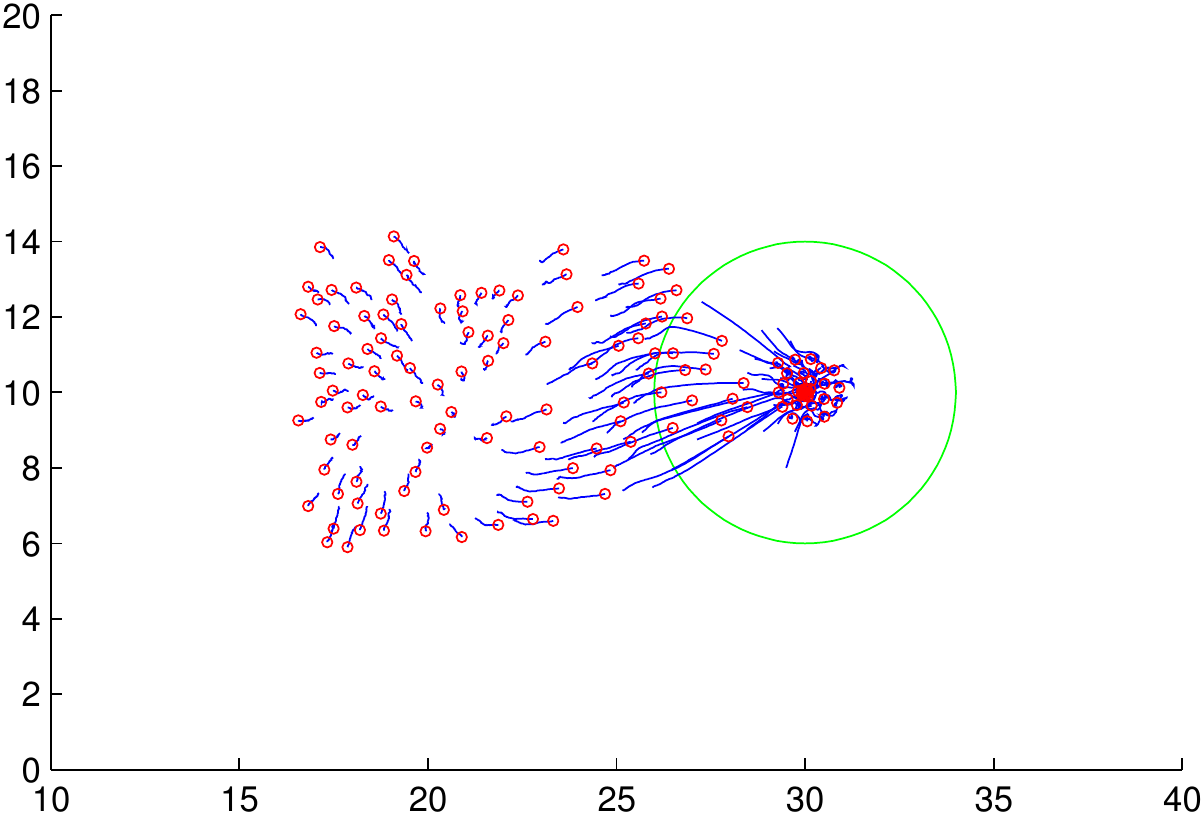}\quad
\includegraphics[width=0.3\textwidth]{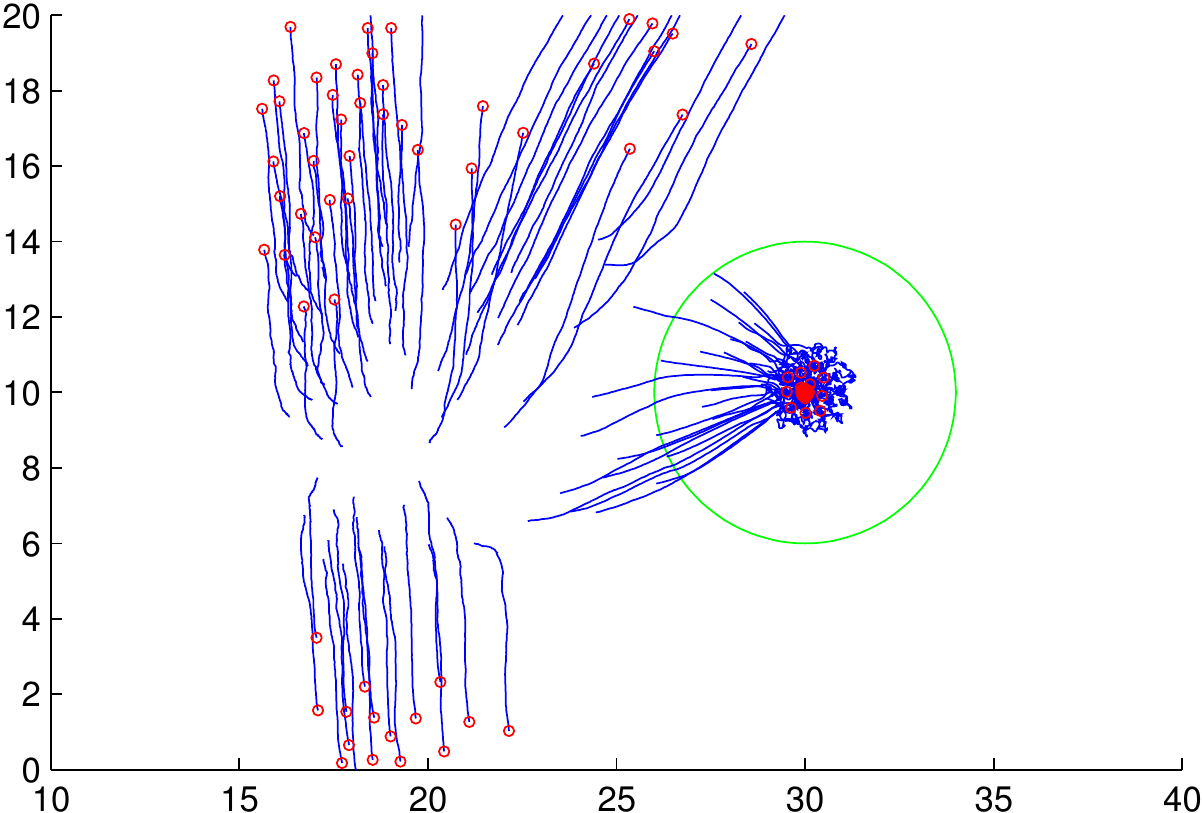}\\ [5mm]
\includegraphics[width=0.3\textwidth]{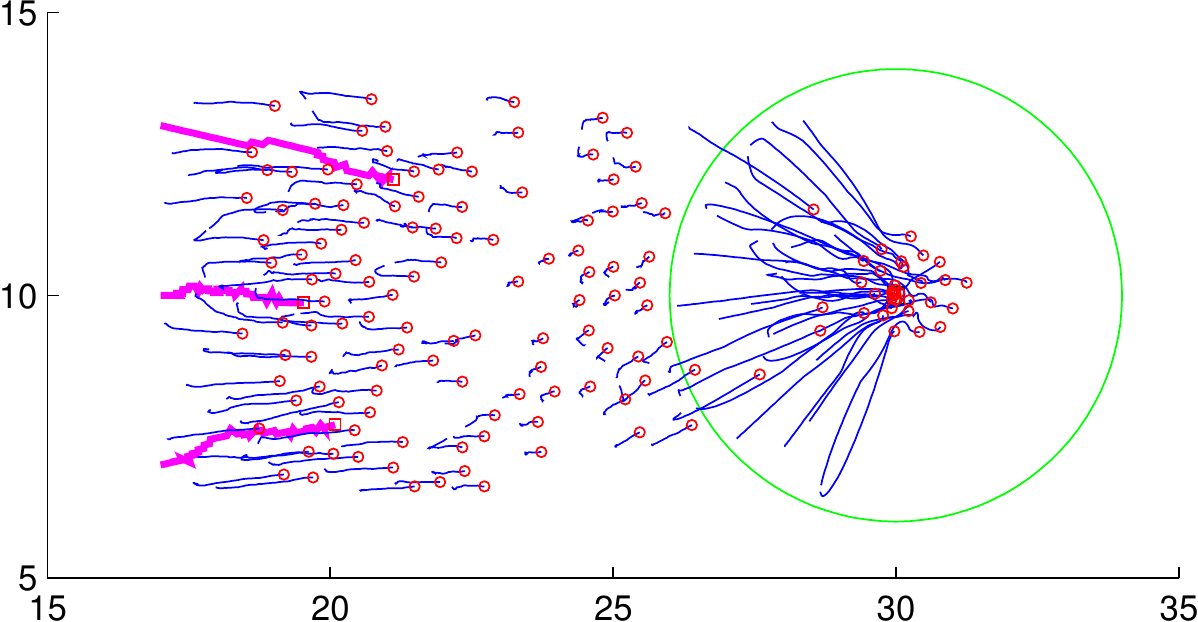}\quad
\includegraphics[width=0.3\textwidth]{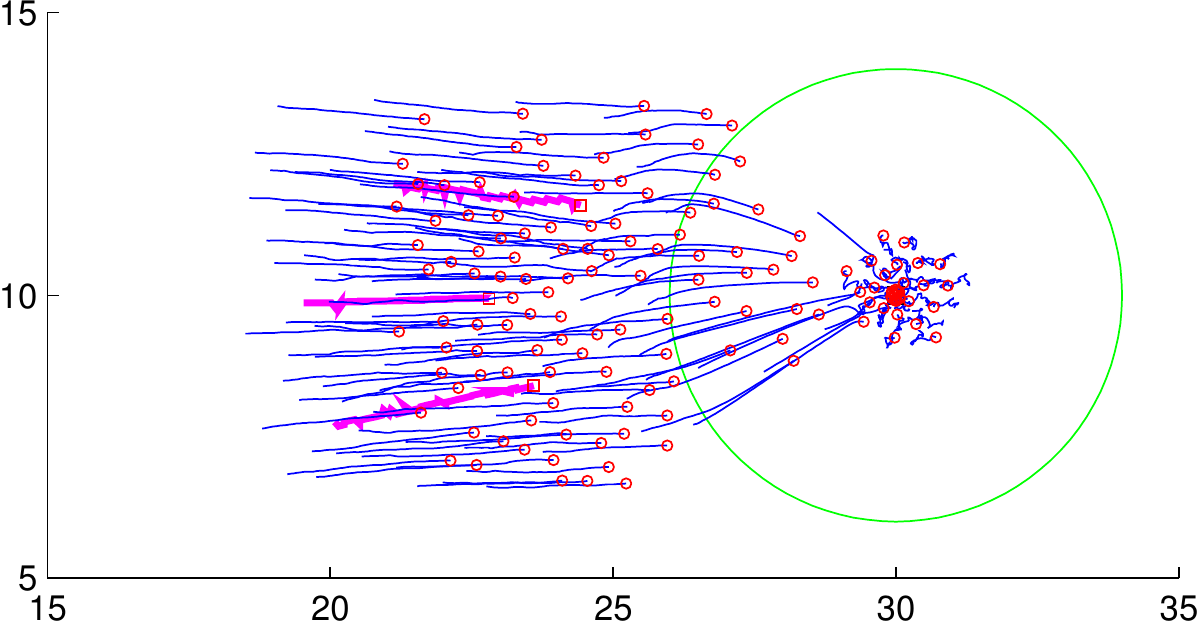}\quad
\includegraphics[width=0.3\textwidth]{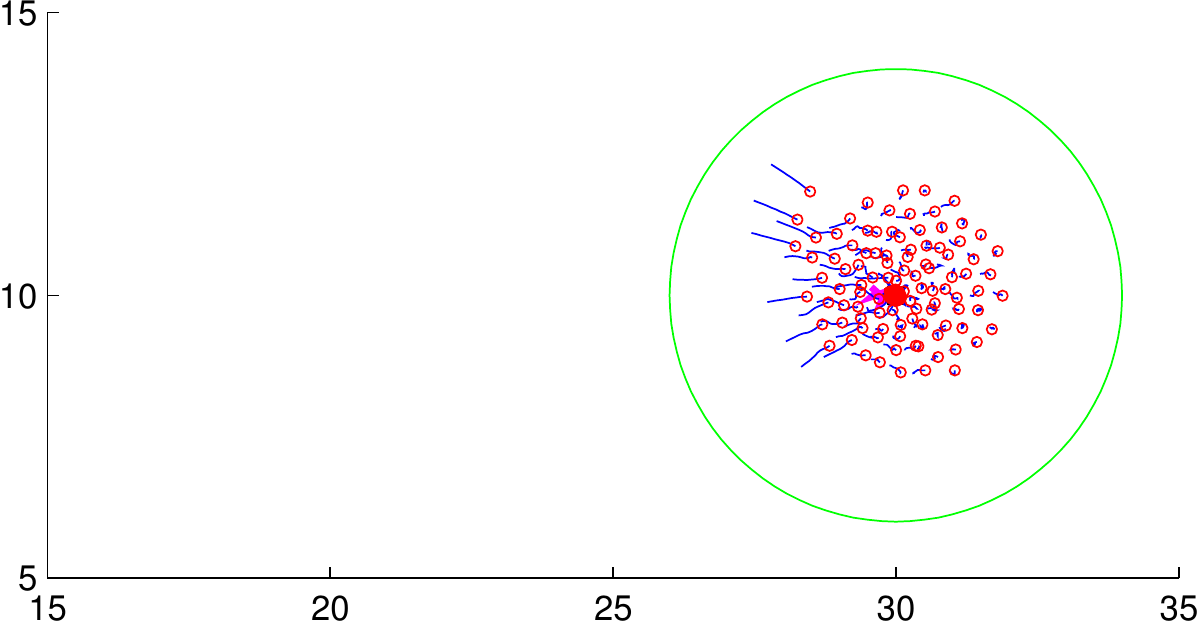}\\ [4mm]
\includegraphics[width=0.3\textwidth]{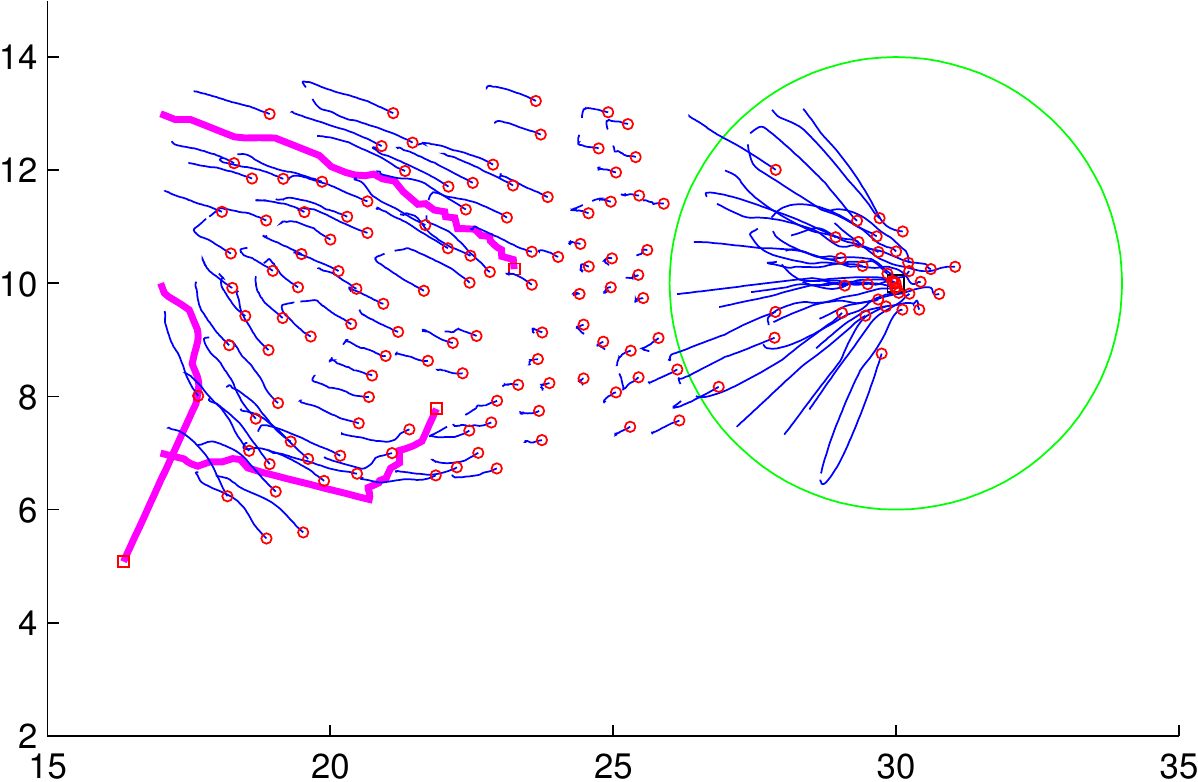}\quad
\includegraphics[width=0.3\textwidth]{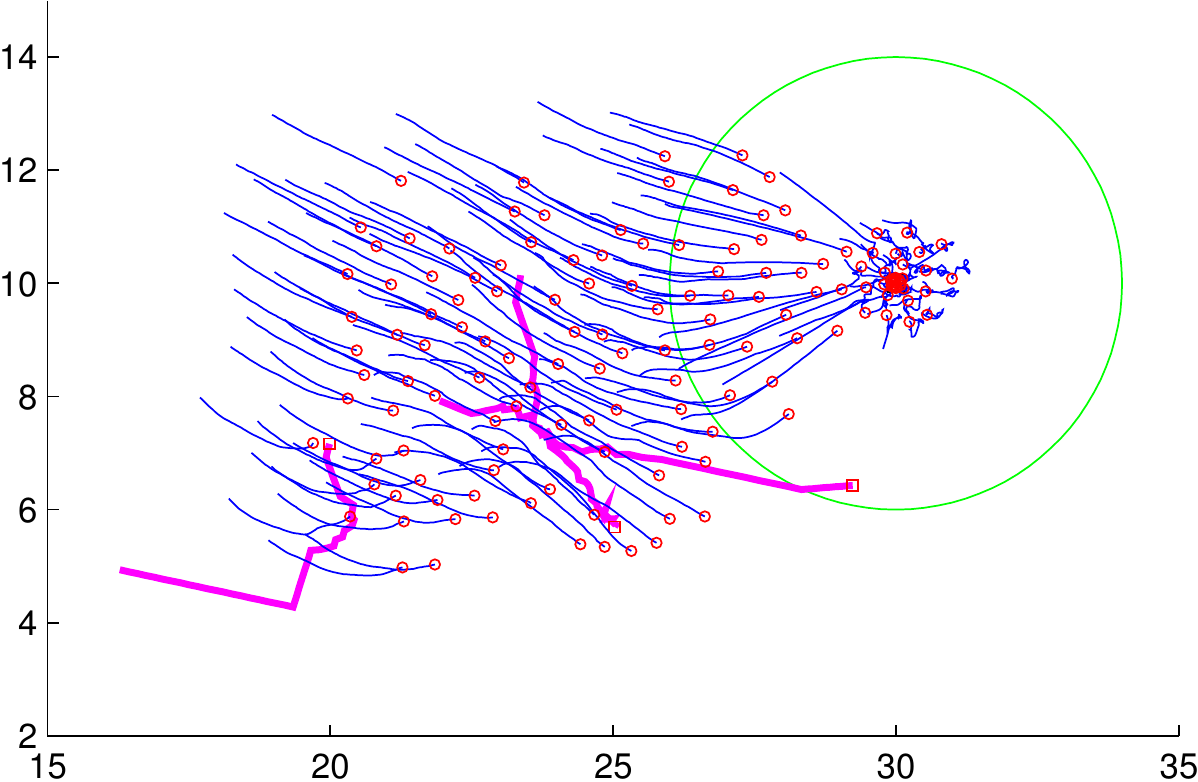}\quad
\includegraphics[width=0.3\textwidth]{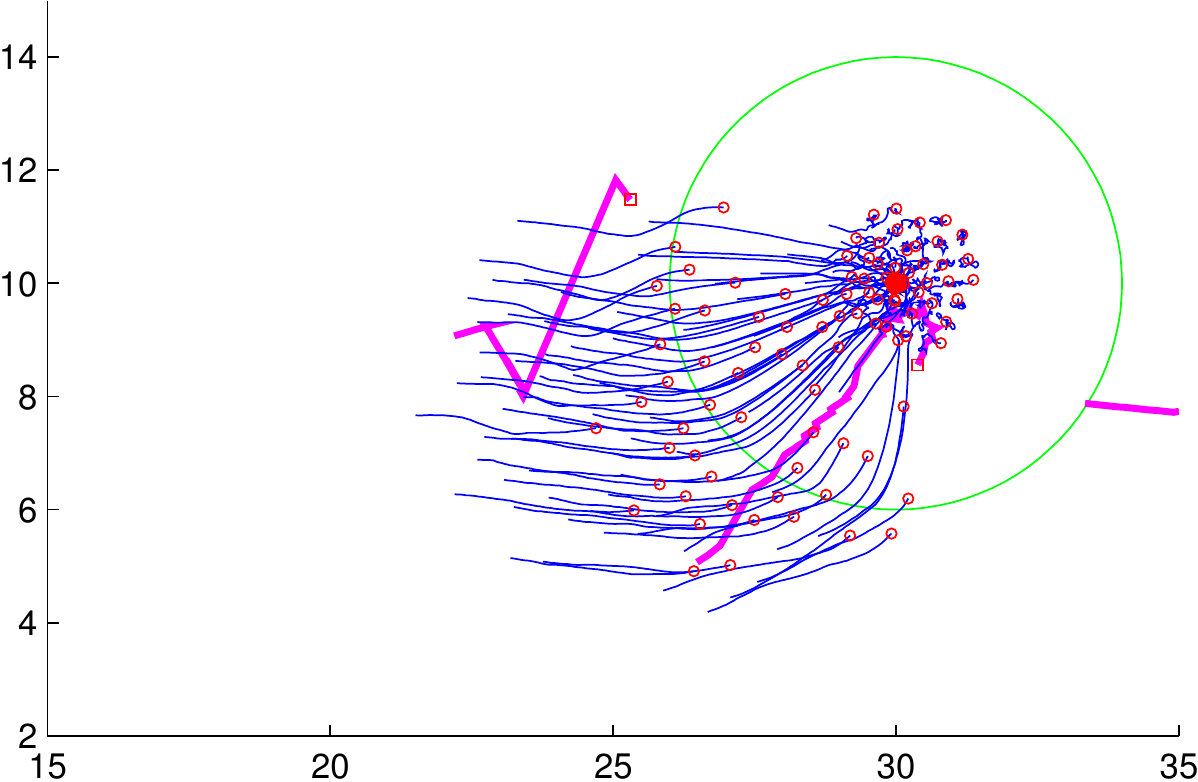}
\caption{Setting 1. Microscopic dynamics. First row: no leaders. Second row: three leaders, go-to-target strategy. Third row: three leaders, optimal strategy (compass search).}
\label{fig:S1_micro}
\end{figure}

Figure \ref{fig:S1_micro}(first row) shows the evolution of the agents computed by the microscopic model, without leaders. Followers having a direct view of the exit immediately point towards it, and some group mates close to them follow thanks to the alignment force. On the contrary, farthest people split in several but cohesive groups with random direction and never reach the exit. 

Figure \ref{fig:S1_micro}(second row) shows the evolution of the agents with three leaders. The leaders' strategy is defined manually. More precisely, at any time the control is equal to the unit vector pointing towards the exit from the current position. Hereafter we refer to this strategy as ``go-to-target''. Note that the final leaders' trajectories are not straight lines because of the additional repulsion force. As it can be seen, the crowd behavior changes completely since, this time, the whole crowd reaches the exit. However followers form a heavy congestion around the exit. The shape of the congestion is circular: this is in line with the results of other social force models as well as physical observation, which report the formation of an ``arch'' near the exits. The arch is correctly substituted here by a full circle due to the absence of walls. Note that the congestion notably delays the evacuation. This suggests that the strategy of the leaders is not optimal and can be improved by an optimization method.

Fig.\ \ref{fig:S1_micro}(third row) shows the evolution of the agents with three leaders and the optimal strategy obtained by the compass search algorithm. Surprisingly enough, the optimizator prescribes that leaders \emph{divert} some pedestrians from the right direction, so as not to steer the whole crowd to the exit at the same time. In this way congestion is avoided and pedestrian flow through the exit is increased.

In this test we have also run the MPC optimization, including a box constraint $u_k(t)\in[-1,1]$. 
We choose $\CclosenessFF=1$, $\CclosenessLF=10^{-5}$, and $\Ccontrol=10^{-5}$. 
MPC results are consistent in the sense that for $\Nmpc=2$, the algorithm recovers a controlled behavior similar to the application of the instantaneous controller (or go-to-target strategy). Increasing the time frame up to $\Nmpc=6$ improves both congestion and evacuation times, but results still remain non competitive if compared to the whole time frame optimization performed with a compass search. 

\begin{figure}[h!]
\centering
\ \includegraphics[scale=0.4]{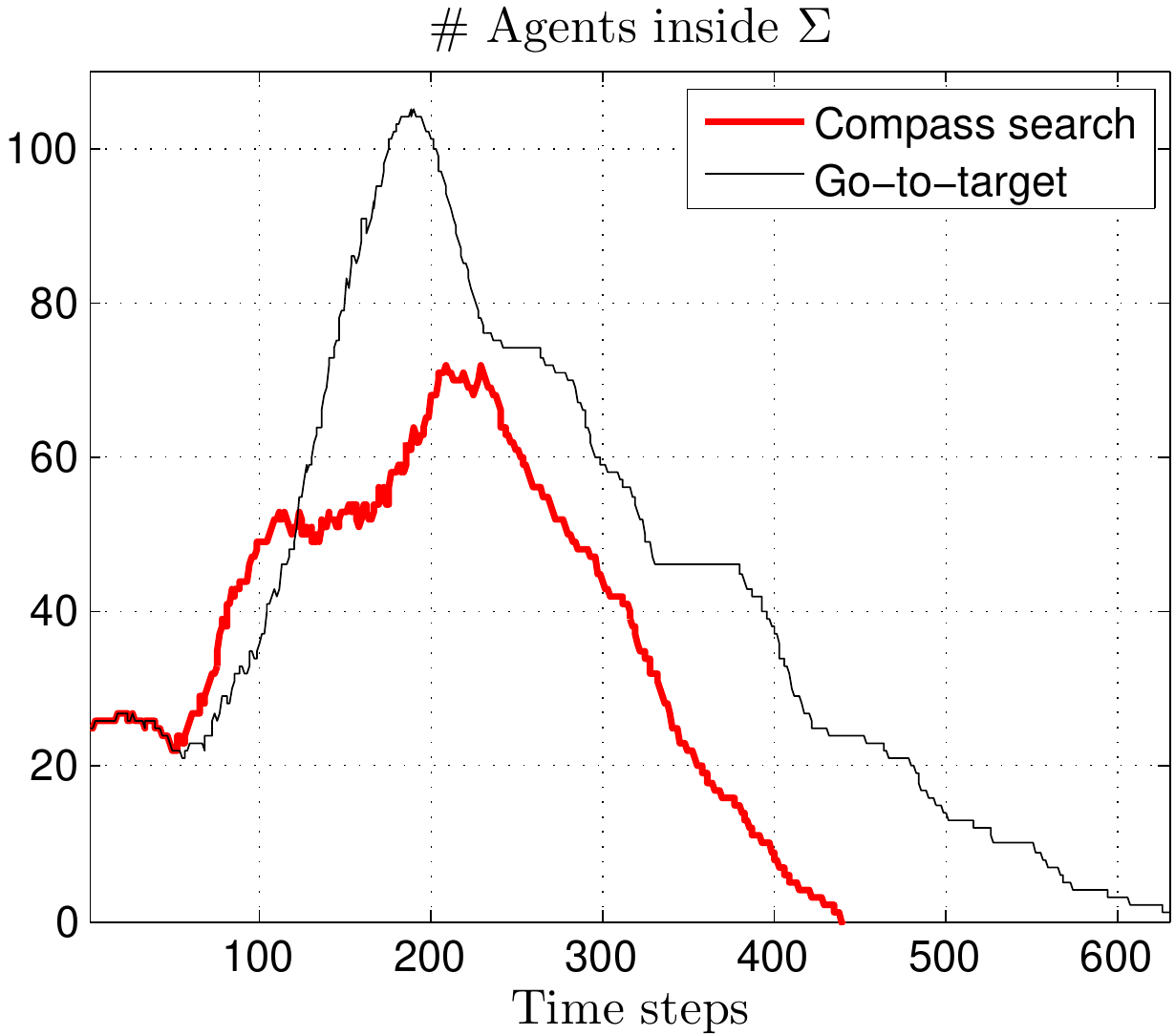}\qquad \ \
\includegraphics[scale=0.4]{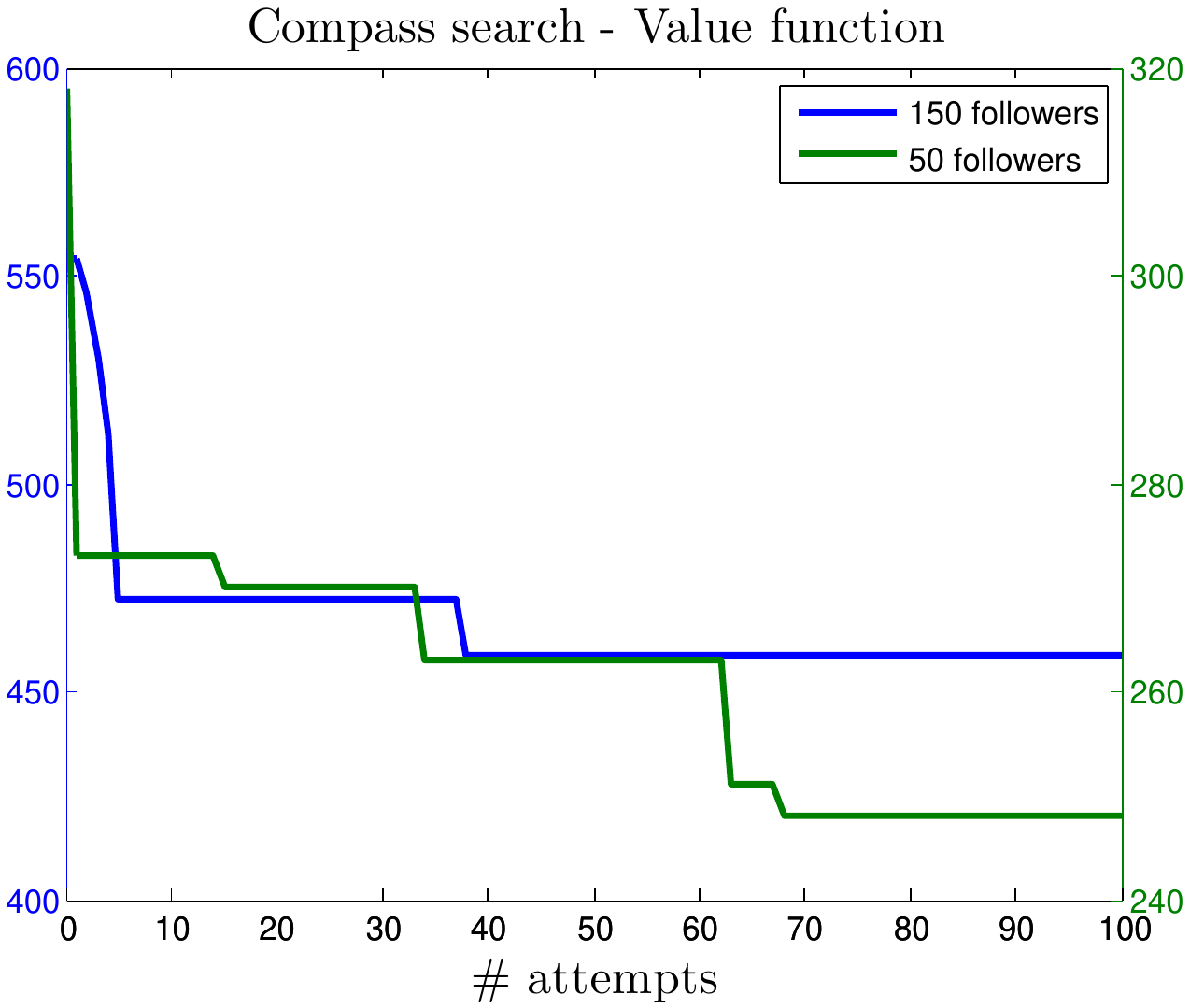}\vskip0.2cm
\includegraphics[scale=0.39]{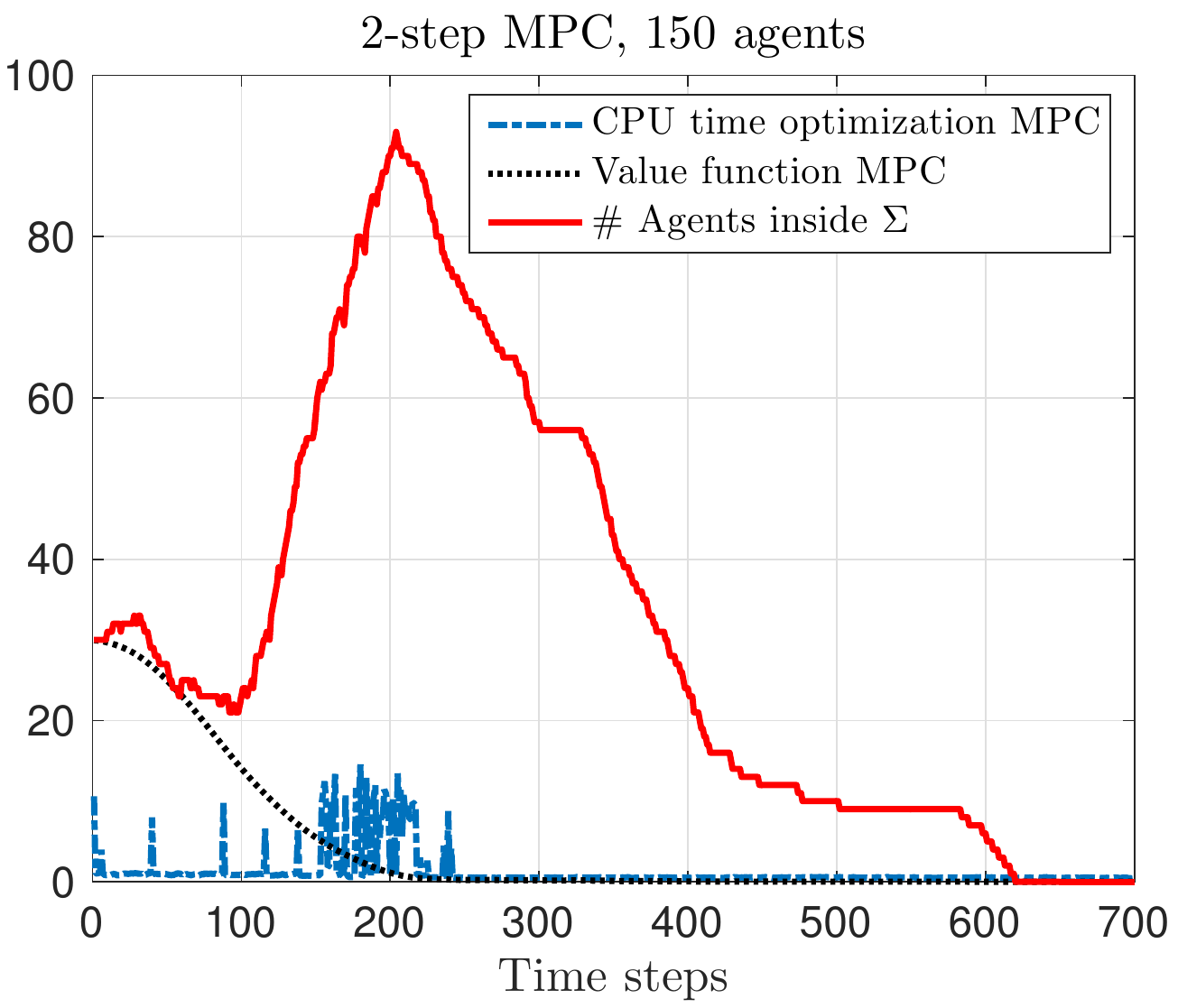}\qquad
\includegraphics[scale=0.39]{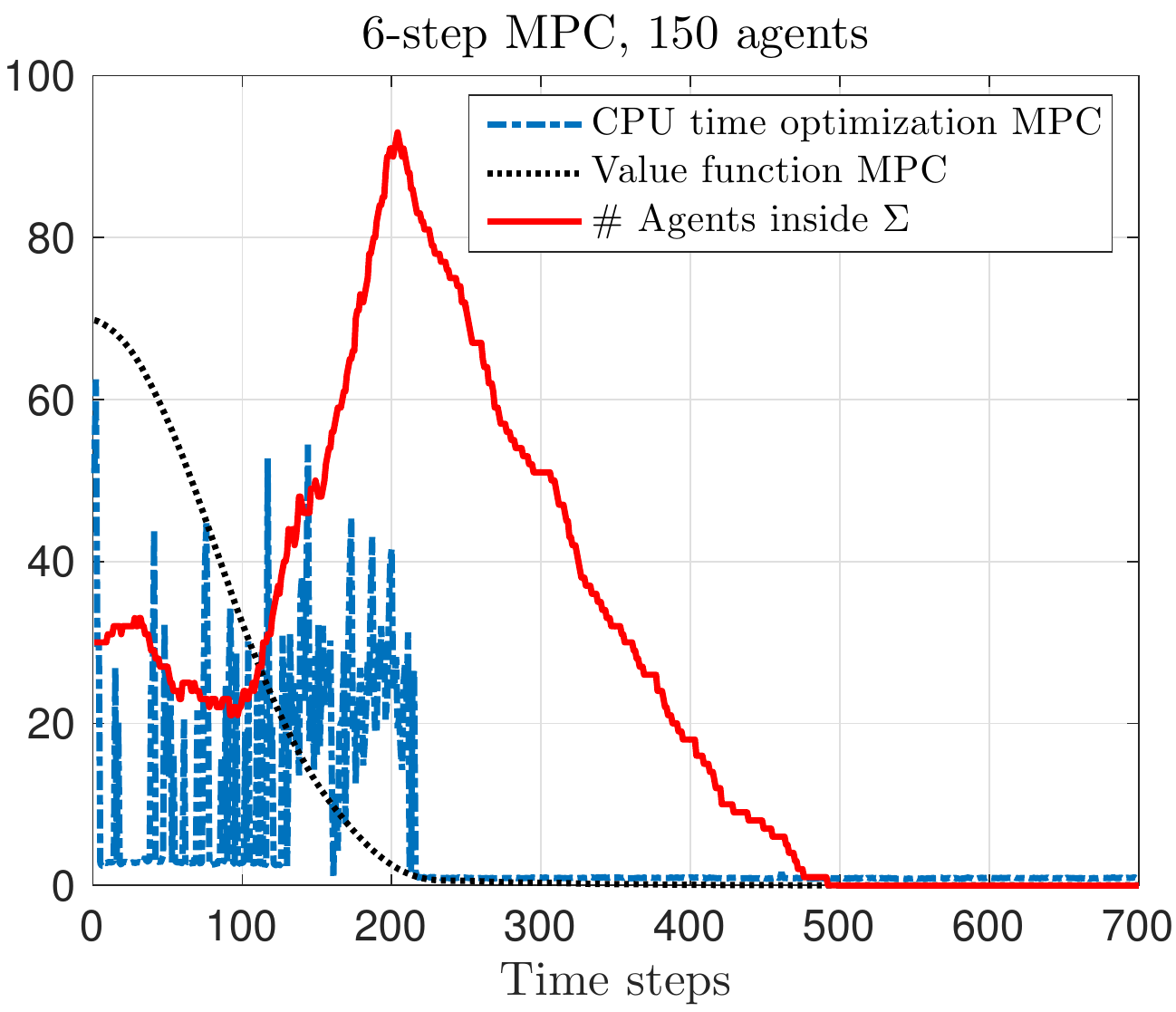}
\caption{Setting 1. Optimization of the microscopic dynamics. Top-left: occupancy of the exit's visibility zone $\Sigma$ as a function of time for optimal strategy (compass search) and go-to-target strategy. Top-right: decrease of the value function \eqref{defJ_evactime} as a function of the iterations of the compass search (for 50 and 150 followers). Bottom: MPC optimization. occupancy of the exit's visibility zone $\Sigma$ as a function of time, CPU time of the optimization call embedded in the MPC solver, and the evolution of the corresponding value (2-step and 6-step).}
\label{fig:S1_hystogram}
\end{figure}

In Fig.\ \ref{fig:S1_hystogram} we compare the occupancy of the exit's visibility zone as a function of time for go-to-target strategy and optimal strategies (compass search, 2-step, and 6-step MPC). We also show the decrease of the value function as a function of attempts (compass search) and time (MPC). Evacuation times are compared in Table \ref{tab:S1_evactimes}. It can be seen that only the long-term optimization strategies are efficient, being able to moderate congestion and clogging around the exit.
\begin{table}[!h]
\caption{Setting 1. Evacuation times (time steps). CS=compass search, IG=initial guess.}
\label{tab:S1_evactimes}
\begin{center}
\begin{tabular}{|c|c|c|c|c|c|}
\hline
 & no leaders & go-to-target & 2-MPC & 6-MPC & CS (IG)\\
\hline
$\NF=50$ & 335 & 297 & 342 & 278 & 248 (318)\\ \hline
$\NF=150$ & $\infty$ & 629 & 619 & 491 & 459 (554)\\
\hline
\end{tabular}
\end{center}
\end{table}

This suggests a quite unethical but effective evacuation procedure, namely misleading some people to a false target and then leading them back to the right one, when exit conditions are safer. Note that in real-life situations, most of the injuries are actually caused by overcompression and suffocation rather than urgency.  

\subsubsection{Mesoscopic model}
\begin{figure}[b!]
\centering
\includegraphics[width=0.3\textwidth]{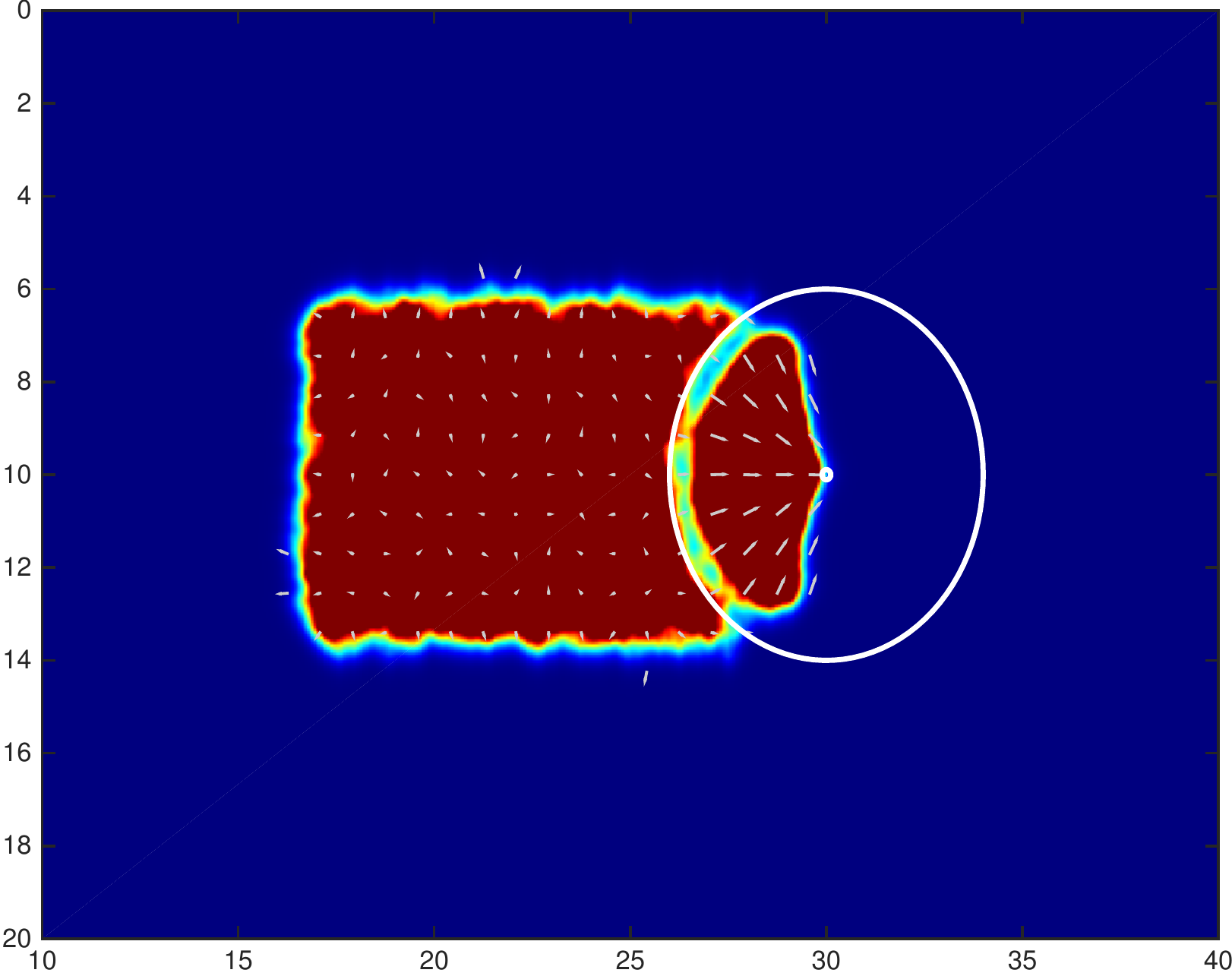}\quad
\includegraphics[width=0.3\textwidth]{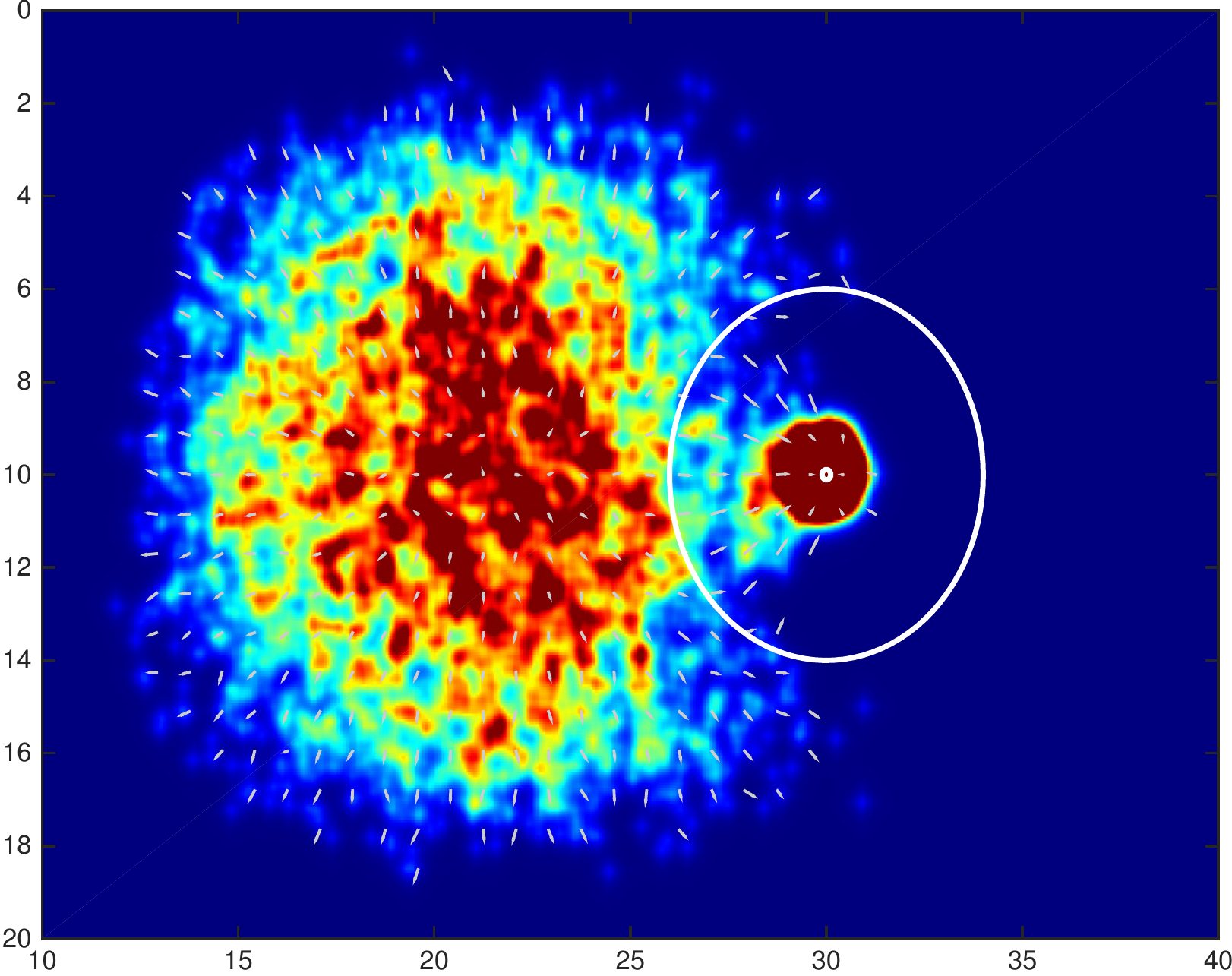}\quad
\includegraphics[width=0.3\textwidth]{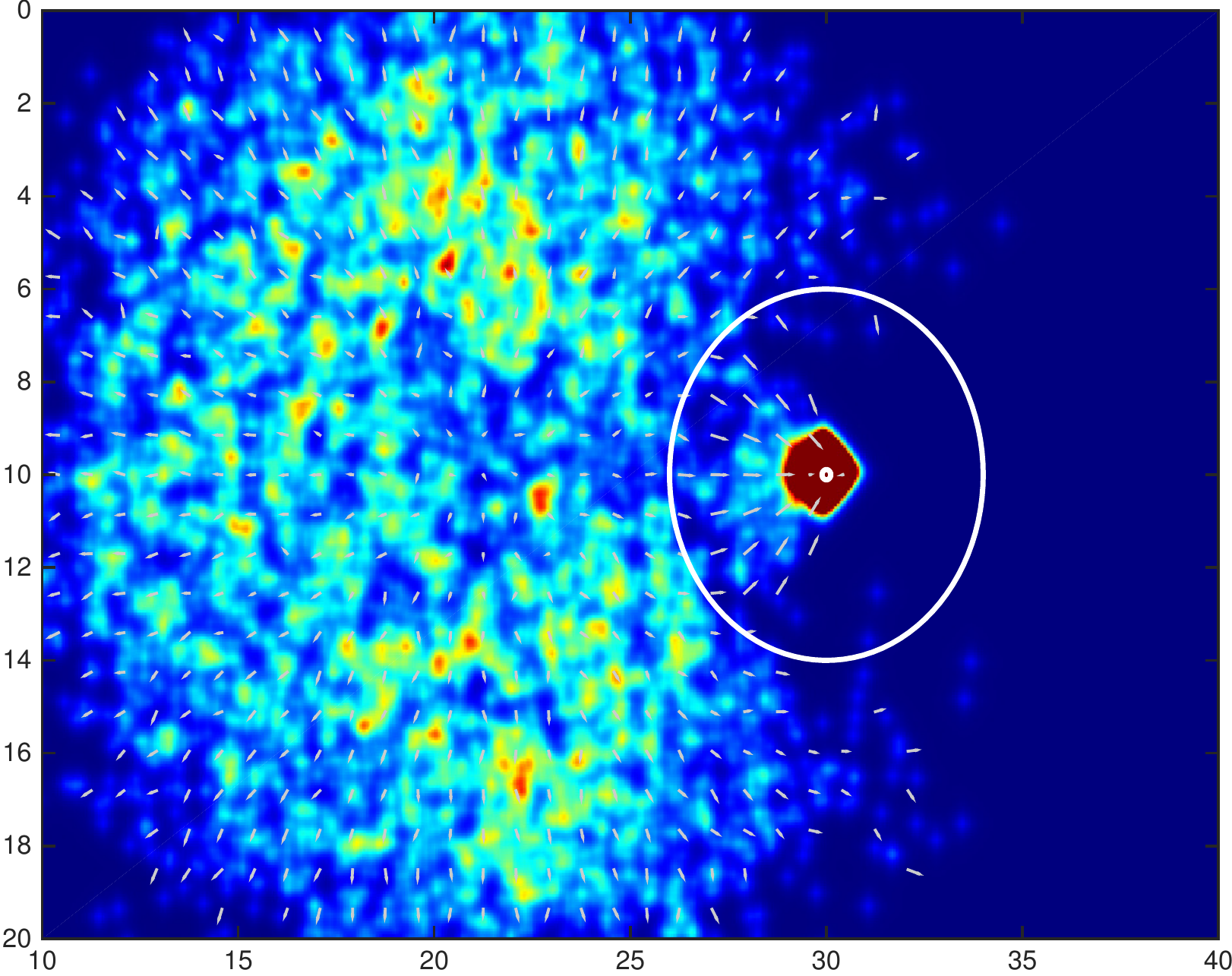}\\
\includegraphics[width=0.3\textwidth]{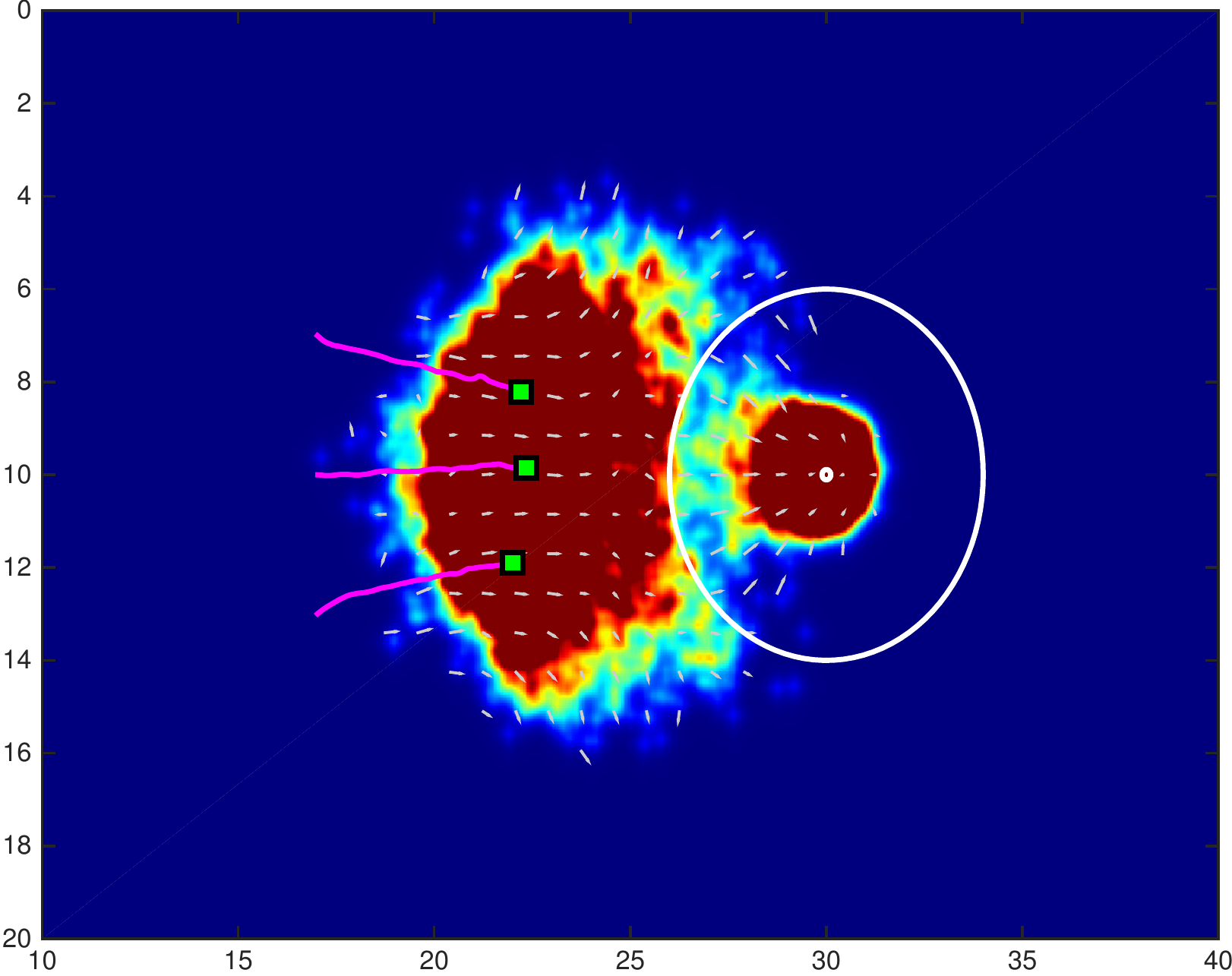}\quad
\includegraphics[width=0.3\textwidth]{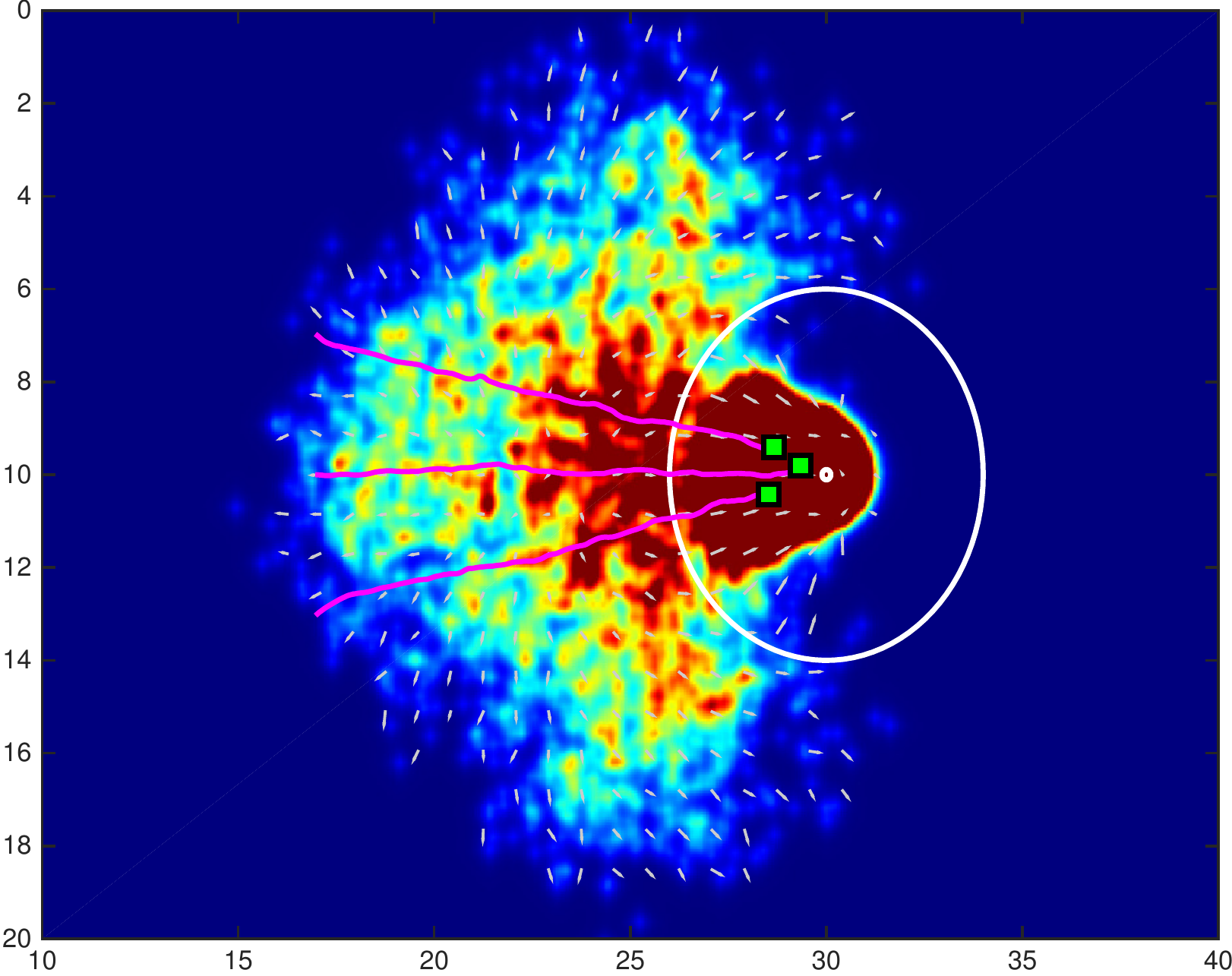}\quad
\includegraphics[width=0.3\textwidth]{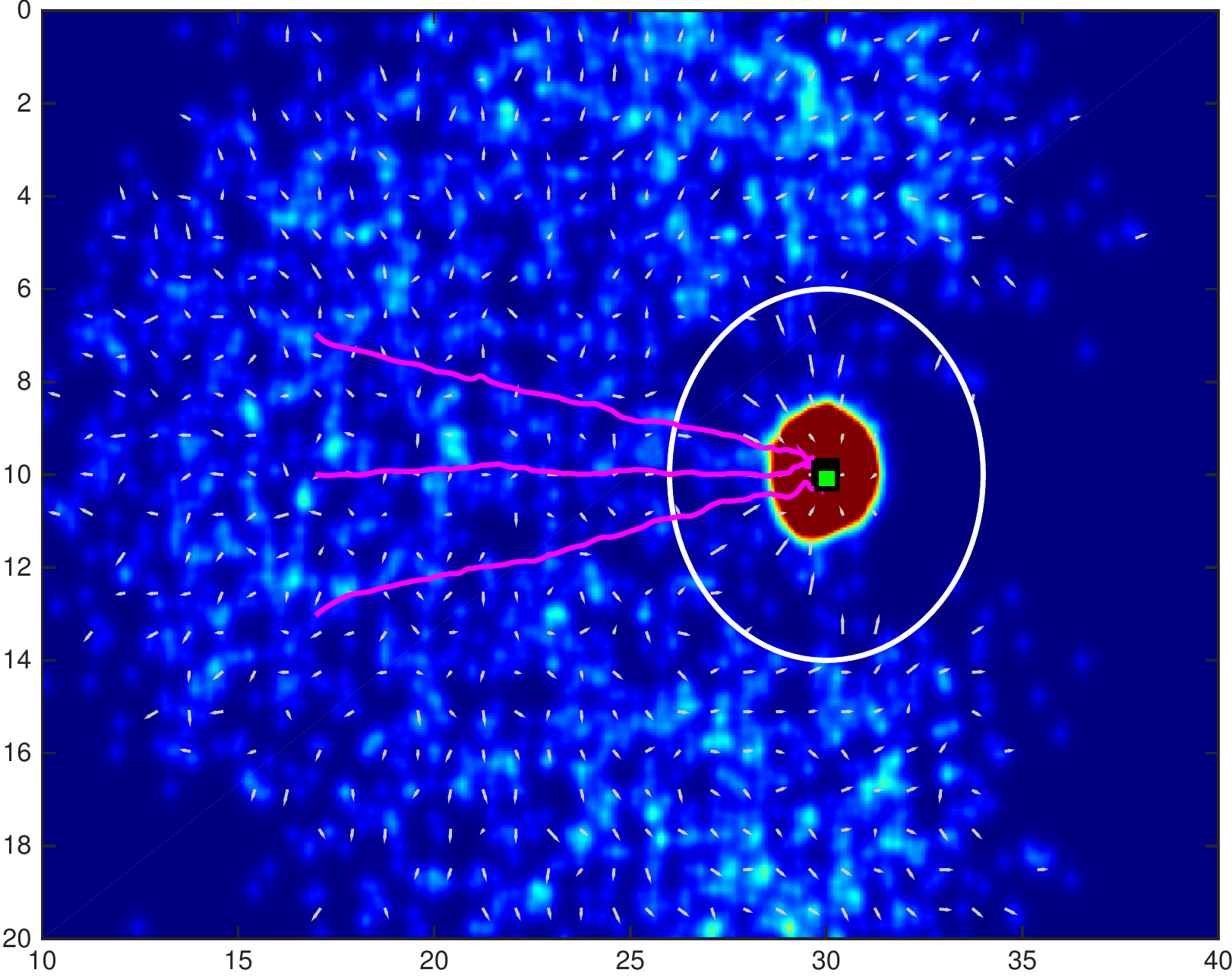}\\
\includegraphics[width=0.3\textwidth]{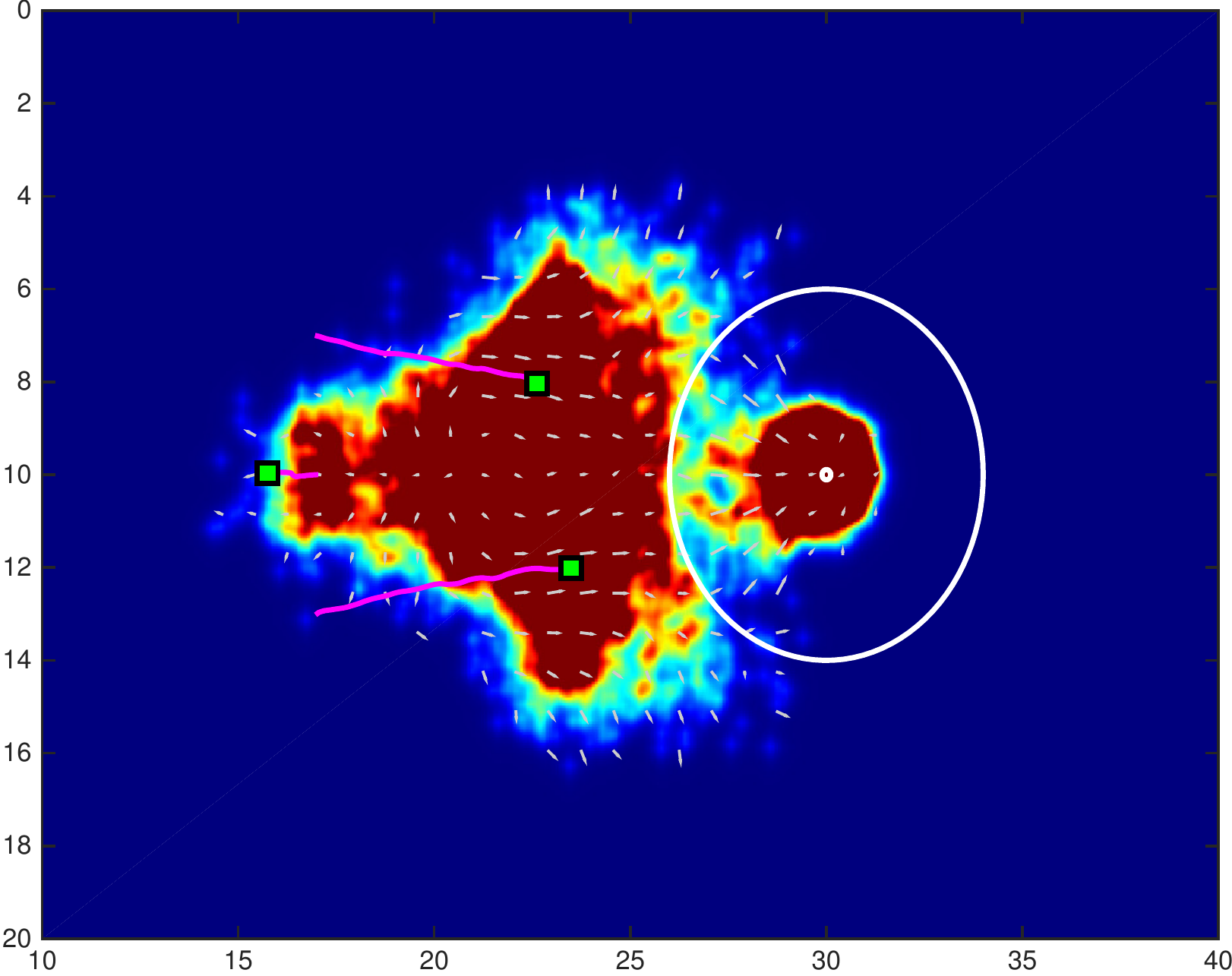}\quad
\includegraphics[width=0.3\textwidth]{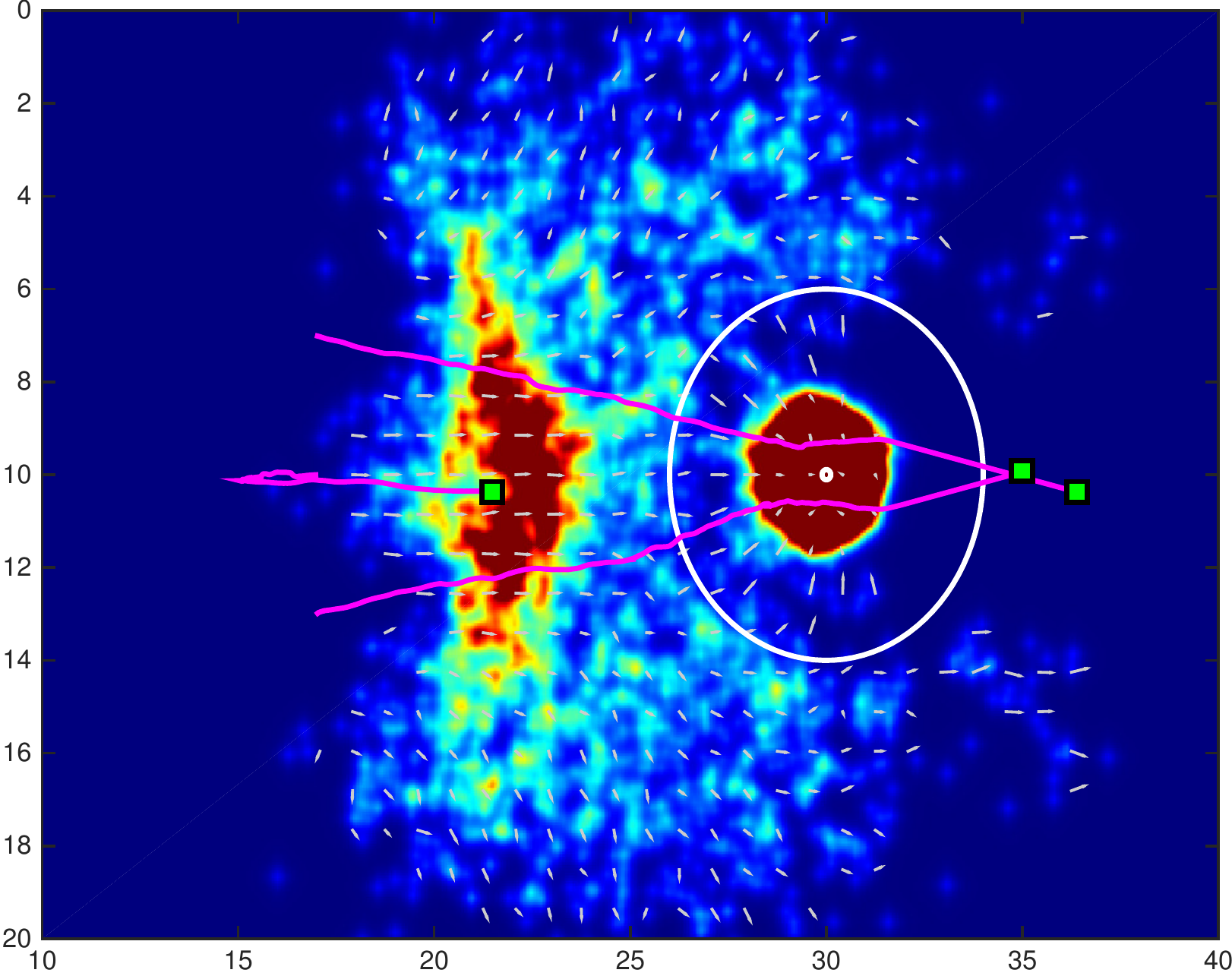}\quad
\includegraphics[width=0.3\textwidth]{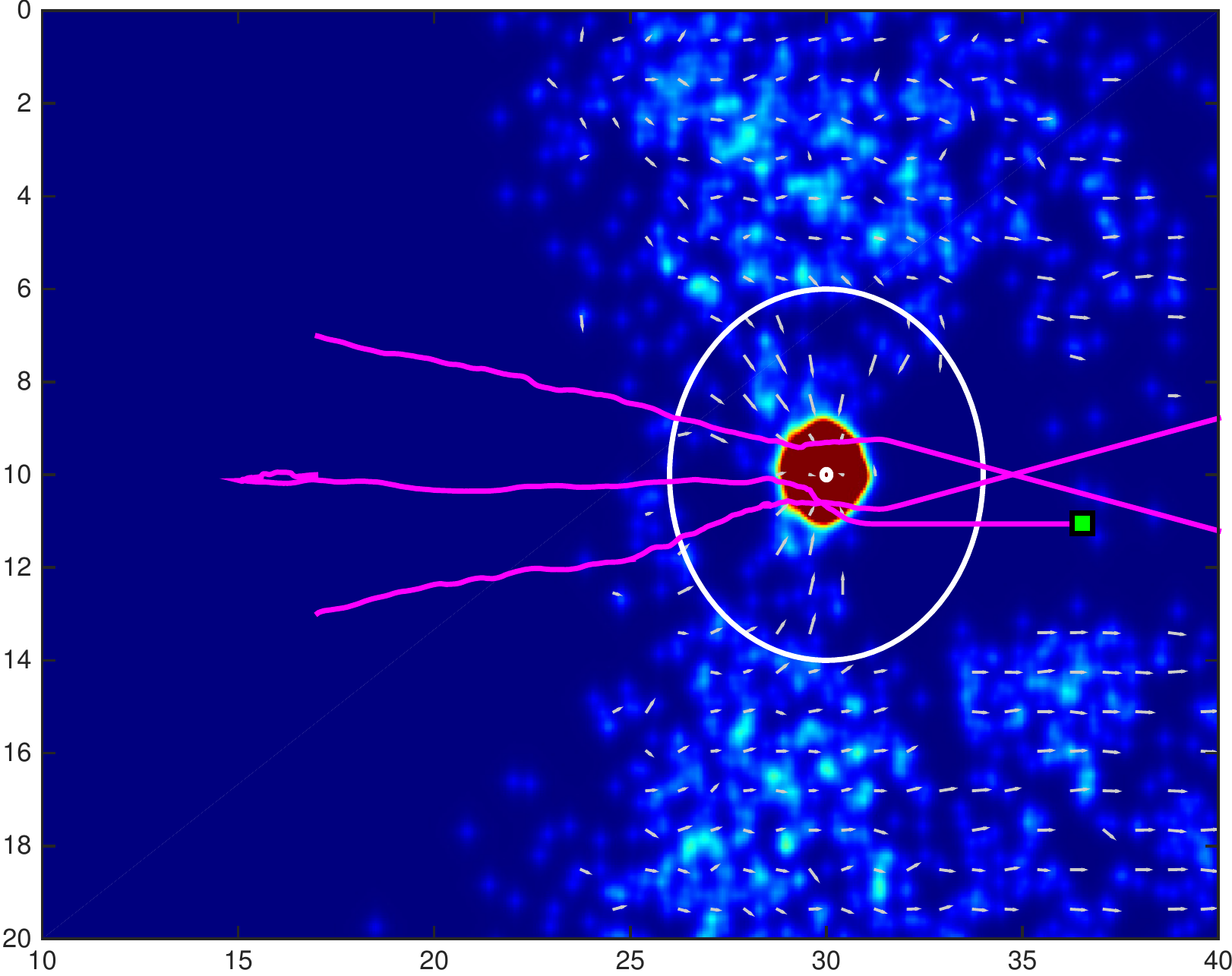}
\caption{Setting 1. Mesoscopic dynamics. First row: no leaders. Second row: three leaders, go-to-target strategy. Third row: three leaders, optimal strategy (compass search).}
\label{fig:S1_meso}
\end{figure}

We consider here the case of a continuous density of followers. Figure \ref{fig:S1_meso}(first row) shows the evolution of the uncontrolled system of followers. Due to the diffusion term and the topological alignment, large part of the mass spreads around the domain and is not able to reach the target exit.

In Fig.\ \ref{fig:S1_meso}(second row) we account the action of three leaders, driven by a go-to-target strategy defined as in the microscopic case. It is clear that also in this case the action of leaders is able to influence the system and promote the evacuation, but the presence of the diffusive term causes the dispersion of part of the continuous density. The result is that part of the mass is not able to evacuate, unlike the microscopic case. 

In order to improve the go-to-target strategy we rely on the compass search, where, differently from the microscopic case, the optimization process accounts the objective functional \eqref{defJ_evacmass}, i.e., the total mass evacuated at final time.
Figure \ref{fig:S1_meso}(third row) sketches the optimal strategy found in this way: on one hand, the two external leaders go directly towards the exit, evacuating part of the density; on the other hand the central leader moves slowly backward, misleading part of the density and only later it moves forwards towards the exit. The efficiency of the leaders' strategy is due in particular by the latter movement of the last leader, which is able to gather the followers' density left behind by the others, and to reduce the occupancy of the exit's visibility area by delaying the arrival of part of the mass. 

\begin{figure}[h!]
\begin{tabular}{ccc}

\includegraphics[width=0.3\textwidth]{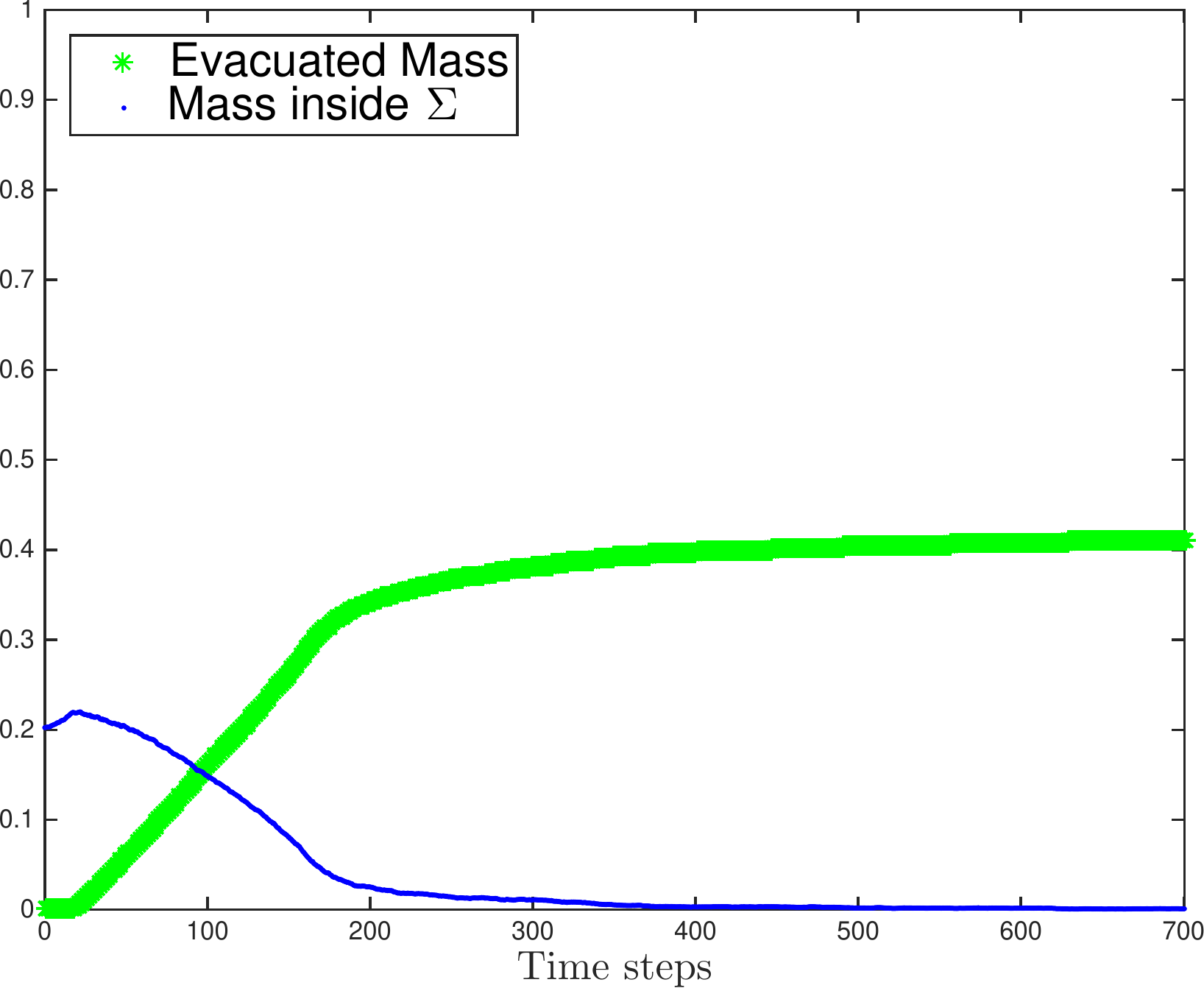}
&
\includegraphics[width=0.3\textwidth]{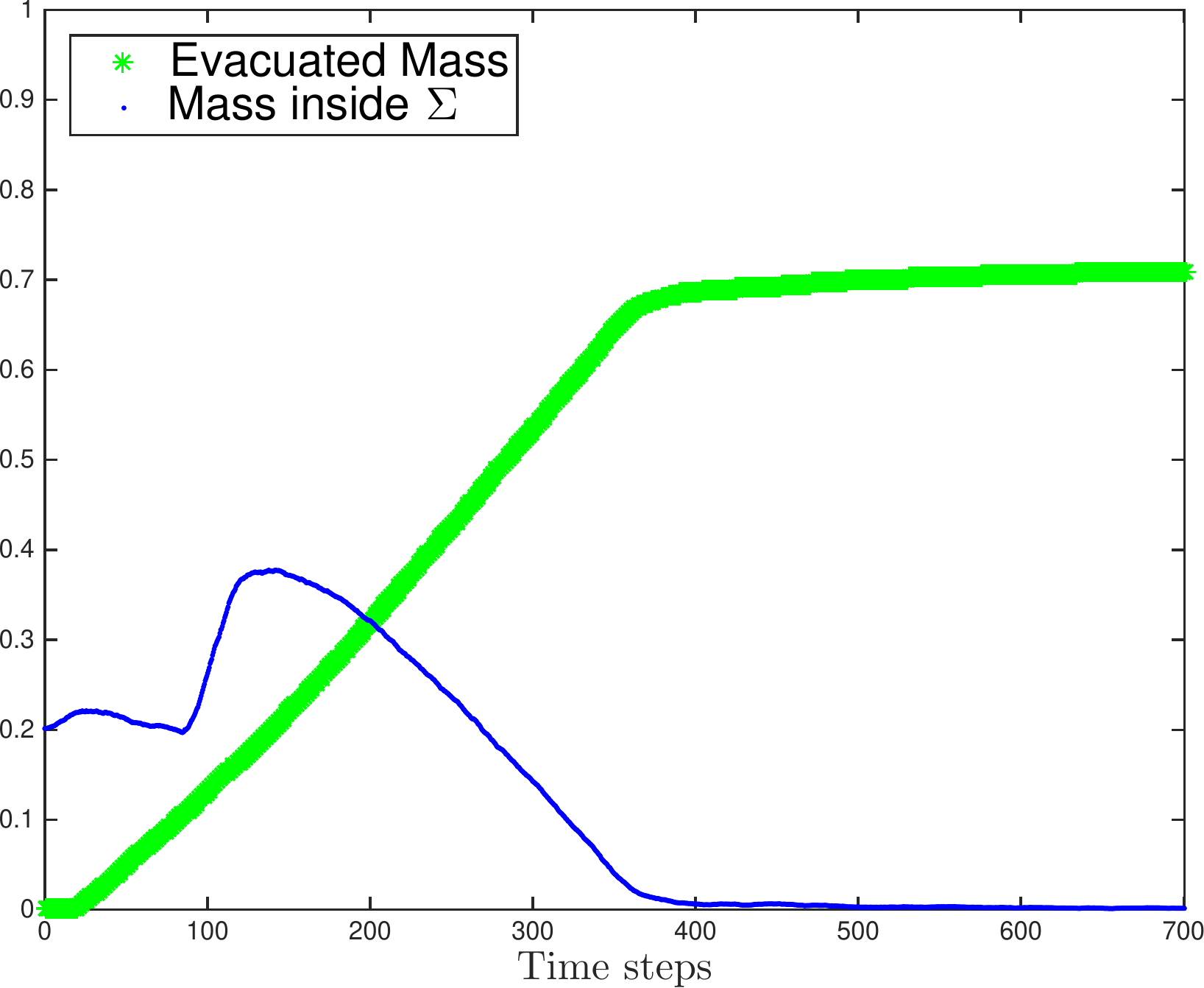}
&
\includegraphics[width=0.3\textwidth]{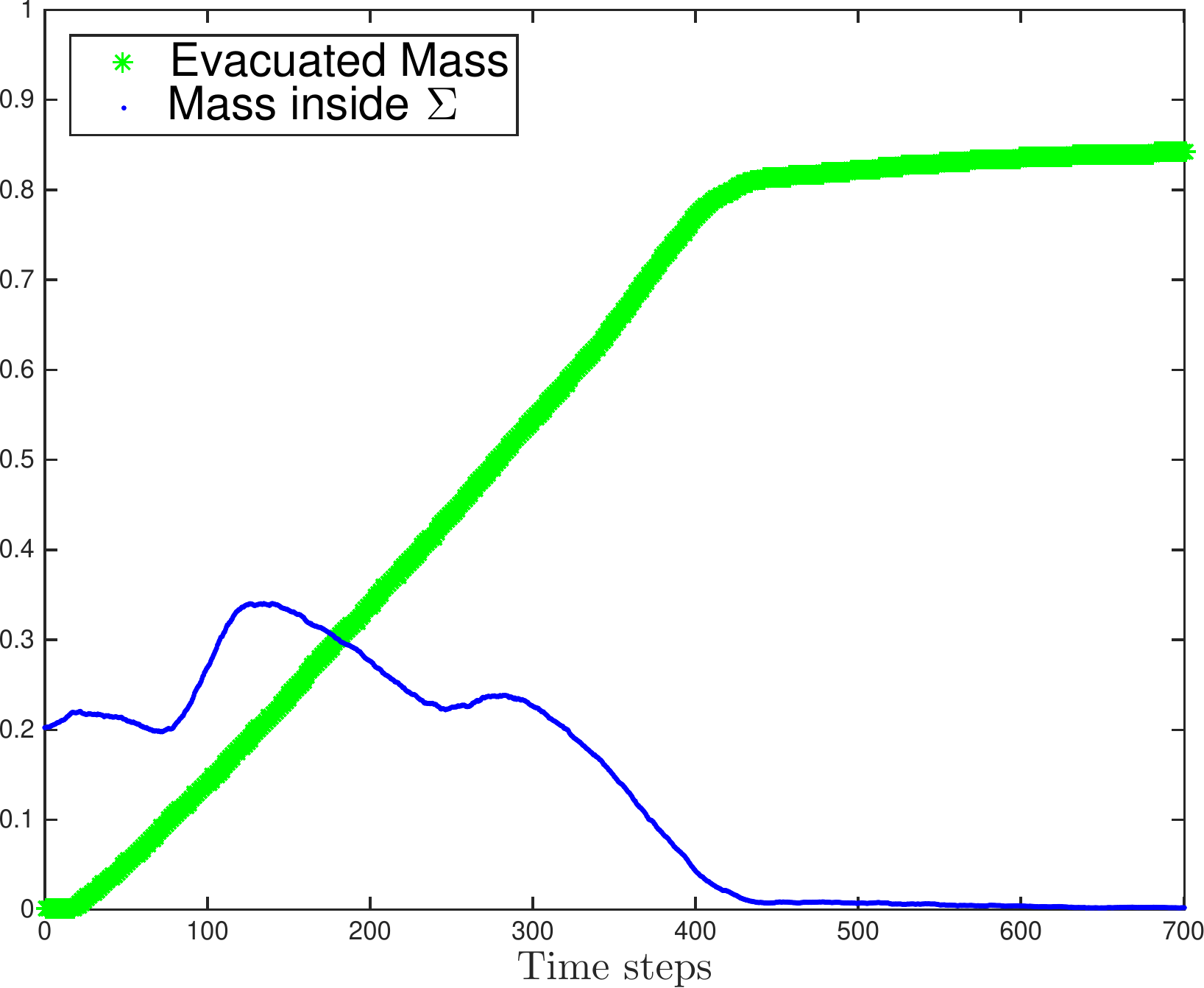}
\\
\end{tabular}
\caption{Setting 1. Occupancy of the exit's visibility zone $\Sigma$, dotted lines, and percentage of evacuated mass, star lines, as function of time. Left: histograms for the case without leaders (percentage of evacuated mass 41.2\%). Center: histograms for leaders moving with the go-to-target strategy (percentage of evacuated mass 71.3\%). Right: histograms for leaders with an optimal strategy (percentage of evacuated mass 85.2\%).}
\label{fig:S1_meso_hystogram}
\end{figure}

In Figure \ref{fig:S1_meso_hystogram} we summarize, for the three numerical experiments, the evacuated mass and the occupancy of the exit's visibility area $\Sigma$ as functions of time. The occupancy of the exit's visibility area shows clearly the difference between the leaders' action: for the go-to-target strategy, the amount of mass occupying  $\Sigma$ concentrates quickly and the evacuation is partially hindered by the clogging effect, as only 71.3\% of the total mass is evacuated. The optimal strategy (obtained after 30 iterations) is able to better distribute the mass arrival in $\Sigma$, and an higher efficiency is reached, evacuating up to 85.2\% of the total mass.
%
%
\begin{remark}
In a realistic context it is impossible to know a priori the initial position  and velocity of a crowd, but they can be expressed as a joint probability density, obtained from statistics on real data. Therefore from the simulation of the kinetic model \eqref{eq:FokkerPlanckStrong}, one can estimate the evolution of these probabilities and testing the efficiency of specific control strategies.
\end{remark}

\subsection{Setting 2}
In the following we test the microscopic model (with compass search optimization) in a more complicated setting, which has also some similarities with the one considered in Section \ref{sec:theexperiment}. The crowd is initially confined in a rectangular room with three walls. In order to evacuate, people must first leave the room and then search for the exit point. We assume that walls are not visible, i.e., people can perceive them only by physical contact. This corresponds to an evacuation in case of null visibility (but for the exit point which is still visible from within $\Sigma$). Walls are handled as in \cite{cristiani2011MMS} (see also \cite{cristiani_obstacles} for an overview of the techniques employed to handle obstacles in pedestrian simulations).

\begin{figure}[h!]
\centering
\includegraphics[width=0.3\textwidth]{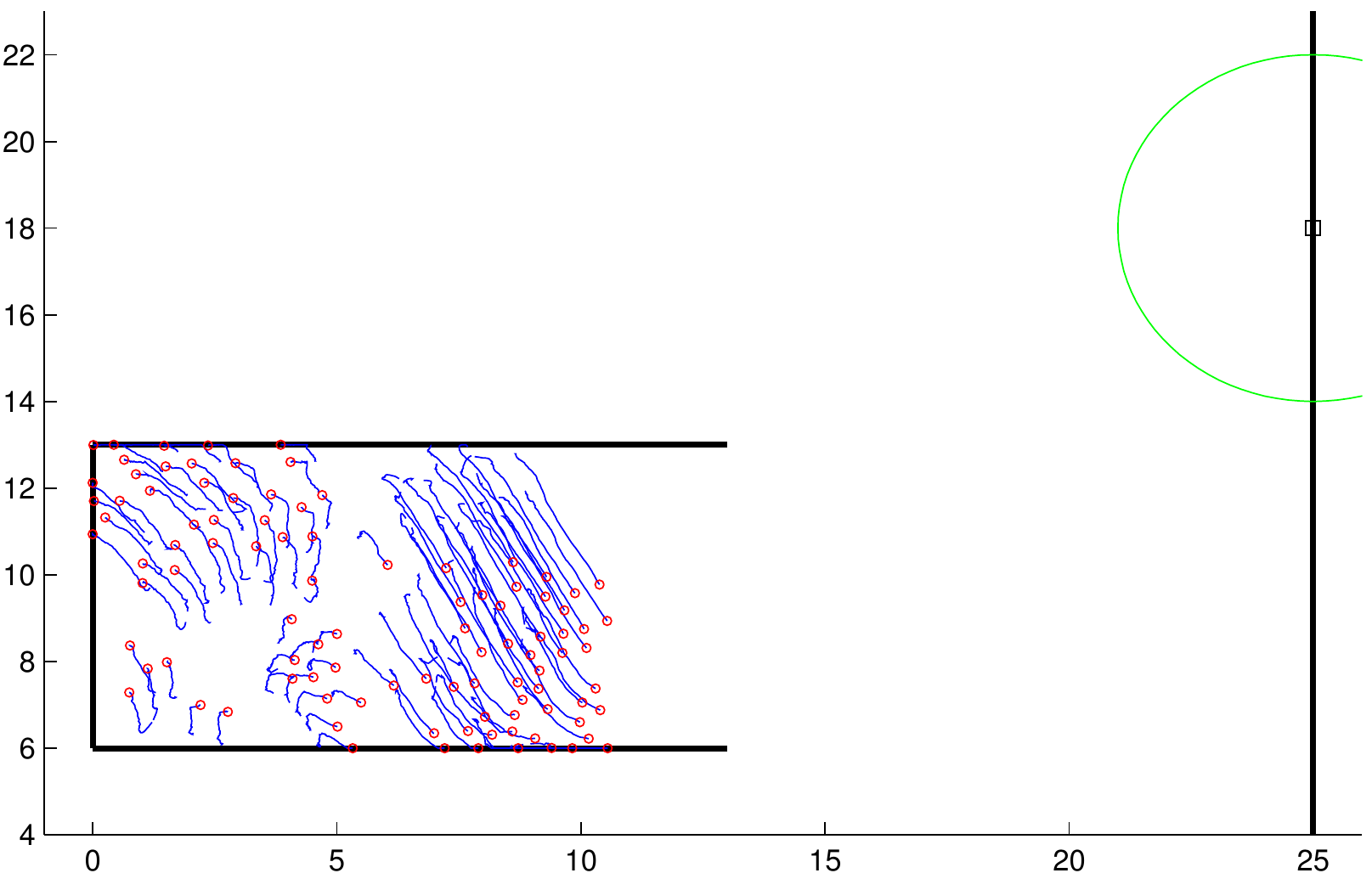}\quad
\includegraphics[width=0.3\textwidth]{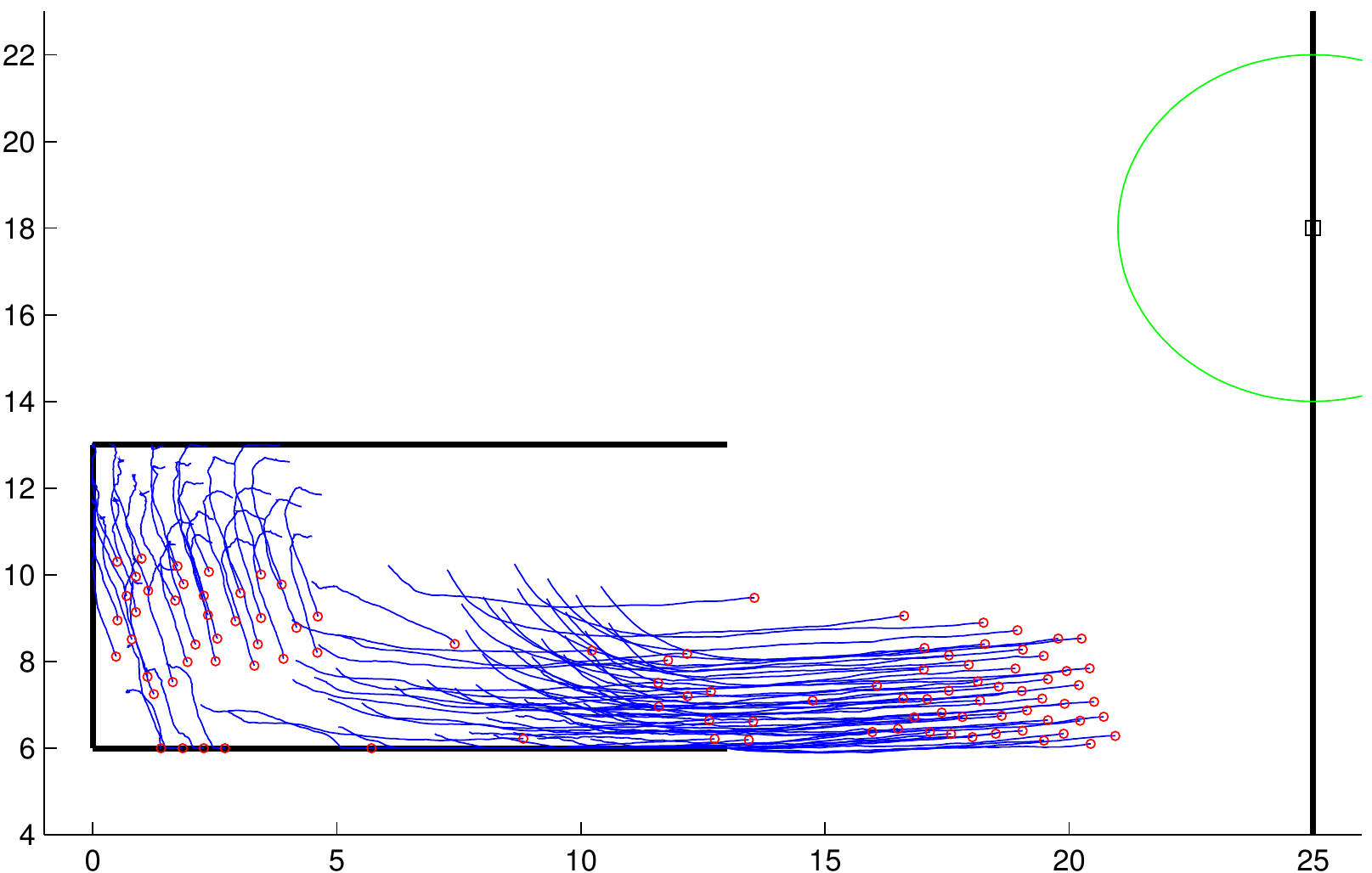}\quad
\includegraphics[width=0.3\textwidth]{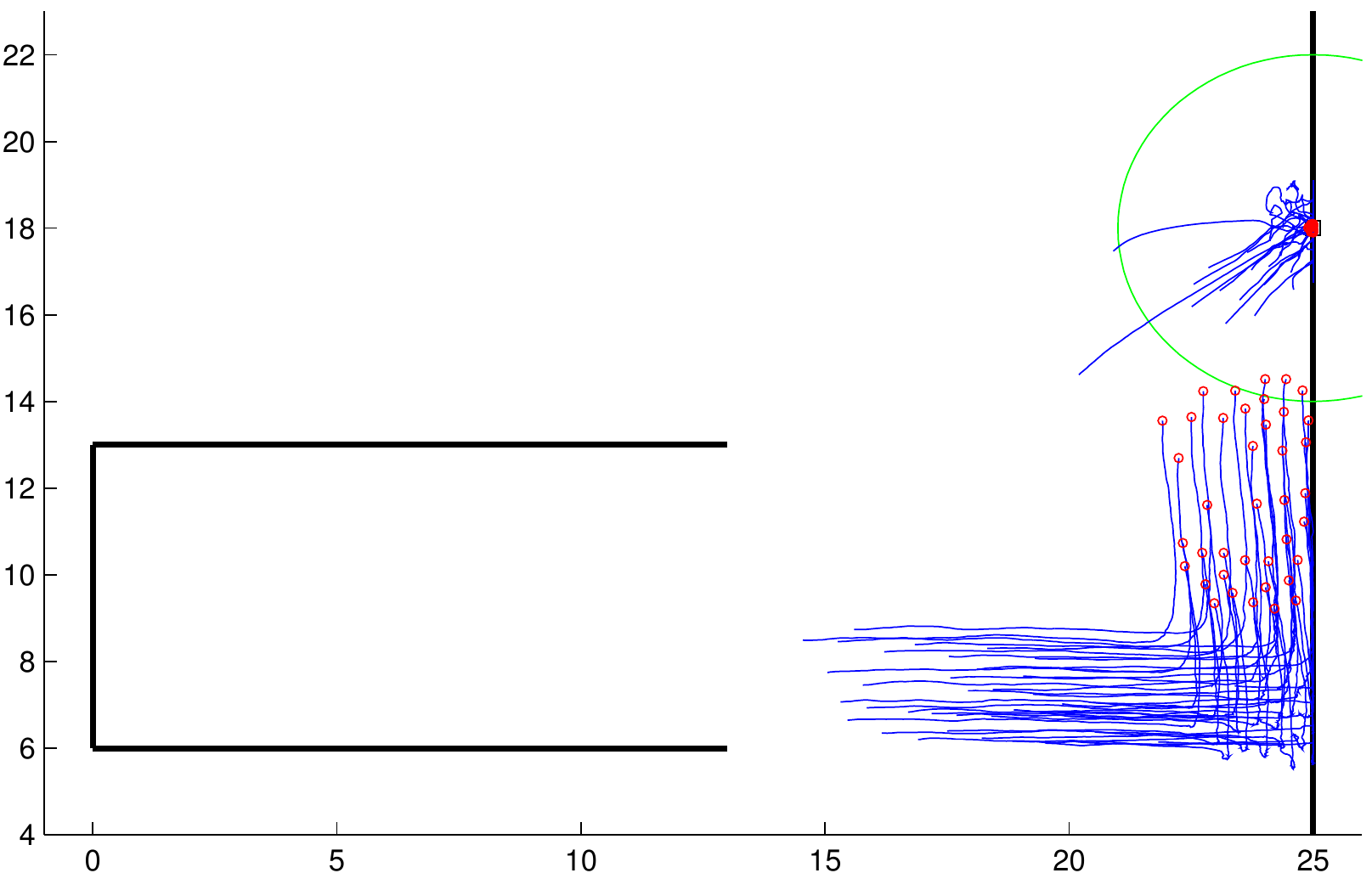}\\ [2mm]
\includegraphics[width=0.3\textwidth]{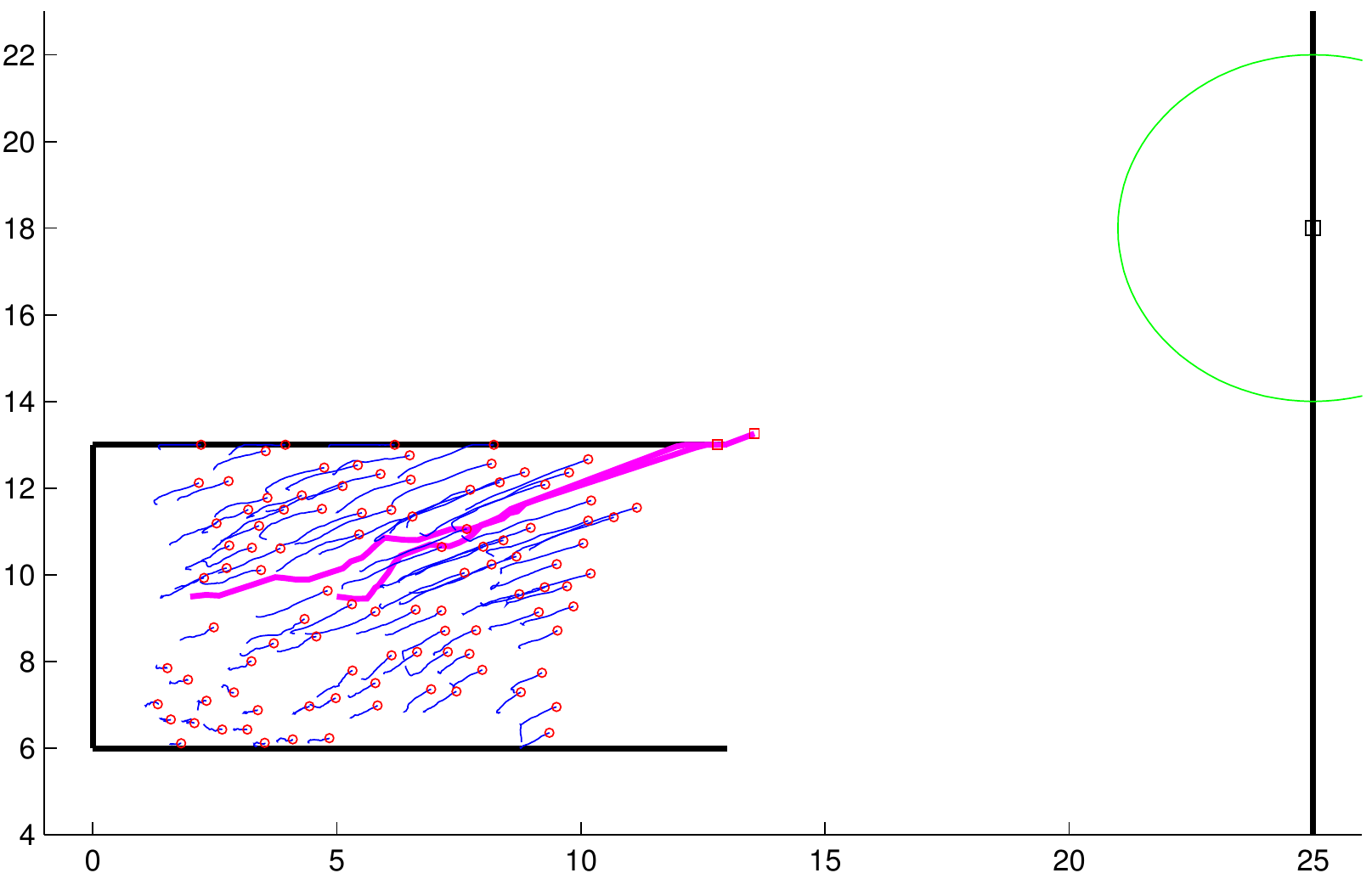}\quad
\includegraphics[width=0.3\textwidth]{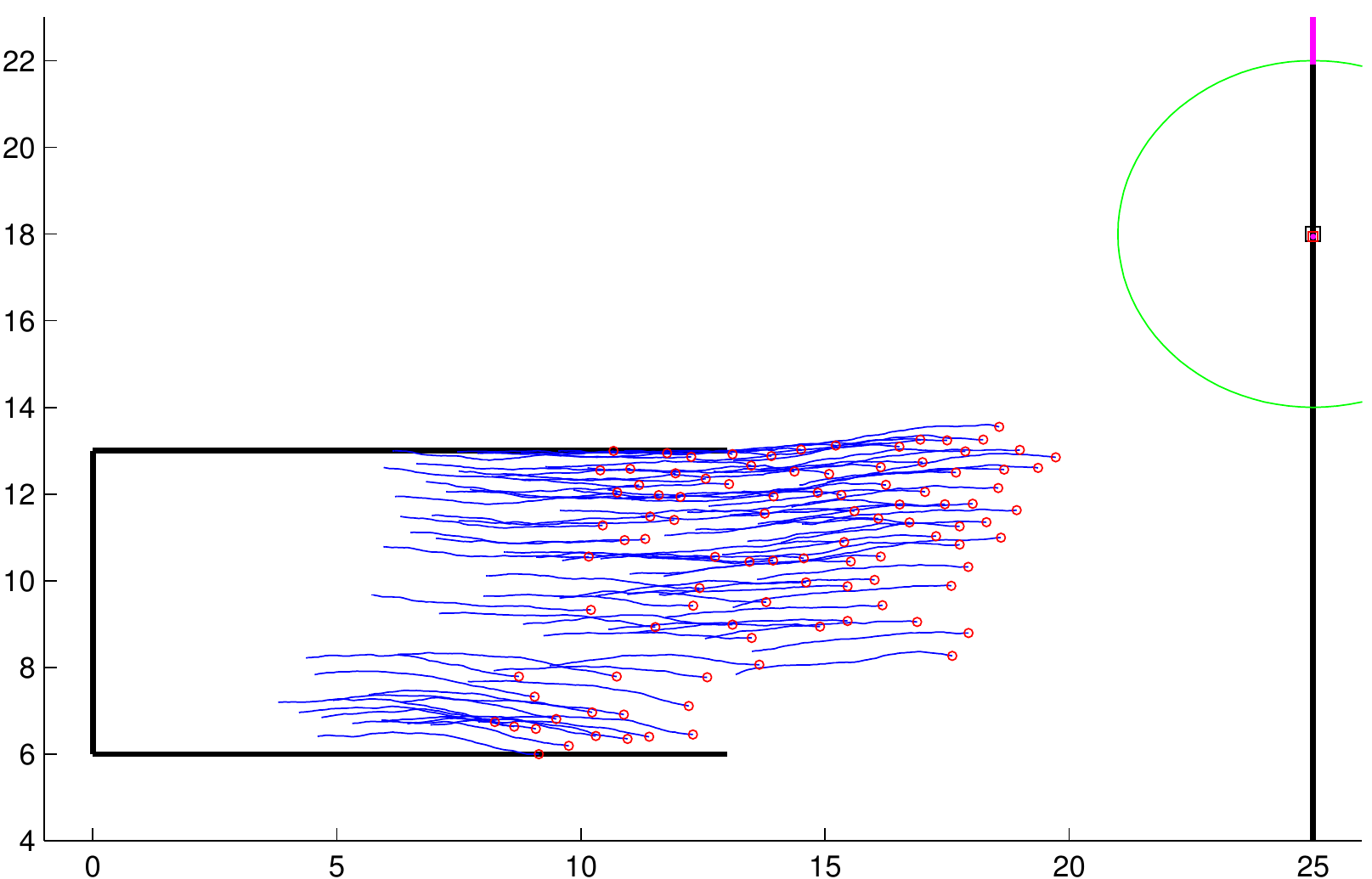}\quad
\includegraphics[width=0.3\textwidth]{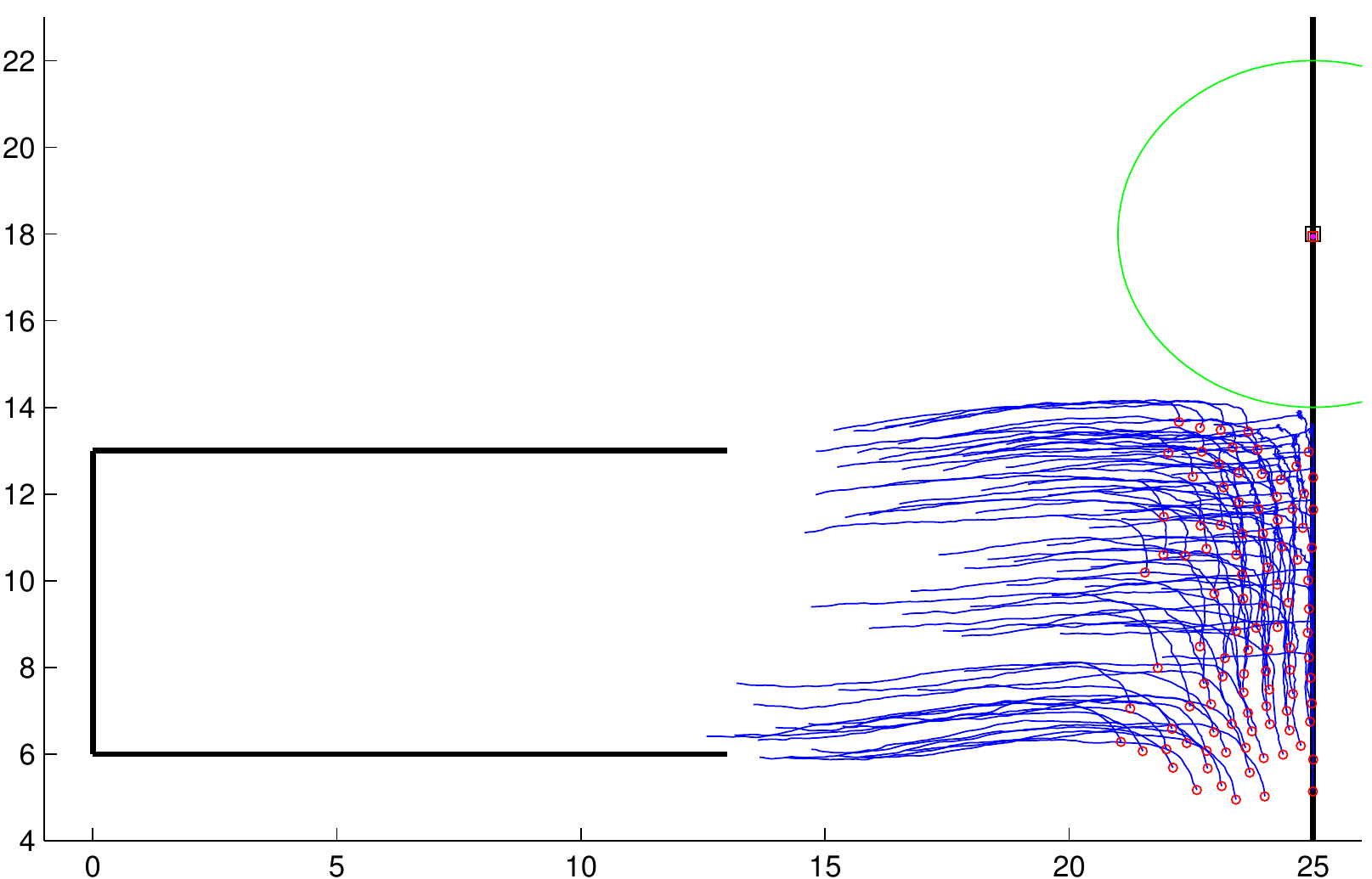}\\ [2mm]
\includegraphics[width=0.3\textwidth]{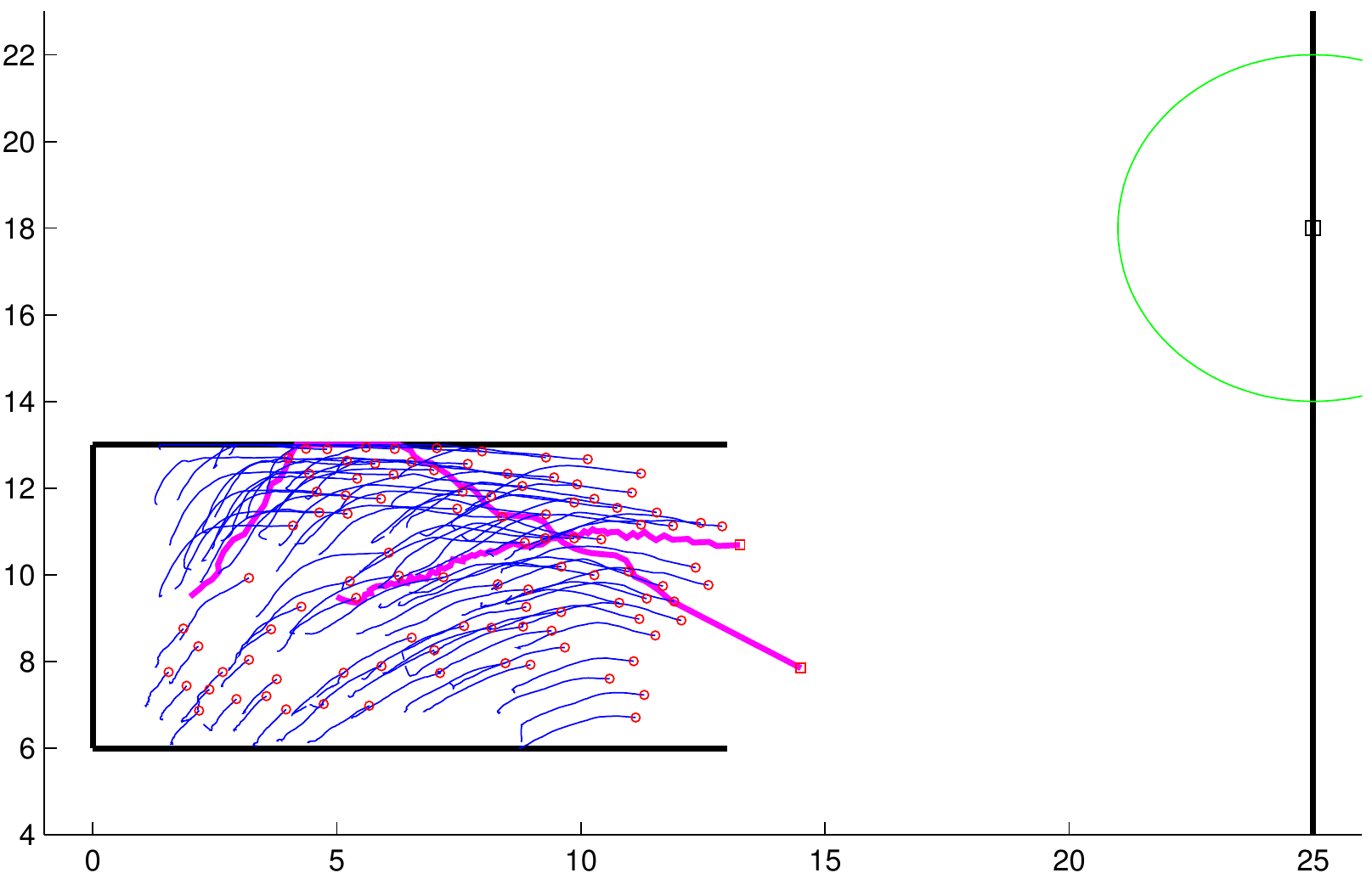}\quad
\includegraphics[width=0.3\textwidth]{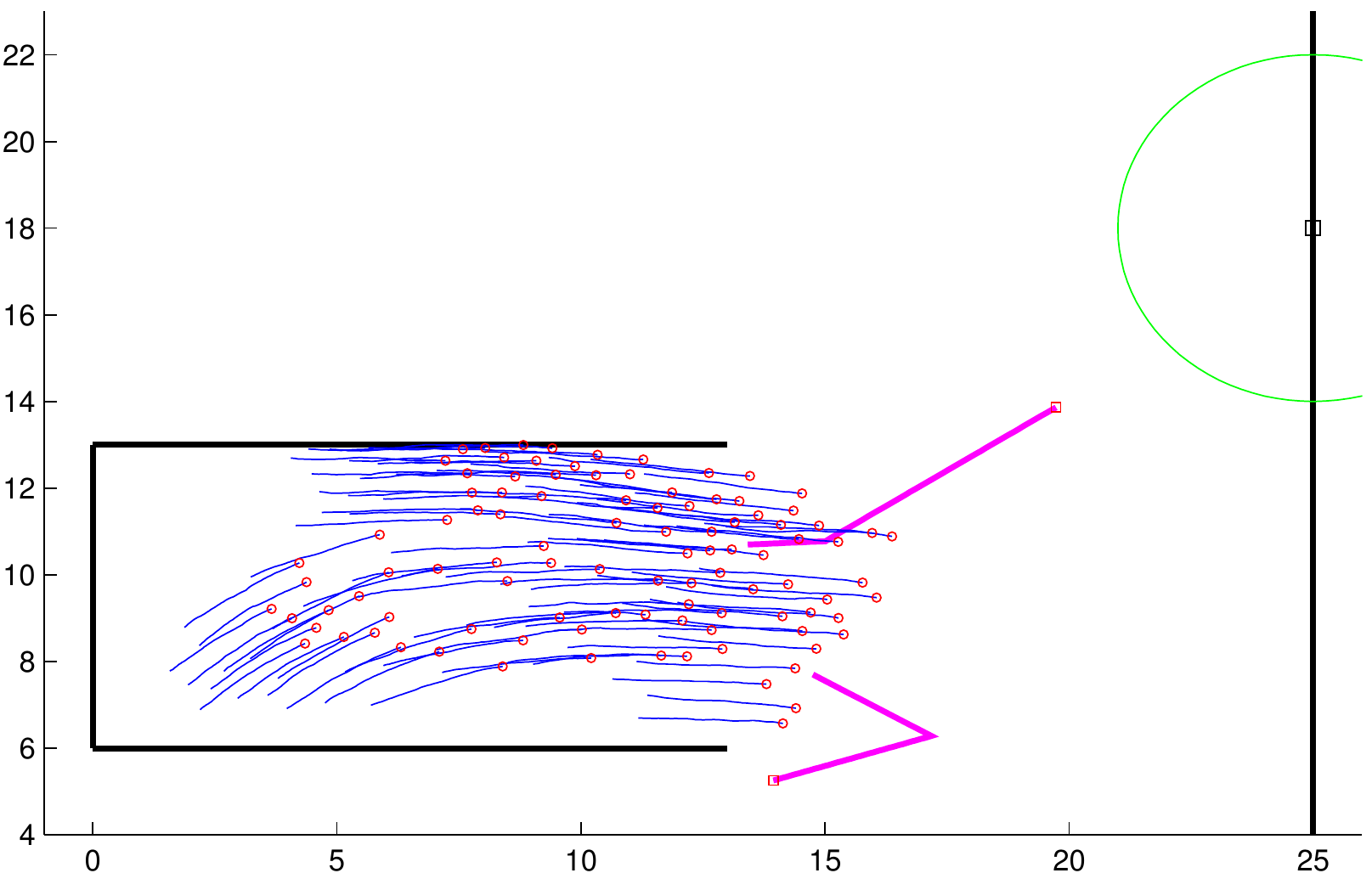}\quad
\includegraphics[width=0.3\textwidth]{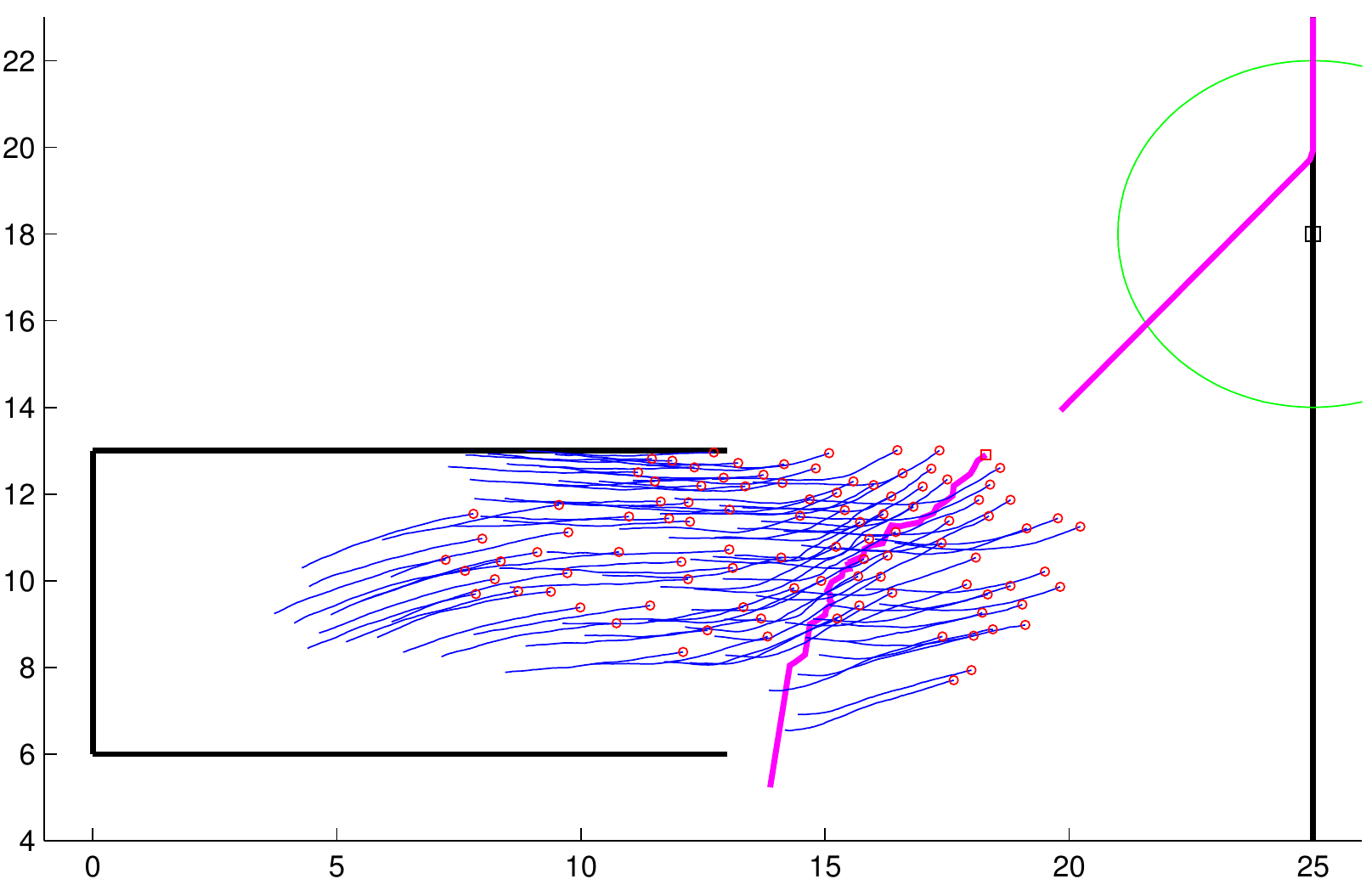}
\caption{Setting 2. Microscopic simulation. 
First row: no leaders (\textit{\href{www.emilianocristiani.it/attach/TUM_test_experiment_no_control_1200.mp4}{video}}).
Second row: two leaders and go-to-target strategy (\textit{\href{www.emilianocristiani.it/attach/TUM_test_experiment_wrongfixedstrategy_oo.mp4}{video}}). 
Third row: two leaders and best strategy computed by the compass search (\textit{\href{www.emilianocristiani.it/attach/TUM_test_experiment_best_control_549.mp4}{video}}).}
\label{fig:S2}
\end{figure}

If no leaders are present, the crowd splits in several groups and most of the people hit the wall, see Fig.\ \ref{fig:S2}(first row). After some attempts the crowd finds the way out, and then it crashes into the right boundary of the domain. Finally, by chance people decide, \emph{en cascade}, to go upward. The crowd leaves the domain in 1162 time steps.

If instead we hide in the crowd two leaders who point fast towards the exit (Fig.\ \ref{fig:S2}(second row)), the evacuation from the room is completed in very short time, but after that, the influence of the leaders vanishes. Unfortunately, this time people decide to go downward after hitting the right boundary, and nobody leaves the domain.  
Slowing down the two leaders helps keeping the leaders' influence for longer time, although it is quite difficult to find a good choice. 

Compass search optimization finds (after 30 iterations) a nice strategy for the two leaders which remarkably improves the evacuation time, see Fig.\ \ref{fig:S2}(third row). One leader behaves similarly to the previous case, while the other diverts the crowd pointing SE, then comes back to wait for the crowd, and finally points NE towards the exit. This strategy allows to bring everyone to the exit in 549 time steps, without bumping anyone against the boundary, and avoiding congestion near the exit. 

\subsection{Setting 3} 
In the following test we propose a different use of the leaders. 
Rather than steering the mass towards the exit, they can be employed to fluidify the evacuation near a door. 
In fact, it is observed that when several pedestrians reach a bottleneck at the same time they slow down each other, possibly coming to a deadlock. As a consequence, none of them can pass through the bottleneck. 
We want to investigate the possibility of using leaders as small \textit{smart obstacles}, which alleviate the high friction among the bodies just moving randomly near the exit.
\begin{figure}[h!]
\centering
\includegraphics[width=0.29\textwidth]{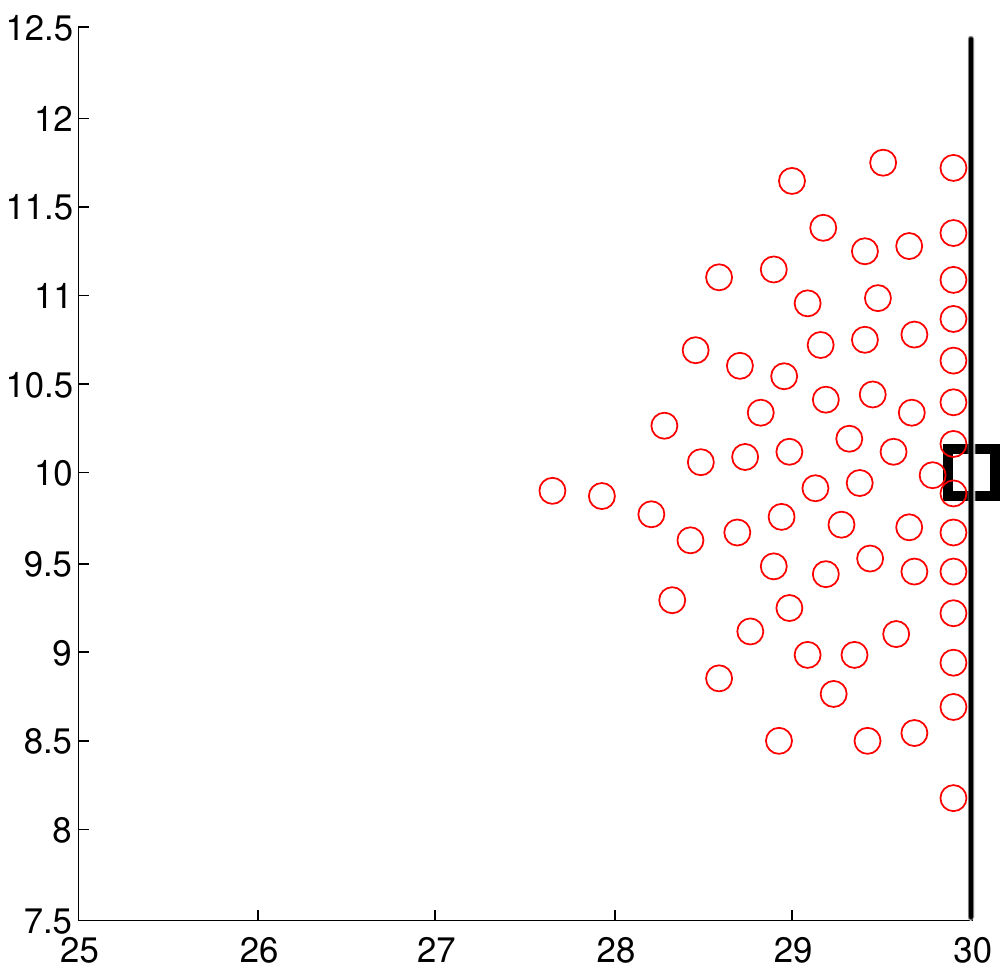}
\includegraphics[width=0.29\textwidth]{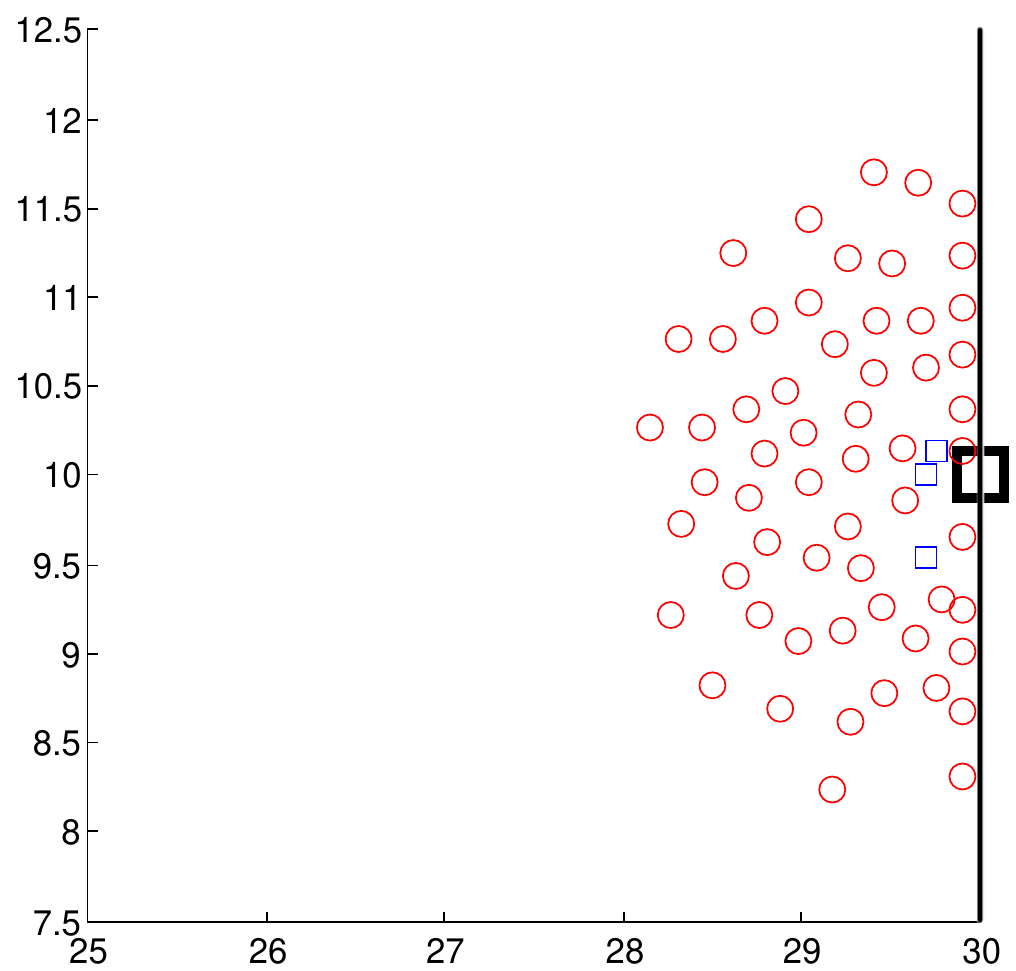}\quad
\includegraphics[width=0.37\textwidth]{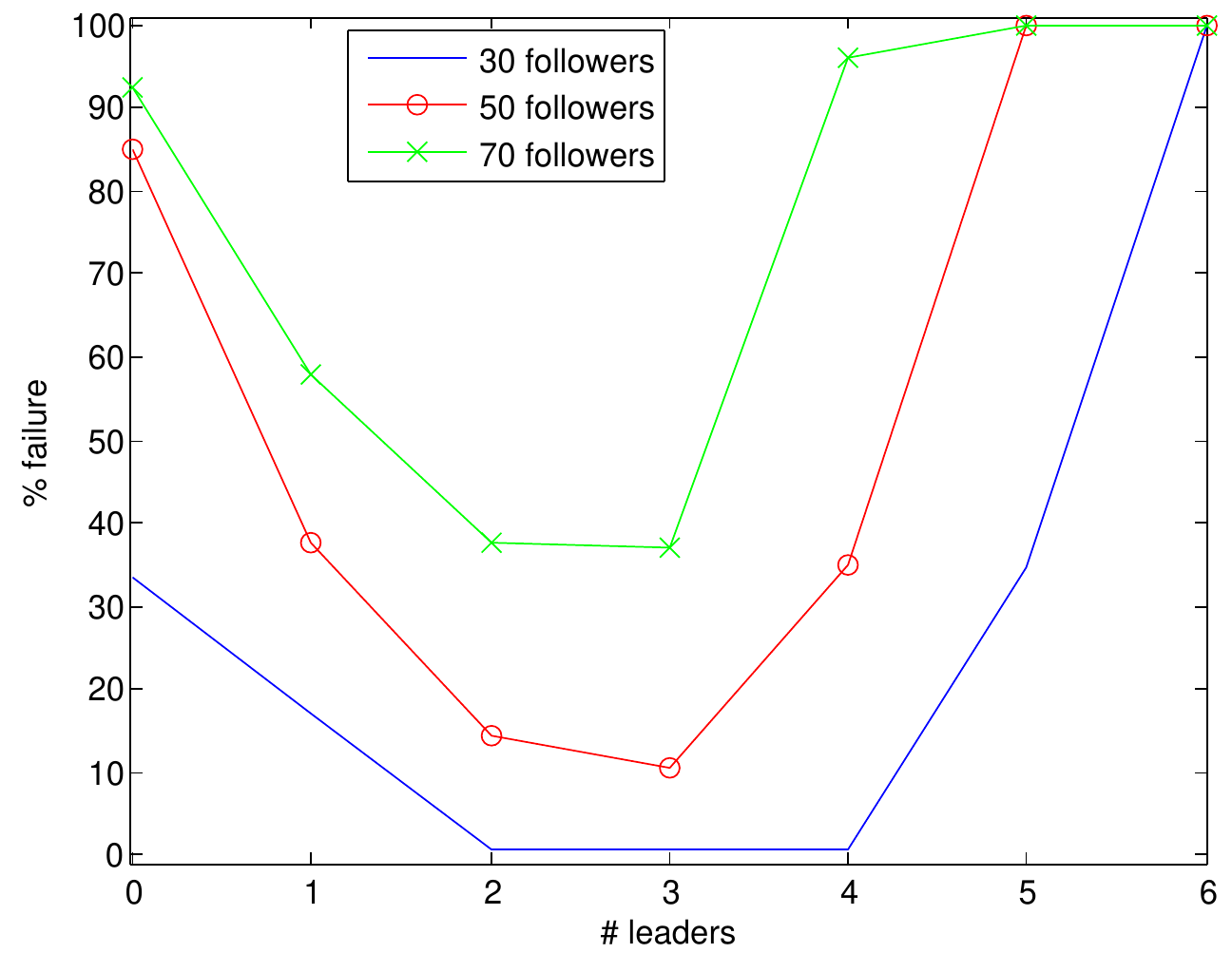}
\caption{Setting 3 (exit is visible from any point). Left: evacuation with no leaders. Center: evacuation with 3 leaders. Right: percentage of evacuation failures as a function of the number of leaders, for $\NF$=30, 50, 70.}
\label{fig:S3}
\end{figure}

We consider a large room with a single exit door located on the right wall, visible from any point ($\Sigma=\R^2$). We nullify the repulsion force perceived by the leaders, see Table \ref{tab:all_parameters}. Most important, in this test we assume that followers have nonzero size  and hard collisions are not possible. This is achieved assuming that they are circular-shaped with diameter 0.25, and keeping frozen those pedestrians who aim at moving too close to another pedestrian in one time step. 
The width of the exit is 0.45. Leaders stay near the exit at any time and move with a random velocity, each component being uniformly distributed in $[-1.3,1.3]$. 
We say that the evacuation succeeded if $\NF-10$ followers leave the room within 2000 time steps, otherwise we say that the evacuation failed. 
Average quantities are computed running the simulation 200 times.  

Fig.\ \ref{fig:S3} shows a typical outcome without (left) and with (center) leaders.  
Fig.\ \ref{fig:S3}-right shows the percentage of evacuation failures as a function of the number of leaders. It is clearly visible that in this scenario the optimal number of leaders is 3, and this number is independent from the number of evacuees. This is expected since leaders act at local level near the exit and do not affect people queuing far from the exit. We noted that the leaders initially act as a barrier for the incoming pedestrians and slow down the evacuation. Later on, instead, they are quite efficient in destroy symmetries and avoid deadlocks. 
As a by-product, we confirmed once again the Braess's paradox, which states that reducing the freedom of choice (e.g., adding obstacles in front of the exit \cite{cristiani2015SIAP, twarogowska2014AMM}) one can actually improve the overall dynamics.
%
%
%
%
%
\section{Validation of the proposed crowd control technique}\label{sec:theexperiment}
The crowd control technique investigated in the previous sections relies on the fact that pedestrians actually exhibit herding behavior in special situations, so that the alignment term in the model is meaningful. The herding behavior is investigated since long time, especially in the psychological literature. To further confirm this behavior in our scenario, we have organized an experiment involving volunteer pedestrians.

The experiment took place on October 1, 2014 at Department of Mathematics of Sapienza -- University of Rome, Rome (Italy), at 1 p.m. Participants were chosen among first-year students. Since courses started just two days before the date of the experiment and the Department has a complex ring-shaped structure, we could safely assume that most of the students were unfamiliar with the environment. 

\paragraph{The experiment} Students were informed to be part of a scientific experiment aiming at studying the behavior of a group of people moving in a partially unknown environment. Their task was to leave the classroom II and reach a certain target inside the Department. Students were asked to accomplish their task as fast as possible, without running and without speaking with others. Participation was completely volunteer and not remunerated (only a celebratory T-shirt was given at later time).

76 students (39 girls, 37 boys) agreed to participate in the experiment. The students were then divided randomly in two groups: group A (42 people) and group B (34 people). The two groups performed the experiment independently one after the other. 
The target was communicated to the students just before the beginning of the experiment. It was the Istituto Nazionale di Alta Matematica (INdAM), which is located just upstairs w.r.t.\ the classroom II. Due to the complex shape of the environment, there are many paths joining classroom II and INdAM. The shortest path requires people to leave the classroom, to go leftward in a little frequented and unfamiliar area (even for experienced students), and climb the stairs. The target can be also reached by going rightward and climb other, more frequented, stairs.

In group B there were 5 incognito students, hereafter referred to as leaders, of the same age of the others. Leaders were previously informed about the goal of the experiment and the location of the target. They were also trained in order to steer the crowd toward the target in minimal time. Nobody recognized them as ``special'' before or during the experiment. Unexpectedly, also in group A there was a girl who knew the target. Therefore, she acted as an \textit{unaware} leader. 

It is important to stress that all the other students continued their usual activities and \emph{participants were not officially recorded by the organizers}. This choice was crucial to get \textit{natural} behavior and meaningful results, cf.\ \cite[Sect.\ 3.4.3]{cristiani2014book}. The price to pay is that we had to extract participants' trajectories by low-resolution videos taken by two observers. Moreover, the area just outside classroom II was rather crowded, introducing a high level of noise in the experiment. Finally, let us mention that some participants broke the rules of the experiment speaking with others and running (confirming the natural behavior...). Our experiment can be compared with those in \cite{dyer2009PTRSB}, although in those experiments group cohesion was imposed as a rule to the participants.

\begin{figure}[h!]
\begin{center}
\includegraphics[height=3.2cm,width=0.8\textwidth]{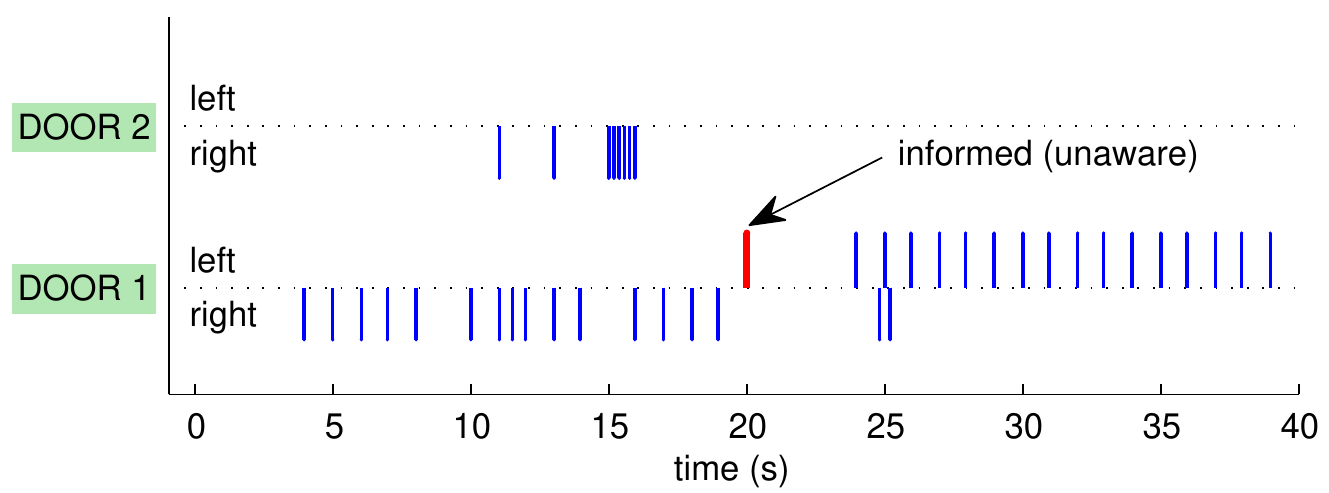}\\
\includegraphics[height=3.2cm,width=0.8\textwidth]{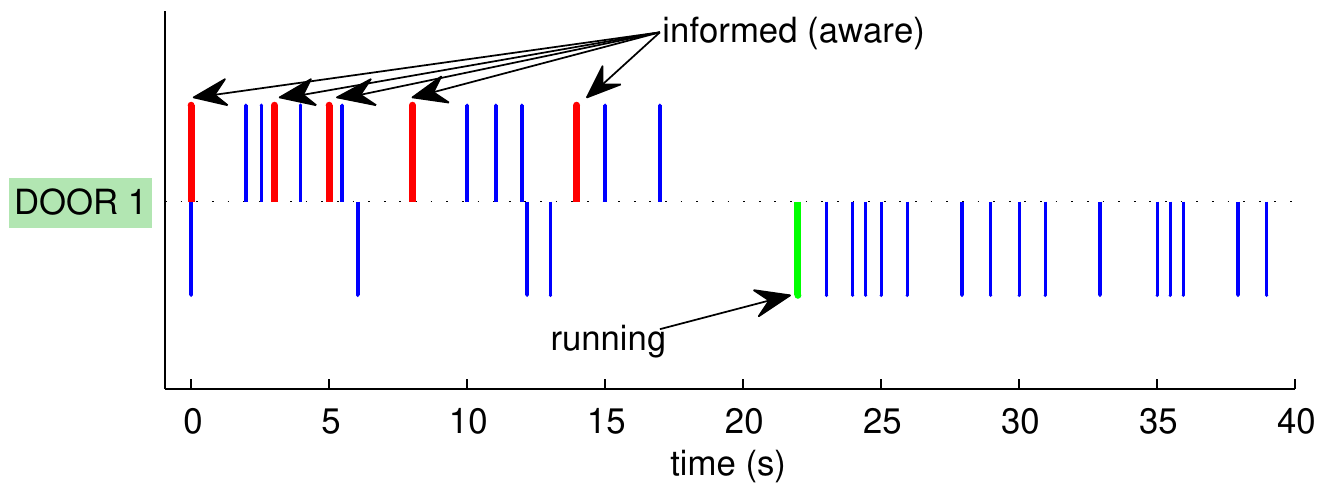}
\end{center}
\caption{Time instants corresponding to the students' right/left decisions in front of the doors for group A (top) and B (bottom).}
\label{fig:results_experiment}
\end{figure}

\paragraph{Results} 
In Figure \ref{fig:results_experiment} we show the history of left/right decisions taken by the students just outside the two doors of the classroom. 
Group A used both exits of the classroom. The 8-people group which used door 2 was uncertain about the direction for a while. Then, once the first two students decided to go rightward, the others decided simultaneously to follow them.
Similarly, the first student leaving from door 1 was uncertain for a few seconds, then he moved rightward and triggered a clear \textit{domino effect}. At $t=20$ s the unaware leader moved to the left, inducing hesitation and mixed behavior in the followers. After that, another domino effect arose. 

Group B used only door 1. Invisible leaders were able to trigger a domino effect but this time 4 people decided unilaterally not to follow them, although they were not informed about the destination. At $t=20$ s, after the passage of the last leader, a girl passed through the door and went leftward. Then, at $t=22$ s, she suddenly began to run toward the right. She first induced hesitation, then triggered a new rightward domino effect. 

\paragraph{Discussion}
Students have shown a tendency to go rightward, being that part of the Department more familiar, but this tendency was greatly overcame by the wish to follow the group mates in front, regardless of their right/left preferences. 
Also, after the experiment, the students admitted that they have been influenced by the leaders because of their clear direction of motion. Interestingly, no more than 5 people have been influenced by a single leader. This is compatible with the topological alignment term used in our model. The small value of $\mathcal N$ here is due to the fact that the space in front of the doors was rather crowded and the visibility was reduced. 

It is also interesting to note that only 4 students reached the target alone. This confirms the tendency not to remain isolated and to form clusters. 
Video recordings also show that (bad or well) informed people behaved in a rather recognizable manner, walking faster, overtaking the others and exhibiting a clear direction of motion, cf.\ \cite{faria2010AB}. No follower has overtaken other people. 

The 5 incognito leaders adopted a good but non optimal strategy, being too close to each other, too much concentrated in the front line, and too fast (cf.\ Fig.\ \ref{fig:S2}(second row)). This was due to the fact that organizers expected a clearer environment in front of the doors and consequently a faster coming out of the people. A better distribution in the crowd would lead to better results, although nothing can be done against a quite noticeable bad informed person. 

\section*{Conclusions and future work}
Both virtual and real experiments suggest that a bottom-up control technique for crowds like the one proposed in the paper is actually feasible and surely deserves further attention. Numerical simulations can give important insights about optimal strategies for leaders, which can be then made ethically acceptable and tested in real situations. The main open issue regards the minimal number of leader needed to reach the desired goal. For that, further theoretical investigation as well as real experiments with a larger number of participants and a clearer environment are needed. The results obtained in Setting 3 should be also experimentally confirmed.

\section*{Acknowledgements} The authors wish to thank M. Fornasier and B. Piccoli for their suggestions and support. They also thank the local organizers of the experiment: A. Tosin, F. S. Priuli S. Cacace, M. Twarogowska; the leaders: M. Di Consiglio, L. Bedarida, M. Arbib, L. Della Torre, M. Sonnino; the head of the Department of Mathematics E. Caglioti and all the students who took part in the experiment. 
Finally, E. Cristiani wishes to dedicate this paper to the memory of the young myrmecologist Dario D'Eustacchio who would have enjoyed seeing people queuing like ants.

\bibliographystyle{siam}
\bibliography{biblio}

\begin{thebibliography}{10}

\bibitem{abdelghany2014EJOR}
{\sc A.~Abdelghany, K.~Abdelghany, H.~Mahmassani, and W.~Alhalabi}, {\em
  Modeling framework for optimal evacuation of large-scale crowded pedestrian
  facilities}, European J. Oper. Res., 237 (2014), pp.~1105--1118.

\bibitem{agnelli2015M3AS}
{\sc J.~P. Agnelli, F.~Colasuonno, and D.~Knopoff}, {\em A kinetic theory
  approach to the dynamics of crowd evacuation from bounded domains}, Math.
  Models Methods Appl. Sci., 25 (2015), pp.~109--129.

\bibitem{albi1401.7798}
{\sc G.~Albi, M.~Herty, and L.~Pareschi}, {\em Kinetic description of optimal
  control problems and applications to opinion consensus}, Commun. Math. Sci.,
  13 (2015), pp.~1407--1429.

\bibitem{albi2013MMS}
{\sc G.~Albi and L.~Pareschi}, {\em Binary interaction algorithms for the
  simulation of flocking and swarming dynamics}, Multiscale Model. Simul., 11
  (2013), pp.~1--29.

\bibitem{albi2013AML}
\leavevmode\vrule height 2pt depth -1.6pt width 23pt, {\em Modeling of
  self-organized systems interacting with a few individuals: from microscopic
  to macroscopic dynamics}, Appl. Math. Lett., 26 (2013), pp.~397--401.

\bibitem{albi1405.0736}
{\sc G.~Albi, L.~Pareschi, and M.~Zanella}, {\em Boltzmann-type control of
  opinion consensus through leaders}, Phil. Trans. R. Soc. A, 372 (2014),
  pp.~20140138/1--18.

\bibitem{andersenhard}
{\sc H.~C. Andersen, J.~D. Weeks, and D.~Chandler}, {\em Relationship between
  the hard-sphere fluid and fluids with realistic repulsion force}, Phys. Rev.
  A, 4 (1971), pp.~1597--1607.

\bibitem{audet2014MPC}
{\sc C.~Audet, K.-C. Dang, and D.~Orban}, {\em Optimization of algorithms with
  {OPAL}}, Math. Prog. Comp., 6 (2014), pp.~233--254.

\bibitem{bellomo2011SR}
{\sc N.~Bellomo and C.~Dogb\'{e}}, {\em On the modeling of traffic and crowds:
  {A} survey of models, speculations, and perspectives}, SIAM Rev., 53 (2011),
  pp.~409--463.

\bibitem{bolley11MMMAS}
{\sc F.~Bolley, J.~A. Ca{\~n}izo, and J.~A. Carrillo}, {\em Stochastic
  mean-field limit: non-{L}ipschitz forces and swarming}, Math. Models Methods
  Appl. Sci., 21 (2011), pp.~2179--2210.

\bibitem{bofo13}
{\sc M.~Bongini and M.~Fornasier}, {\em Sparse stabilization of dynamical
  systems driven by attraction and avoidance forces}, Netw. Heterog. Media, 9
  (2014), pp.~1--31.

\bibitem{bonginijunge2014sparse}
{\sc M.~Bongini, M.~Fornasier, O.~Junge, and B.~Scharf}, {\em Sparse control of
  alignment models in high dimension}, Netw. Heterog. Media, 10 (2015),
  pp.~647--697.

\bibitem{bongini2015conditional}
{\sc M.~Bongini, M.~Fornasier, and D.~Kalise}, {\em {(Un)conditional} consensus
  emergence under perturbed and decentralized feedback controls}, Discrete
  Contin. Dyn. Syst. Ser. A, 35 (2015), pp.~4071--4094.

\bibitem{FornasierSolombRossiBongini}
{\sc M.~Bongini, M.~Fornasier, F.~Rossi, and F.~Solombrino}, {\em Mean-field
  {P}ontryagin maximum principle}.
\newblock Preprint arXiv:1504.02236, 2015.

\bibitem{borzi2015M3ASa}
{\sc A.~Borz\`i and S.~Wongkaew}, {\em Modeling and control through leadership
  of a refined flocking system}, Math. Models Methods Appl. Sci., 25 (2015),
  pp.~255--282.

\bibitem{butail2013PLOSONE}
{\sc S.~Butail, T.~Bartolini, and M.~Porfiri}, {\em Collective response of
  zebrafish shoals to a free-swimming robotic fish}, PLoS ONE, 8 (2013),
  p.~e76123.

\bibitem{CCR11}
{\sc J.~A. Ca{\~n}izo, J.~A. Carrillo, and J.~Rosado}, {\em A well-posedness
  theory in measures for some kinetic models of collective motion}, Math.
  Models Methods Appl. Sci., 21 (2011), pp.~515--539.

\bibitem{caponigro2013MCRF}
{\sc M.~Caponigro, M.~Fornasier, B.~Piccoli, and E.~Tr\'elat}, {\em Sparse
  stabilization and optimal control of the {C}ucker-{S}male model}, Math.
  Control Relat. Fields, 3 (2013), pp.~447--466.

\bibitem{carrillo2009double}
{\sc J.~A. Carrillo, M.~R. D'Orsogna, and V.~Panferov}, {\em Double milling in
  self-propelled swarms from kinetic theory}, Kinet. Relat. Models, 2 (2009),
  pp.~363--378.

\bibitem{carrillo2010particle}
{\sc J.~A. Carrillo, M.~Fornasier, G.~Toscani, and F.~Vecil}, {\em Particle,
  kinetic, and hydrodynamic models of swarming}, in Mathematical modeling of
  collective behavior in socio-economic and life sciences, Springer, 2010,
  pp.~297--336.

\bibitem{carrillo1501.07054}
{\sc J.~A. Carrillo, S.~Martin, and M.-T. Wolfram}, {\em An improved version of
  the {H}ughes model for pedestrian flow}, Math. Models Methods Appl. Sci., 26
  (2016), pp.~671--697.

\bibitem{cercignani1994}
{\sc C.~Cercignani, R.~Illner, and M.~Pulvirenti}, {\em The mathematical theory
  of dilute gases}, Springer, 1994.

\bibitem{cirillo2013PhysA}
{\sc E.~N.~M. Cirillo and A.~Muntean}, {\em Dynamics of pedestrians in regions
  with no visibility - {A} lattice model without exclusion}, Physica A, 392
  (2013), pp.~3578--3588.

\bibitem{colombo2012JNS}
{\sc R.~Colombo and M.~L\'ecureux-Mercier}, {\em An analytical framework to
  describe the interactions between individuals and a continuum}, J. Nonlinear
  Sci., 22 (2012), pp.~39--61.

\bibitem{cordier2005kinetic}
{\sc S.~Cordier, L.~Pareschi, and G.~Toscani}, {\em On a kinetic model for a
  simple market economy}, Journal of Statistical Physics, 120 (2005),
  pp.~253--277.

\bibitem{couzin2005N}
{\sc I.~D. Couzin, J.~Krause, N.~R. Franks, and S.~A. Levin}, {\em Effective
  leadership and decision-making in animal groups on the move}, Nature, 433
  (2005), pp.~513--516.

\bibitem{cristiani_obstacles}
{\sc E.~Cristiani and D.~Peri}, {\em Handling obstacles in pedestrian
  simulations: models and optimization}.
\newblock Preprint arXiv:1512.08528.

\bibitem{cristiani2011MMS}
{\sc E.~Cristiani, B.~Piccoli, and A.~Tosin}, {\em Multiscale modeling of
  granular flows with application to crowd dynamics}, Multiscale Model. Simul.,
  9 (2011), pp.~155--182.

\bibitem{cristiani2014book}
\leavevmode\vrule height 2pt depth -1.6pt width 23pt, {\em Multiscale Modeling
  of Pedestrian Dynamics}, Modeling, Simulation \& Applications, Springer,
  2014.

\bibitem{cristiani2015SIAP}
{\sc E.~Cristiani, F.~S. Priuli, and A.~Tosin}, {\em Modeling rationality to
  control self-organization of crowds: an environmental approach}, SIAM J.
  Appl. Math., 75 (2015), pp.~605--629.

\bibitem{duan2014SR}
{\sc H.~Duan and C.~Sun}, {\em Swarm intelligence inspired shills and the
  evolution of cooperation}, Sci. Rep., 4 (2014), p.~5210.

\bibitem{during2009PRSA}
{\sc B.~D\"uring, P.~Markowich, J.-F. Pietschmann, and M.-T. Wolfram}, {\em
  {B}oltzmann and {F}okker-{P}lanck equations modelling opinion formation in
  the presence of strong leaders}, Proc. R. Soc. A, 465 (2009), pp.~3687--3708.

\bibitem{dyer2009PTRSB}
{\sc J.~R.~G. Dyer, A.~Johansson, D.~Helbing, I.~D. Couzin, and J.~Krause},
  {\em Leadership, consensus decision making and collective behaviour in
  humans}, Phil. Trans. R. Soc. B, 364 (2009), pp.~781--789.

\bibitem{faria2010AB}
{\sc J.~J. Faria, J.~R.~G. Dyer, C.~R. Tosh, and J.~Krause}, {\em Leadership
  and social information use in human crowds}, Animal Behaviour, 79 (2010),
  pp.~895--901.

\bibitem{fornasier2014PTRSA}
{\sc M.~Fornasier, B.~Piccoli, and F.~Rossi}, {\em Mean-field sparse optimal
  control}, Phil. Trans. R. Soc. A, 372 (2014), p.~20130400.

\bibitem{fornasier1306.5913}
{\sc M.~Fornasier and F.~Solombrino}, {\em Mean-field optimal control}, ESAIM
  Control Optim. Calc. Var., 20 (2014), pp.~1123--1152.

\bibitem{guo2012TRB}
{\sc R.-Y. Guo, H.-J Huang, and S.~C. Wong}, {\em Route choice in pedestrian
  evacuation under conditions of good and zero visibility: experimental and
  simulation results}, Transportation Res. B, 46 (2012), pp.~669--686.

\bibitem{Tadmor2008}
{\sc S.~Y. Ha and E.~Tadmor}, {\em From particle to kinetic and hydrodynamic
  descriptions of flocking}, Kinet. Relat. Models, 1 (2008), pp.~415--435.

\bibitem{halloy2007S}
{\sc J.~Halloy, G.~Sempo, G.~Caprari, C.~Rivault, M.~Asadpour, F.~T\^ache,
  I.~Sa\"id, V.~Durier, S.~Canonge, J.~M. Am\'e, C.~Detrain, N.~Correll,
  A.~Martinoli, F.~Mondada, R.~Siegwart, and J.~L. Deneubourg}, {\em Social
  integration of robots into groups of cockroaches to control self-organized
  choices}, Science, 318 (2007), pp.~1155--1158.

\bibitem{han2006JSSC}
{\sc J.~Han, M.~Li, and L.~Guo}, {\em Soft control on collective behavior of a
  group of autonomous agents by a shill agent}, Jrl. Syst. Sci. \& Complexity,
  19 (2006), pp.~54--62.

\bibitem{han2013PLOSONE}
{\sc J.~Han and L.~Wang}, {\em Nondestructive intervention to multi-agent
  systems through an intelligent agent}, PLoS ONE, 8 (2013), p.~e61542.

\bibitem{helbing2001RMP}
{\sc D.~Helbing}, {\em Traffic and related self-driven many-particle systems},
  Rev. Mod. Phys., 73 (2001), pp.~1067--1141.

\bibitem{johansson2007}
{\sc A.~Johansson and D.~Helbing}, {\em Pedestrian flow optimization with a
  genetic algorithm based on boolean grids}, in Pedestrian and Evacuation
  Dynamics 2005, N.~Waldau, P.~Gattermann, H.~Knoflacher, and M.~Schreckenberg,
  eds., Springer-Verlag Berlin Heidelberg, 2007, pp.~267--272.

\bibitem{kachroo2008book}
{\sc P.~Kachroo, S.~J. Al-nasur, S.~A. Wadoo, and A.~Shende}, {\em Pedestrian
  dynamics. Feedback control of crowd evacuation}, Understanding Complex
  Systems, Springer-Verlag, Berlin Heidelberg, 2008.

\bibitem{karper2015hydrodynamic}
{\sc T.~K. Karper, A.~Mellet, and K.~Trivisa}, {\em Hydrodynamic limit of the
  kinetic {C}ucker--{S}male flocking model}, Math. Models Methods Appl. Sci.,
  25 (2015), pp.~131--163.

\bibitem{MayneRawlingsRao2000aa}
{\sc D.~Q. Mayne, J.~B. Rawlings, C.~V. Rao, and P.~O.~M. Scokaert}, {\em
  Constrained model predictive control: stability and optimality}, Automatica
  J. IFAC, 36 (2000), pp.~789--814.

\bibitem{motsch2011JSP}
{\sc S.~Motsch and E.~Tadmor}, {\em A new model for self-organized dynamics and
  its flocking behavior}, J. Stat. Phys., 144 (2011), pp.~923--947.

\bibitem{PT:13}
{\sc L.~Pareschi and G.~Toscani}, {\em Interacting multi-agent systems.
  {K}inetic equations \& {M}onte {C}arlo methods}, Oxford University Press,
  USA, 2013.

\bibitem{parisi2005PhysA}
{\sc D.~R. Parisi and C.~O. Dorso}, {\em Microscopic dynamics of pedestrian
  evacuation}, Physica A, 354 (2005), pp.~606--618.

\bibitem{toscani2006kinetic}
{\sc G.~Toscani}, {\em Kinetic models of opinion formation}, Comm. Math. Sci.,
  4 (2006), pp.~481--496.

\bibitem{toscanibellomoenskog}
{\sc G.~Toscani and N.~Bellomo}, {\em The {E}nskog-{B}oltzmann equation in the
  whole space {$\mathbb R^3$}: some global existence, uniqueness and stability
  results}, Comput. Math. Applic., 13 (1987), pp.~851--859.

\bibitem{twarogowska2014AMM}
{\sc M.~Twarogowska, P.~Goatin, and R.~Duvigneau}, {\em Macroscopic modeling
  and simulations of room evacuation}, Appl. Math. Model., 38 (2014),
  pp.~5781--5795.

\bibitem{vicsek1995novel}
{\sc T.~Vicsek, A.~Czir{\'o}k, E.~Ben-Jacob, I.~Cohen, and O.~Shochet}, {\em
  Novel type of phase transition in a system of self-driven particles}, Phys.
  Rev. Lett., 75 (1995), pp.~1226--1229.

\bibitem{vicsek2012PR}
{\sc T.~Vicsek and A.~Zafeiris}, {\em Collective motion}, Phys. Rep., 517
  (2012), pp.~71--140.

\bibitem{Vil02}
{\sc C.~Villani}, {\em Handbook of Mathematical Fluid Dynamics}, vol.~1,
  Elsevier, 2002, ch.~A review of mathematical topics in collisional kinetic
  theory.

\bibitem{wang2015PhysA}
{\sc J.~Wang, L.~Zhang, Q.~Shi, P.~Yang, and X.~Hu}, {\em Modeling and
  simulating for congestion pedestrian evacuation with panic}, Physica A, 428
  (2015), pp.~396--409.

\end{thebibliography}
\end{document}